\documentclass[12pt]{amsart}
\font\emailfont=cmtt10

\usepackage{amsmath,amsthm,amsfonts,amscd,flafter,epsf}

\newcommand\commentable[1]{#1}
\newcommand\Rk{\mathrm{rk}}
\newcommand\Id{\mathrm{Id}}

\newcommand{\Tors}{\mathrm{Tors}}
\newcommand{\rk}{\mathrm{rk}}
\newcommand{\HF}{HF}

\newtheorem{theorem}{Theorem}[section]
\newtheorem{prop}[theorem]{Proposition}
\newtheorem{cor}[theorem]{Corollary}

\newtheorem{lemma}[theorem]{Lemma}

\newtheorem{defn}[theorem]{Definition}

\newtheorem{remark}[theorem]{Remark}

\def\endproof{\relax\ifmmode\expandafter\endproofmath\else
  \unskip\nobreak\hfil\penalty50\hskip.75em\hbox{}\nobreak\hfil\bull
  {\parfillskip=0pt \finalhyphendemerits=0 \bigbreak}\fi}
\def\endproofmath$${\eqno\bull$$\bigbreak}
\def\bull{\vbox{\hrule\hbox{\vrule\kern3pt\vbox{\kern6pt}\kern3pt\vrule}\hrule}}

\newcommand{\Q}{\mathbb{Q}}
\newcommand{\R}{\mathbb{R}}

\newcommand{\C}{\mathbb{C}}

\newcommand{\Z}{\mathbb{Z}}

\newcommand{\OneHalf}{\frac{1}{2}}

\newcommand{\CP}[1]{{\mathbb{CP}}^{#1}}

\newcommand{\Zmod}[1]{\Z/{#1}\Z}

\newcommand{\Ker}{\mathrm{Ker}}
\newcommand{\CoKer}{\mathrm{Coker}}

\newcommand{\Image}{\mathrm{Im}}

\newcommand{\cm}{\cdot}

\newcommand{\ModSWfour}{\mathcal{M}}
\newcommand{\ModFlow}{\ModSWfour}

\newcommand{\SpinC}{{\mathrm{Spin}}^c}
\newcommand{\Spin}{{\mathrm{Spin}}}

\newcommand\sgn{\mathrm{sgn}}

\newcommand\Wedge{\Lambda}

\newcommand\abuts\Rightarrow
\newcommand\Sym{\mathrm{Sym}}

\newcommand\Laurent{\mathbb L}
\newcommand\renEuler{\widehat\chi}
\newcommand\spinccan{\mathfrak k}

\newcommand\HFpRed{\HFp_{\red}}

\newcommand\liftGr{\widetilde\gr}
\newcommand\chiTrunc{\chi^{\mathrm{trunc}}}
\newcommand\chiRen{\widehat\chi}
\newcommand\Fred[1]{F^{\mathrm{red}}_{#1}}

\newcommand\HFinfty{\HFinf}

\newcommand\x{\mathbf x}
\newcommand\w{\mathbf w}

\newcommand\y{\mathbf y}

\newcommand\ModSphere{\ModFlow\left({\mathbb S}\longrightarrow 
\Sym^{g-1}(\Sigma_{1})\times \Sym^2(\Sigma_{2})\right)}
\newcommand\ModSpheres\ModSphere
\newcommand\CF{CF}

\newcommand\CFa{\widehat{CF}}
\newcommand\CFp{\CFb}
\newcommand\CFm{\CF^-}

\newcommand\HFleq{\HF^{\leq 0}}

\newcommand\HFred{\HF_{\rm red}}

\newcommand{\red}{\mathrm{red}}

\newcommand\HFp{\HFb}
\newcommand\HFpm{HF^{\pm}}
\newcommand\HFm{\HF^-}
\newcommand\CFinf{CF^\infty}
\newcommand\HFinf{HF^\infty}
\newcommand\CFb{CF^+}
\newcommand\HFa{\widehat{HF}}
\newcommand\HFb{HF^+}
\newcommand\gr{\mathrm{gr}}
\newcommand\Mas{\mu}
\newcommand\UnparModSp{\widehat \ModSp}
\newcommand\UnparModFlow\UnparModSp
\newcommand\Mod\ModSp

\newcommand{\cald}{{\mathcal D}}

\newcommand\PD{\mathrm{PD}}

\newcommand{\spinc}{\mathfrak s}

\newcommand{\spinct}{\mathfrak t}

\newcommand\ModMaps{\mathcal M}
\newcommand\ModSp\ModMaps

\newcommand\Ta{{\mathbb T}_{\alpha}}
\newcommand\Tb{{\mathbb T}_{\beta}}
\newcommand\Tc{{\mathbb T}_{\gamma}}
\newcommand\Td{{\mathbb T}_{\delta}}

\newcommand\alphas{\mbox{\boldmath$\alpha$}}

\newcommand\betas{\mbox{\boldmath$\beta$}}
\newcommand\gammas{\mbox{\boldmath$\gamma$}}
\newcommand\deltas{\mbox{\boldmath$\delta$}}

\newcommand\PerDom{\mathcal P}

\newcommand\uCF{\underline\CF}
\newcommand\uCFinf{\uCF^\infty}
\newcommand\uHF{\underline{\HF}}
\newcommand\uHFp{\underline{\HF}^+}

\newcommand\uCFp{\underline{\CFp}}

\newcommand\uHFa{\underline{\HFa}}

\newcommand\uHFinf{\uHF^\infty}

\newcommand\Fm[1]{F^{-}_{#1}}
\newcommand\Fp[1]{F^{+}_{#1}}
\newcommand\Fpm[1]{F^{\pm}_{#1}}
\newcommand\Fa[1]{{\widehat F}_{#1}}
\newcommand\Finf[1]{F^{\infty}_{#1}}
\newcommand\Fleq[1]{F^{\leq 0}_{#1}}

\newcommand\liftalpha{\widetilde\alpha}
\newcommand\liftbeta{\widetilde\beta}
\newcommand\sj{\mathfrak j}
\newcommand\MT{t}

\newcommand\KOver{K_+}
\newcommand\KUnder{K_-}

\newcommand\prefspinc{\mathfrak u}
\newcommand\spincu{\mathfrak u}
\newcommand\SpinCCobord{\theta^c}
\newcommand\Field{\mathbb F}

\newcommand\Fc{F^\circ}

\newcommand\HFc{\HF^\circ}

\newcommand\Dual{\mathcal D}
\newcommand\Duality\Dual

\newcommand\liftTa{{\widetilde{\mathbb T}}_\alpha}
\newcommand\liftTb{{\widetilde{\mathbb T}}_\beta}
\newcommand\liftTc{{\widetilde{\mathbb T}}_\gamma}

\newcommand\liftTcPr{{\widetilde{\mathbb T}}_\gamma'}

\newcommand\Tor{\mathrm{Tor}}

\newcommand\TcPr{\Tc'}

\newcommand\Spider{\sigma}
\newcommand\EulerMeasure{\widehat\chi}

\newcommand\Knot{\mathbb K}
\newcommand\Link{\mathbb L}

\newcommand\liftSym{{\widetilde \Sym}^g(\Sigma)_\xi}
\newcommand\liftCFinf{\widetilde{CF}^\infty}
\newcommand\liftz{\widetilde z}
\newcommand\liftCFp{\widetilde{CF}^+}

\newcommand\liftHFp{\widetilde{HF}^+}
\newcommand\liftx{\widetilde\x}
\newcommand\lifty{\widetilde\y}
\newcommand\liftSigma{\widetilde\Sigma}
\newcommand\liftalphas{\widetilde \alphas}
\newcommand\liftbetas{\widetilde \betas}
\newcommand\liftgammas{{\widetilde\gammas}}
\newcommand\liftdeltas{{\widetilde\deltas}}
\newcommand\liftgammasPr{\liftgammas'}
\newcommand\liftGamma{\widetilde{\Gamma}}
\newcommand\liftdeltasPr{\liftdeltas'}
\newcommand\Tz{{\widetilde{\mathbb T}}_\zeta}
\newcommand\liftTheta{\widetilde\Theta}

\newcommand\liftF[1]{{\widetilde F}_{#1}}
\newcommand\spinck{\mathfrak k}

\title[{Absolutely Graded Floer Homologies}]
{Absolutely graded Floer homologies and intersection forms 
for four-manifolds with boundary}
\author[Peter Ozsv{\'a}th]{Peter Ozsv\'ath}
\address{Department of
Mathematics,  Princeton University, New Jersey 08540 \newline
\indent{\emailfont{petero@math.princeton.edu}}}
\thanks{PSO was supported by NSF grant number DMS 9971950 and a Sloan 
Research Fellowship}

\author[Zolt{\'a}n Szab{\'o}]{Zolt{\'a}n Szab{\'o}} 
\address{Department of
Mathematics,  Princeton University, New Jersey 08540 \newline
\indent{\emailfont{szabo@math.princeton.edu}}}
\thanks{ZSz was supported by NSF grant number DMS 0107792
and a Packard Fellowship}

\begin{document}

\begin{abstract}  
In~\cite{HolDiskFour}, we introduced absolute gradings
on the three-manifold invariants developed in~\cite{HolDisk} and
\cite{HolDiskTwo}. Coupled with the surgery long exact sequences,
we obtain a number of three- and four-dimensional applications of this
absolute grading including strengthenings of the ``complexity bounds''
derived in~\cite{HolDiskTwo}, restrictions on knots whose surgeries
give rise to lens spaces, and calculations of $\HFp$ for a variety of
three-manifolds. Moreover, we show how the structure of $\HFp$
constrains the exoticness of definite intersection forms for smooth
four-manifolds which bound a given three-manifold.  In addition to
these new applications, the techniques also provide alternate proofs
of Donaldson's diagonalizability theorem and the Thom conjecture for
$\CP{2}$.
\end{abstract} 

\maketitle
\section{Introduction}

In~\cite{HolDisk} and \cite{HolDiskTwo}, we introduced Floer-homology
invariants ($\HFa$, $\HFpm$, and $\HFinf$) for oriented
three-manifolds equipped with $\SpinC$
structures. In~\cite{HolDiskTwo}, we derived a number of exact
sequences relating $\HFp$ (and also $\HFa$) for three-manifolds which differ by
surgeries.  In~\cite{HolDiskFour}, we introduced the corresponding
four-dimensional theory, where cobordisms $W$ between three-manifolds
induce maps $\Fa{}$, $\Fpm{}$, $\Finf{}$ between the corresponding
Floer homology theories on the boundaries. One by-product of this
four-dimensional theory is an absolute rational lift to the relative
$\Z$ grading on the homology groups for a three-manifold $Y$ endowed
with a torsion $\SpinC$ structure (i.e. one whose first Chern class is
a torsion element of $H^2(Y;\Z)$).

The present article presents a number of three- and four-dimensional
applications of the interplay
between the absolute gradings and the surgery exact sequences.

\subsection{The correction term $d(Y)$}

After two introductory sections (Section~\ref{sec:Review} where we
review some of the terminology from the earlier articles, and
Section~\ref{sec:GradExact}, where we pin down the relationship
between the absolute gradings and the exact sequence), we give the
first application of the absolute gradings: a ``correction term'' for
three-manifolds. For simplicity, we presently assume that $Y$ is a
rational homology three-sphere, equipped with a $\SpinC$ structure
$\spinct$. With the help of the absolute grading, we define a
numerical invariant $d(Y,\spinct)$ for $Y$, which is the minimal
degree of any non-torsion class in $\HFp(Y,\spinct)$ coming from
$\HFinf(Y,\spinct)$. This invariant is analogous to a gauge-theoretic
invariant of Fr{\o}yshov (see~\cite{Froyshov}).

The basic properties of $d$ are conveniently stated with the help of a
result which we state after introducing some terminology.

\begin{defn}
The three-dimensional $\SpinC$ homology bordism group $\SpinCCobord$ is
the set of equivalence classes of 
pairs $(Y,\spinct)$ where $Y$ is a rational homology three-sphere,
and $\spinct$ is a $\SpinC$ structure over $Y$, and the equivalence
relation identifies $(Y_1,\spinct_1)\sim (Y_2,\spinct_2)$ if there is a
(connected, oriented, smooth) cobordism $W$ from $Y_1$ to $Y_2$ with
$H_i(W,\Q)=0$ for $i=1$ and $2$, which can be endowed with a $\SpinC$
structure $\spinc$ whose restrictions to $Y_1$ and $Y_2$ are $\spinct_1$
and $\spinct_2$ respectively. The connected sum operation endows this
set with the structure of an Abelian group (whose unit is $S^3$
endowed with its unique $\SpinC$ structure). 
\end{defn}

There is a classical homomorphism $$\rho\colon
\theta^c\longrightarrow \Q/\Z$$ (see for instance~\cite{APSII}),  defined
as follows. Consider a rational homology three-sphere $(Y,\spinct)$,
and let $X$ be any four-manifold equipped with a $\SpinC$ structure
$\spinc$ with $\partial X \cong Y$ and $\spinc|\partial X \cong
\spinct$. Then $$\rho(Y,\spinct)\equiv \frac{c_1(\spinc)^2-\sgn(X)}{4}
\pmod{2\Z}$$ where $\sgn(X)$ denotes the signature of the intersection
form of $X$.

\begin{theorem}
\label{intro:HomologyBordism}
The numerical invariant $d(Y,\spinct)$ descends to give a group homomorphism
$$d\colon \SpinCCobord\longrightarrow \Q$$
which is a lift of $\rho$. Moreover, 
$d$ is  invariant under conjugation; i.e.
$d(Y,\spinct)=d(Y,{\overline \spinct})$.
\end{theorem}

Additivity of the correction term under connected sums, and its
behaviour under orientation reversal, is established
Section~\ref{sec:CorrTerm}.  Its rational homology cobordism
invariance is established only after more of the four-dimensional theory is
developed in Section~\ref{sec:DefiniteForms}.  When $b_1(Y)>0$,
torsion $\SpinC$ structures can be endowed with a collection of
correction terms. In Subsection~\ref{subsec:CorrTermBig}, we discuss
this construction in the case where $H_1(Y;\Z)\cong \Z$.

\subsection{Regularized Euler characteristics}

When $Y$ is an integral homology three-sphere,
$\HFp(Y)$ is infinitely
generated. (We drop the $\SpinC$ structure from the notation, since
it is unique.) To obtain a finitely generated group, we must pass to the
quotient $$\HFpRed(Y)=\HFp(Y)/\HFinf(Y).$$
We
have the following relationship between the Euler characteristic of
$\HFpRed(Y)$, the correction term for $Y$, and Casson's invariant
$\lambda(Y)$,

\begin{theorem}
\label{intro:EulerCasson}
Let $Y$ be an integer homology three-sphere. Then, 
the following relationship holds between
Casson's invariant, the Euler characteristic of $\HFp$ (thought
of as a $\Z$-module), 
and the correction term:
$$
\lambda(Y)=\chi(\HFpRed(Y))-\OneHalf d(Y),
$$ where here Casson's invariant is normalized so that, if
$\Sigma(2,3,5)$ denotes the Poincar\'e homology sphere, oriented as
the boundary of the negative-definite $E8$ plumbing, then
$\lambda(\Sigma(2,3,5))=-1$.
\end{theorem}

\subsection{\bf{Invisible three-manifolds.}}
When $Y$ is an integer homology three-sphere, we define its
{\em complexity} $N(Y)$ to be the rank of $\HFred(Y)$. Let $K\subset Y$ be a
knot in $Y$, $r$ be a rational number, and let $Y_r(K)$ be the integer
homology three-sphere obtained by $r$ surgery on $Y$ along $K$.  In
Theorem~\ref{HolDiskTwo:thm:Complexity} of~\cite{HolDiskTwo}, we gave
lower bounds on the sum of complexities of $Y$ and $Y_{1/n}(K)$ (where
$n$ is a whole number), in terms of the Alexander polynomial of
$K$. Using $d(Y)$ in conjunction with $N(Y)$, these bounds can be
strengthened, according to the following:

\begin{theorem}
\label{intro:RenormalizedComplexityFrac} 
Let $Y$ be an integral homology three-sphere and $K\subset Y$ be a
knot, whose symmetrized Alexander polynomial is $$\Delta_K(T)= a_0(K)+ \sum_{i=1}^d
a_i(K)\cm (T^{i}+T^{-i}),$$ then for each positive integer $n$,
there is a bound $$ n\left(|\MT_0(K)|+
2\sum_{i=1}^d |\MT_i(K)|\right) \leq N(Y)+\frac{d(Y)}{2} +
N(Y_{1/n}(K))-\frac{d(Y_{1/n}(K))}{2},$$
where
\begin{equation}
\label{eq:defBi}
\MT_i(K)=  \sum_{j=1}^d j a_{|i|+j}(K).
\end{equation}
\end{theorem}

An integral homology three-sphere $Y$ for which both $N(Y)=0=d(Y)$ is
called {\em invisible}, since these invariants do not distinguish it
from $S^3$.  Of course (in view of Theorem~\ref{intro:EulerCasson}),
all invisible homology three-spheres have $\lambda(Y)=0$.  We prove in
Section~\ref{sec:Complexity} that the set of invisible three-manifolds
$Y$ is closed under connected sum and orientation reversal.

Theorem~\ref{intro:RenormalizedComplexityFrac} has the following immediate consequence:

\begin{cor}
If $K\subset Y$ is a knot with non-trivial Alexander polynomial in an
invisible three-manifold, then any non-trivial surgery on $K$ gives a
three-manifold which is not invisible.
\end{cor}

\subsection{\bf{Surgeries giving lens spaces.}}

The graded Floer homologies can also be used to place restrictions on knots
$K\subset S^3$ whose surgeries give lens spaces. For example, we have the following:

\begin{theorem}
\label{intro:SmallLenses}
If $K\subset S^3$ is a knot with the property that some integral surgery along
$K$ with integral coefficient $p$ with $|p|\leq 4$ gives a lens space,
then $\HFp(Y_0(K))\cong \HFp(S^2\times S^1)$. In particular, the Alexander polynomial
of $K$ is trivial.
\end{theorem}

\begin{remark}
The ``cyclic surgery theorem'' of
Culler-Gordon-Luecke-Shalen~\cite{CyclicSurgery} ensures that if $K$
is not a torus knot, and some surgery along $K$ gives a lens space,
then that surgery is integral. It is conjectured in~\cite{Gordon} that
the trivial knot is the only knot $K$ with the property that some
surgery on $K$ gives a lens space $L(p,1)$ with $|p|\leq 4$.
\end{remark}

A $+5$-surgery on the right-handed trefoil gives the lens space
$-L(5,1)$, so care must be taken in generalizing
Theorem~\ref{intro:SmallLenses}. The following result holds for
general $p$:

\begin{theorem}
\label{intro:PLens}
Suppose that $K\subset S^3$ is a knot with the property that $+p$
surgery on $K$ gives the lens space $L(p,1)$ (with its canonical
orientation). Then, $\HFinf(S^3_0(K))\cong\HFinf(S^2\times S^1)$ as
absolutely graded groups; in particular, the Alexander polynomial of
$K$ is trivial.
\end{theorem}

In fact, Gordon has conjectured (see~\cite{GordonConjecture}) that 
if $+p$ surgery on a knot $K\subset S^3$ gives $L(p,1)$, then $K$ is
the unknot. Moreover,
Berge~\cite{Berge} has conjectured a complete list of knots which give
rise to lens spaces.

Theorems~\ref{intro:SmallLenses} and \ref{intro:PLens} are proved in
Section~\ref{sec:Lens}, where they are given as a corollaries to a
theorem, Theorem~\ref{thm:PSurgeryLens}, which, for any fixed $p$ and
$q$, gives strong restrictions on knots $K\subset S^3$ for which $+p$
surgery gives the lens space $L(p,q)$. Specifically, if $K$ is such a
knot, the theorem shows how $\HFp$ of the zero-surgery (and hence the
Alexander polynomial of the knot) is determined, up to a finite
indeterminacy, by the correction terms for $L(p,q)$. The indeterminacy
comes from certain allowed permutations of these correction terms. The
terms themselves are calculated in Proposition~\ref{prop:dLens}
(see~\cite{HolSomePlumb} for a more conceptual interpretation of these
numbers). For another illustration of these constraints, we tabulate
at the end of this paper (c.f Subsection~\ref{subsec:LensAlex}) the
possible Alexander polynomials of knots in $S^3$ for which some
positive integral surgery $p\leq 26$ gives a lens space.

Turned around, Theorem~\ref{thm:PSurgeryLens} also gives obstructions
to realizing lens spaces as integral surgeries on a knot in the
three-sphere. For
example, these techniques show that the lens space $L(22,3)$ -- and
indeed, any lens spaces in the family $L(2k(3+8k),2k+1)$ where $k$ is
any positive integer not divisible by $4$ -- is not obtained by
integral surgery on a single knot $K\subset S^3$. This obstruction is
``$S^3$-specific''. An argument of Fintushel and Stern gives a
necessary and sufficient condition for $L(p,q)$ to be realized as
integral surgery on a knot in a {\em homology} sphere: the condition
is that $\pm q$ be a square modulo $p$, a condition which all the lens
spaces in this family satisfy.
(This application is described in Proposition~\ref{prop:FamiliesOfLenses}.)

\subsection{\bf{Some calculations.}}
We give several calculations of $N(Y)$ and $d(Y)$, including
calculations of both invariants for integer homology spheres obtained
as surgeries on the $(p,q)$ torus knot. Indeed, we calculate all of
$\HFp$ for the Brieskorn spheres $\Sigma(2,3,6n\pm1)$. Moreover, we
calculate $\HFp(E_n)$ where $E_n$ is the three-manifold obtained by
$1/n$ surgery with $n>0$ on the figure eight knot in $S^3$
(c.f. Proposition~\ref{prop:FigureEight}): $$\HFp_k(E_n)\cong
\left\{\begin{array}{ll}
\Z &   {\text{if $k\equiv 0\pmod{2}$ and $k\geq 0$}} \\
\Z^{n} & {\text{if $k=-1$}} \\
0 & {\text{otherwise}}
\end{array}
\right..$$ 

We also give the following simple calculation for the invariants of $T^3$:

\begin{prop}
\label{intro:T3}
Let $T^3$ denote the three-dimensional torus. Then, we have $H_1(T^3;\Z)$-module isomorphisms:
\begin{eqnarray*}
\HFa(T^3)&\cong &H^2(T^3;\Z)\oplus H^1(T^3;\Z), \\
\HFp(T^3)&\cong & \Big(H^2(T^3;\Z)\oplus H^1(T^3;\Z)\Big)\otimes_\Z \Z[U^{-1}], \\
\HFinf(T^3)&\cong & \Big(H^2(T^3;\Z)\oplus H^1(T^3;\Z)\Big)\otimes_\Z \Z[U,U^{-1}].
\end{eqnarray*}
The absolute grading is symmetric, in the sense that
$\liftGr(H^2(T^3;\Z)\subset \HFa(T^3))=1/2$, 
$\liftGr(H^1(T^3;\Z)\subset \HFa(T^3))=-1/2$.
\end{prop}

Further calculations are given for surgeries on certain pretzel knots,
via surgery long exact sequences which mirror the skein relations
for calculating Alexander polynomials.

\subsection{\bf{Intersection forms of definite four-manifolds}.}

By analyzing the maps $F^\infty_{W,\spinc}$ for negative-definite
cobordisms, we give an alternative proof of Donaldson's
diagonalizability theorem:

\begin{theorem} (Donaldson)
\label{intro:Donaldson}
If $X$ is a smooth, closed, oriented four-manifold with definite
intersection form, then the form of $X$ is diagonalizable over $\Z$.
\end{theorem}

Indeed, applying the same idea to four-manifolds which bound a given integer homology
three-sphere $Y$, we obtain the following analogue of a theorem of Fr{\o}yshov.

Let $Y$ be an integral homology three-sphere, and $X$ be a four-manifold which bounds $Y$. There is an 
intersection form
$$Q_X \colon (H_2(X;\Z)/\Tors) \otimes (H_2(X;\Z)/\Tors) \longrightarrow \Z.$$
A characteristic vector for this intersection form is an element $\xi\in H_2(X;\Z)/\Tors$
with
$$\xi\cdot v \equiv v\cm v \pmod{2}$$
for each $v\in H_2(X;\Z)/\Tors$. (We write $\xi\cdot\eta$ for
$Q_X(\xi,\eta)$, and $\xi^2$ for $Q_X(\xi,\xi)$.)

\begin{theorem}
\label{intro:IntForm}
Let $Y$ be an integral homology three-sphere, then 
for each negative-definite four-manifold $X$ which bounds $Y$, we have the inequality:
$$
	\xi^2 + \rk(H^2(X;\Z)) \leq 4d(Y),
$$
for each characteristic vector $\xi$.
\end{theorem}

This theorem has an obvious generalization to rational homology
three-spheres, and, indeed, generalizations for three-manifolds with
$b_1(Y)>0$ discussed in Section~\ref{subsec:BOneBigForms}. When
$H_1(Y;\Z)\cong \Z$, the inequality can be used as obstruction to
realizing $Y$ as the boundary of an integral homology $S^2\times D^2$.
Moreover, when applied to certain circle bundles over two-manifolds,
the intersection form bounds provide another proof of the Thom
conjecture for $\CP{2}$ first proved by
Kronheimer-Mrowka~\cite{KMthom} and
Morgan-Szab{\'o}-Taubes~\cite{MSzT},
c.f. Subsection~\ref{subsec:Thom}. The proof given here is analogous
to a Seiberg-Witten proof given recently by Strle,
see~\cite{Strle}. 

In Section~\ref{sec:Examples}, we close with some other applications
of the intersection form bounds, as combined with the calculations
from Section~\ref{sec:SampleCalculations}. Specifically, we show
how the results constrain the intersection forms of four-manifolds
which bound certain Seifert fibered spaces. In particular, in 
Subsection~\ref{subsec:NotSurgeryOnKnot}, we 
exhibit a three-manifold whose first homology is isomorphic to $\Z$ and 
which can
be expressed as surgery on a certain two-component link in $S^3$, but
which cannot be expressed as surgery on a single knot.  Finally, we 
close with a table of allowed Alexander polynomials of knots giving rise
to (small) lens spaces.

\vskip.2cm
\noindent{\bf Acknowledgements:} It is our pleasure to thank 
Cameron Gordon,  Robion Kirby, Tomasz Mrowka, and Andr\'as Stipsicz for their
encouragement during the preparation of this work. We would especially like
to thank Cameron Gordon for some very useful suggestions on an early
version of this paper.

\section{Floer homologies and maps between them}
\label{sec:Review}

In~\cite{HolDiskFour}, we showed how the constructions
from~\cite{HolDisk} and \cite{HolDiskTwo} can give rise to a naturally
associated assignment which associates to a three-manifold $Y$
equipped with a $\SpinC$ structure $\spinct$, four relatively
$\Zmod{c_1(\spinct)\cup H^1(Y;\Z)}$-graded Abelian groups $\HFa$,
$\HFp$, $\HFm$, and $\HFinf$, the latter three of which are also
endowed with the action by a polynomial algebra $\Z[U]$, where
multiplication by $U$ lowers relative degree by $2$. All four
groups are acted upon by the exterior algebra $\Wedge^*
H_1(Y;\Z)/\Tors$.

In fact, these groups are related by functorially associated long exact sequences:
\begin{equation}
\label{eq:HFaExactSequence}
\begin{CD}
... @>>>\HFa(Y,\spinct)@>{\widehat\iota}>>\HFp(Y,\spinct)@>{U}>> \HFp(Y,\spinct)@>>> ...
\end{CD}
\end{equation}
and 
\begin{equation}
\label{eq:HFinfExactSequence}
\begin{CD}
... @>>>\HFm(Y,\spinct)@>{\iota}>>\HFinf(Y,\spinct)@>{\pi}>> \HFp(Y,\spinct)@>>> ...
\end{CD}
\end{equation}
The second long exact sequence (which we use more often) will
frequently be abbreviated $\HFc(Y,\spinct)$. Note that maps in the
first sequence are equivariant under the action of $\Wedge^*(Y)/\Tors$,
while those for the second are equivariant under the action of
$\Z[U]\otimes_\Z \Wedge^*(Y)/\Tors$.

When $W$ is a (smooth) cobordism from a three-manifold $Y_1$ to $Y_2$,
equipped with a $\SpinC$ structure $\spinc$ whose restrictions to the
two boundary components are $\spinct_1$ and $\spinct_2$ respectively, then
there are induced maps between the long exact sequences: $$
\begin{CD}
... @>>>\HFm(Y_1,\spinct_1)@>{\iota}>>\HFinf(Y_1,\spinct_1)@>{\pi}>> \HFp(Y_1,\spinct_1)@>>> ... \\
&& @V{\Fm{W,\spinc}}VV @V{\Finf{W,\spinc}}VV @V{\Fp{W,\spinc}}VV \\
... @>>>\HFm(Y_2,\spinct_2)@>{\iota}>>\HFinf(Y_2,\spinct_2)@>{\pi}>> \HFp(Y_2,\spinct_2)@>>> ...\\
\end{CD}
$$
and
$$
\begin{CD}
... @>>>\HFa(Y_1,\spinct_1)@>{\widehat\iota}>>\HFp(Y_1,\spinct_1)@>{U}>> \HFp(Y_1,\spinct_1)@>>> ... \\
&& @V{\Fa{W,\spinc}}VV @V{\Fp{W,\spinc}}VV @V{\Fp{W,\spinc}}VV \\
... @>>>\HFa(Y_2,\spinct_2)@>{\widehat\iota}>>\HFp(Y_2,\spinct_2)@>{U}>> \HFp(Y_2,\spinct_2)@>>> ...\\
\end{CD}
$$ The maps $\Fa{W,\spinc}$, $\Fpm{W,\spinc}$, and $\Finf{W,\spinc}$
are uniquely defined up to multiplication by an overall sign $\pm
1$. 

These maps also admit an action by $\Wedge^*H_1(W)/\Tors$.  This
action is related to the analogous actions over the three-manifolds,
according to the following formula. Fix arbitrary elements
$\gamma_i\in H_1(Y_i)/\Tors$ for $i=1,2$, and let $\gamma\in
H_1(W)/\Tors$ be the homology class obtained as
$\gamma=j_1(\gamma_1)-j_2(\gamma_2)$, where for $i=1,2$, $j_i$ denotes the natural inclusion
$j_i\colon H_1(Y_i)\longrightarrow H_1(X)$.  Then, $$\pm (\gamma\otimes
F_{W,\spinc})(\xi)=F_{W,\spinc}(\gamma_1\cm \xi)-\gamma_2 \cm
F_{W,\spinc}(\xi).$$ In particular, if $\gamma_1$ and $\gamma_2$ are homologous in $W$, then
$$F_{W,\spinc}(\gamma_1\cm\xi)=\gamma_2 \cm F_{W,\spinc}(\xi).$$

Suppose now that $Y$ is equipped with a torsion $\SpinC$ structure. In
Theorem~\ref{HolDiskFour:thm:AbsGrade} of~\cite{HolDiskFour}, we showed that
$\HFc(Y,\spinct)$, and also $\HFa(Y,\spinct)$ can be given an absolute $\Q$ grading $\liftGr$
which lifts the relative $\Z$ grading. It is uniquely characterized by
the following properties:
\begin{itemize}
\item ${\widehat \iota}$, $\iota$, and $\pi$ above preserve the absolute grading
\item $\HFa(S^3)$ is supported in absolute grading zero
\item if $W$ is a cobordism
from $Y_1$ to $Y_2$, and $\xi\in \HFinf(Y_1,\spinct_1)$, then
\begin{equation}
\label{eq:DimensionShiftFormula}
\liftGr(F_{W,\spinc}(\xi))-\liftGr(\xi) = 
\frac{c_1(\spinc)^2-2\chi(W)-3\sigma(W)}{4},
\end{equation}
where $\spinct_i=\spinc|Y_i$ for $i=1,2$.
\end{itemize}
The existence and characterization of this
absolute grading is the key result we use here from~\cite{HolDiskFour}.
\section{Exact sequences and absolute gradings}
\label{sec:GradExact}

Recall that in~\cite{HolDiskTwo}, surgery long exact sequences for
$\HFp$ and $\HFa$ were established. We give now a graded refinement,
focusing mainly on the case of $\HFp$. The
key point is the relationship between the maps in the long exact
sequences and the maps induced by the cobordisms obtained by surgeries
on knots.

More concretely, recall that when $K$ is a knot in an integral
homology three-sphere, Theorem~\ref{HolDiskTwo:thm:ExactOne} of~\cite{HolDiskTwo}
gives a long exact sequence $$
\begin{CD}
...@>>>\HFp(Y)@>{F_1}>> \HFp(Y_0) @>{F_2}>>\HFp(Y_1)@>{F_3}>>\HFp(Y)@>>>...,
\end{CD}
$$ where $Y_0=Y_0(K)$ and $Y_1=Y_1(K)$ are the three-manifolds
obtained by $0$-surgery and $+1$-surgery on $Y$ along $K$ (using its
canonical framing). Here,
$$\HFp(Y_0)=\bigoplus_{\spinct\in\SpinC(Y_0)}\HFp(Y_0,\spinct).$$ The
maps $F_1$ and $F_2$ were defined by counting pseudo-holomorphic
triangles in a Heegaard triple. An easy comparison between the
definition of the maps in~\cite{HolDiskTwo} and the maps associated to
two-handles (Subsection~\ref{HolDiskFour:subsec:TwoHandles} of~\cite{HolDisk}) shows that 
$F_1$ and $F_2$ are sums of maps associated to cobordisms; i.e.
\begin{eqnarray*}
F_1= \sum_{\spinc\in W_0}\pm F_{W_0,\spinc},&{\text{and}}&
F_2=\sum_{\spinc\in W_1}\pm F_{W_1,\spinc},
\end{eqnarray*}
where $W_0$ is the cobordism from $Y$ to $Y_0$ defined by attaching a
single two-handle to $Y$ (with $0$-framing), and $W_1$ is the
cobordism from $Y_0$ to $Y_1$ defined by attaching a single two-handle
to $Y_0$. (The signs are chosen as in the proof of
Theorem~\ref{HolDiskTwo:thm:ExactOne} of~\cite{HolDiskTwo} to make the
sequence exact.)  Observe that $\SpinC$ structures on both $W_0$ and
$W_1$ are uniquely determined by their restrictions to $Y_0$.

Note that the map $F_3$ arises as the coboundary map of an associated
long exact sequence, though one could alternatively realize it, too,
as a sum of maps belonging to cobordisms
(compare~\cite{BraamDonaldson}). However, for our present purposes, it
suffices to understand properties of the maps $F_1$ and $F_2$.

\begin{lemma}
\label{lemma:DegShiftExact}
Let $K\subset Y$ be a knot in an integral homology three-sphere, and
let $\spinct_0$ denote the $\SpinC$ structure over $Y_0$ with trivial
first Chern class. In the exact sequence 
$$
\begin{CD}
...@>>>\HFp(Y)@>{F_1}>> \HFp(Y_0) @>{F_2}>>\HFp(Y_1)@>{F_3}>>\HFp(Y)@>>>...,
\end{CD}
$$
the component of $F_1$ mapping into $\HFp(Y_0,\spinct_0)$ (now thought of as absolutely $\Q$-graded) has degree $-1/2$, 
the restriction of $F_2$ to $\HFp(Y_0,\spinct_0)$ has degree $-1/2$.
\end{lemma}

\begin{proof}
In view of the above remarks, the component of $F_1$ landing in the
$\HFp(Y_0,\spinct_0)$ summand is the map induced by the cobordism $W_0$,
endowed with the $\SpinC$ structure $\spinc_0$ with trivial first
Chern class. This cobordism is obtained by a single one-handle
addition (so its $\chi=0$), and its signature $\sigma=0$. Thus, by
Equation~\eqref{eq:DimensionShiftFormula}, the result follows. The map
$F_2$ is also obtained by a single two-handle addition, with signature zero.
\end{proof}

\begin{remark}
Note that the above discussion could have been made using $\HFa$ in
place of $\HFp$, with only notational changes.
\end{remark}

When $Y$ is an integral homology three-sphere, we can compare the
absolute grading of $\HFc(Y)$ as defined in
Section~\ref{HolDiskFour:sec:AbsGrade} of~\cite{HolDiskFour} with the
absolute $\Zmod{2}$ grading as defined in~\cite{HolDiskTwo}.

Observe first that for an integral homology three-sphere, the absolute
grading of Section~\ref{HolDiskFour:sec:AbsGrade}
of~\cite{HolDiskFour} is actually a $\Z$ lift (rather than a $\Q$
lift) of the usual relative grading. To see this, let $W$ be a a
cobordism containing only two-handles from $S^3$ to $Y$ (used in the
definition of the absolute grading), and observe that
$$\frac{c_1(\spinc)^2-2\chi(W)-3\sgn(W)}{4}\equiv
\frac{c_1(\spinc)^2-\sgn(W)}{4}, \pmod{\Z}$$ (where $\sgn(W)$ is the signature of the intersection form of $W$) which is integral since $c_1(\spinc)$ is a characteristic vector
for the definite, unimodular form $H^2(W;\Z)$.

Next, recall that in general $\HFc(Y,\spinct)$ comes equipped with an absolute
$\Zmod{2}$ grading, defined in
Subsection~\ref{HolDiskTwo:subsec:AbsoluteGradings}
of~\cite{HolDiskTwo}. When $Y$ is an integral homology three-sphere,
this is the $\Zmod{2}$ grading characterized by the property that
$$\chi(\HFa(Y))=1.$$ (Since there is only one $\SpinC$ structure in
this case, we drop it from the notation.)

An equivalent (but technically more useful) formulation of this
grading follows from the calculation of $\HFinf$ of an integer
homology three-sphere, Theorem~\ref{HolDiskTwo:thm:HFinfGen}
of~\cite{HolDiskTwo}, where it is shown that $\HFinf(Y)\cong
\Z[U,U^{-1}]$. We then define the $\Zmod{2}$ degree so that
$\HFinf(Y)$ is supported in even degree.

\begin{prop}
\label{prop:ZModTwoGrading}
If $Y$ is an integral homology three-sphere, then the parity of the
absolute $\Z$ grading defined above is the absolute $\Zmod{2}$ grading
defined in~\cite{HolDiskTwo}. 
\end{prop}

\begin{proof}
Both gradings (mod $2$) agree on $Y$, if and only if they agree on
$Y_1$. This is easily seen from the long exact sequence connecting
$\HFp(Y)$, $\HFp(Y_0)$, and $\HFp(Y_1)$, according to which the
composite map from $Y$ to $Y_1$ shifts the absolute $\Z$-degree by
$-1$ (modulo two); it also shifts the absolute $\Zmod{2}$ degree. The
result then follows from the model case $S^3$, together with fact that
every integer homology can be constructed from $S^3$ by a sequence of
$\pm 1$ surgeries.
\end{proof}

Another case which will interest us is that of three-manifolds $Y_0$
with $H_1(Y_0;\Z)\cong\Z$. For such a manifold, the absolute
$\Zmod{2}$ grading from~\cite{HolDiskTwo}
is defined so that the image of the action by $H_1(Y_0;\Z)$
on $\HFinf(Y_0)$ has degree zero modulo $2$.

\begin{prop}
\label{prop:ZModTwoGradingBOne}
Let $Y_0$ be a three-manifold with $H_1(Y_0;\Z)\cong \Z$, and let $\spinct_0$
be the $\SpinC$ structure with $c_1(\spinct_0)=0$.
Then the parity of $\liftGr+\OneHalf$ is the absolute $\Zmod{2}$ grading given above.
\end{prop}

\begin{proof}
By surgering out the one-dimensional homology, we can fit $Y_0$ into a
long exact sequence for $Y$, $Y_0$, and $Y_1$, where $Y$ and $Y_1$ are
integer homology three-spheres. The result then follows from the
grading in the long exact sequence
(Lemma~\ref{lemma:DegShiftExact}) together with
Proposition~\ref{prop:ZModTwoGrading} above.
\end{proof}

We now turn to the graded refinement of the long exact sequence, after
introducing some notation.  Let $Y$ be a closed, oriented
three-manifold and $a\in\R$, and let $\spinc_0$ be a torsion $\SpinC$
structure. Then, we let $$\HFp_{\leq a}(Y,\spinc_0)\subset
\HFp(Y,\spinc_0)$$ denote the subgroup generated by homogeneous
elements $\xi$ with degree $\liftGr(\xi)\leq a$. Let ${\mathfrak
T}\subset \SpinC(Y)$ denote the subset of torsion $\SpinC$ structures,
then $$\HFp_{\leq a}(Y)=\left(\bigoplus_{\spinc_0\in{\mathfrak T}}
\HFp_{\leq a}(Y,\spinc_0)\right)\oplus
\left(\bigoplus_{\spinct\in \SpinC(Y)-{\mathfrak T}}\HFp(Y,\spinct)\right).$$

We have the following:

\begin{theorem}
\label{thm:TruncExact}
Let $K\subset Y$ be a knot in an integral homology three-sphere. Then, for all sufficiently large $n$, we 
have an exact sequence
$$
\begin{CD}
...@>>>\HFp_{\leq 2n+1}(Y)@>{F_1}>> 
\HFp_{\leq 2n+1}(Y_0) @>{F_2}>>\HFp_{\leq 2n+1}(Y_1)@>{F_3}>>...
\end{CD}
$$
\end{theorem}

\begin{proof}
Observe that $\Ker F_1 \subset \HFp(Y)$ is a finitely generated
$\Z$-module, since for all sufficently large $k$, the $\spinc_0$
component of $F_1$ carries $\HFp_{2k}(Y)$ isomomorphically to
$\HFp_{2k-1/2}(Y_0,\spinc_0)$, and $\HFp_{2k+1}(Y)=0$ (we are using
here the $\HFinf$ characterization of the $\Zmod{2}$ grading and
Proposition~\ref{prop:ZModTwoGrading} to compare it with the absolute
$\Z$ grading). Moreover, this kernel is the image of
$F_{3}(\HFp(Y_1))$.  It follows that we can find an $n$ large enough
that $$F_3(\HFp(Y_1))=F_3(\HFp_{\leq 2n+1}(Y_1)).$$

Moreover, since $\bigoplus_{\spinct\neq\spinc_0}\HFp(Y_0,\spinct)$ is
finitely generated, we have that 
$F_2$ maps $\HFp_{\leq 2n+1}(Y_0)$ into $\HFp_{\leq 2n+1}(Y_1)$,
for all sufficiently large $n$.

The theorem then follows from the
Lemma~\ref{lemma:DegShiftExact}.
\end{proof}

\subsection{More exact sequences: fractional surgeries}

There are analogues of the refined exact sequence (Theorem~\ref{thm:TruncExact}) which
holds for fractional surgeries, as well. 

Recall that in Theorem~\ref{HolDiskTwo:thm:ExactFrac}
of~\cite{HolDiskTwo}, we established an exact sequence for fractional
surgeries of the form 
\begin{equation}
\label{eq:FractionalExact}
\begin{CD}
...@>>>\HFp(Y)@>{\liftF{1}}>> \uHFp(Y_0;\Zmod{q}) @>{\liftF{2}}>>
\HFp(Y_{1/q})@>{\liftF{3}}>>\HFp(Y)@>>>...,
\end{CD}
\end{equation}
where $Y$ is an integral homology sphere, and $q$ is any real
number. In the above sequence,
$$
\uHFp(Y_0;\Zmod{q})=\bigoplus_{\spinct\in\SpinC(Y_0)} \uHFp(Y_0,\spinct;\Zmod{q})
$$ is a sum of $\HFp$-groups twisted by some representation
$$H_1(Y_0;\Z)\longrightarrow \Zmod{q}.$$ Observe that the
$\spinct_0$-component of $\uHFp(Y_0,\Zmod{q})$
(where, as usual, $c_1(\spinct_0)=0$) can be
given an absolute $\Q$-grading (inherited from the absolute grading of
the untwisted group).

Defining $\uHFp_{\leq 2n+1}(Y_0;\Zmod{q})$ with respect to this
grading on the torsion component, we have the following truncated exact sequence
(analogous to 
Theorem~\ref{thm:TruncExactFrac}):

\begin{theorem}
\label{thm:TruncExactFrac}
Let $K\subset Y$ be a knot in an integral homology three-sphere. Then,
for all sufficiently large $n$, there is an exact sequence
$$
\begin{CD}
...@>>>\HFp_{\leq 2n+1}(Y)@>{\liftF{1}}>> 
\uHFp_{\leq 2n+1}(Y_0,\Zmod{q}) @>{\liftF{2}}>>\HFp_{\leq 2n+1}(Y_{1/q})@>{\liftF{3}}>>...
\end{CD}
$$
\end{theorem}

The proof of the above theorem is complicated by the fact that in the
proof of the fractional surgeries exact sequence
(Theorem~\ref{HolDiskTwo:thm:ExactFrac} of~\cite{HolDiskTwo}), the
maps $\liftF{2}$ and $\liftF{3}$ were constructed by counting
pseudo-holomorphic triangles (and thus they correspond to maps of
cobordisms), while the map $\liftF{1}$ is an induced coboundary
map. Thus, our aim is to give a construction of $\liftF{1}$ by
counting pseudo-holomorphic triangles. This construction diverges
slightly from the context set up in
Subsection~\ref{HolDiskFour:subsec:TwoHandles} of~\cite{HolDiskFour},
but the following analogue of Lemma~\ref{lemma:DegShiftExact} still
holds (and from this, Theorem~\ref{thm:TruncExactFrac} follows
easily):

\begin{prop}
\label{prop:DegShiftExactFrac}
Let $K\subset Y$ be a knot in an integral homology three-sphere. 
There is an exact sequence
$$
\begin{CD}
...@>>>\HFp(Y)@>{\liftF{1}}>> \uHFp(Y_0;\Zmod{q}) @>{\liftF{2}}>>\HFp(Y_{1/q})@>{\liftF{3}}>>\HFp(Y)@>>>...,
\end{CD}
$$
with the property that
the component of $\liftF{1}$ mapping into $\HFp(Y_0,\spinct_0;\Zmod{q})$ has degree $-1/2$, 
the restriction of $\liftF{2}$ to $\uHFp(Y_0,\spinct_0;\Zmod{q})$ has degree $-1/2$.
\end{prop} 

To prove the above proposition, we find it convenient to introduce a
construction closely related to twisted coefficients, and to verify
that the associated maps satisfy the naturality properties needed to
establish exactness for the surgery sequence. We return to the proof
of Proposition~\ref{prop:DegShiftExactFrac} and
Theorem~\ref{thm:TruncExactFrac} at the end of the present section,
after this lengthy digression.

\subsection{Covering spaces}
We introduce Floer homology groups associated to covering spaces of
the symmetric product, and establish enough naturality properties
for the above applications (Theorem~\ref{thm:TruncExactFrac}).

Fix a two-manifold $\Sigma$.
A  one-dimensional cohomology class $\xi\in
H^1(\Sigma,\Zmod{n})\cong H^1(\Sym^{g}(\Sigma),\Zmod{n})$
describes an $n$-fold cyclic covering space of $\Sym^g(\Sigma)$, which
we denote by 
$$\Pi_\xi\colon \liftSym\longrightarrow \Sym^g(\Sigma).$$  
We assume that the covering space is
connected, which is equivalent to the condition that there is 
a homology class in $\Sigma$ whose pairing with $\xi$ is one.
Observe that $\xi$ also gives an $n$-fold covering of $\Sigma$, which we
denote by $\liftSigma$, and there is a branched covering space
$$\Sym^g(\liftSigma)\longrightarrow \liftSym.$$

Next, fix a Heegaard triple $(\Sigma,\alphas,\betas,z)$, and let
$\liftalphas$ and $\liftbetas$ be $g$-tuples of simple, closed curves
in $\liftSigma$, with the property that the $\Pi_\xi\liftalpha_i$ is
homotopic to some multiple of $\alpha_i$ and $\Pi_\xi\liftbeta_i$ is
homotopic to some multiple of $\beta_i$.  There are induced tori
$\liftTa$ and $\liftTb$ inside $\liftSym$.

We can define chain complexes analogous to the $\CFm$, $\CFinf$,
$\CFp$, and $\CFa$. For instance, we can define
${\liftCFinf}(\liftalphas,\liftbetas,z)$ generated by pairs
$[\liftx,i]$ where $i\in\Z$ and $\liftx\in\liftTa\cap\liftTb$ so that
$\spinc_z(\Pi(\liftx))$, with $${\widetilde\partial} [\liftx,i] =
\sum_{\lifty\in\liftTa\cap\liftTb}
\sum_{{\widetilde\phi}\in\pi_2(\liftx,\lifty)} \#(\ModFlow(\widetilde\phi))
\cm {[\lifty, i-n_z(\Pi\circ\phi)]}.$$
Here, of course, $\pi_2(\liftx,\lifty)$ denotes the space of homotopy
classes of Whitney disks in $\liftSym$ which connect $\liftx$ to
$\lifty$.  We can also define its quotient complex
${\liftCFp}(\liftalphas,\liftbetas)$, where the integer is
required to be non-negative. 
In fact, for notational simplicity, we
will focus on this case.

The verification that these are, in fact, complexes follows by
modifying the discussion in of~\cite{HolDisk}. We outline the main
steps.

First, of course, we must use families of almost-complex structures
on the lift $\liftSym$ to obtain genericity.  

Next, observe that a Whitney disk $\phi$ in
$\liftSym$ has a corresponding domain $\cald(\phi)$ in $\Sigma$, 
whose boundary lies in the images of the
$\liftalphas$ and $\liftbetas$. Holomorphic curves in $\liftSym$ give
rise to branched covers of the disk, which map holomorphically into
$\Sigma$, as before. The energy bounds of~\cite{HolDisk} follow
directly.

Next, since the map from $\liftSym$ to $\Sym^{g}(\Sigma)$ is a
covering space, we see that
$\pi_2(\liftSym)\cong\pi_2(\Sym^{g}(\Sigma))$, and also that
$T\liftSym\cong \Pi_\xi^* T\Sym^{g}(\Sigma)$.  
Thus, the dimension counts from~\cite{HolDisk} can be used to show that
generic choices of $\liftalphas$ and $\liftbetas$, two-dimensional
moduli spaces miss the holomorphic spheres (which could otherwise
spoil ${\widetilde\partial}^2=0$).

For boundary degenerations, observe that our homological hypotheses
show that for any $\phi\in\pi_2({\widetilde\x},{\widetilde\x})$, 
$\cald(\phi)$ consists of some number of copies of $\Sigma$. 

It then as before that  that the only types of boundary degenerations
which can occur (in the proof that $\partial^2=0$)
are those in the homotopy class 
$O_{\widetilde\x}+S\in \pi_2({\widetilde\x},{\widetilde\x})$. 
As before, an $\sj$-holomorphic representative for this homotopy class
is generically injective, and hence, spaces of boundary degenerations
are smooth for generic perturbations of the almost-complex structure
induced over $\liftSym$. The total number of these can be calculated
by degenerating the base $\Sigma$, to see that the signed count is zero.

Finally, when there are relations between the homology classes in the
span of $\liftalphas$ and those of $\liftbetas$, then we need
corresponding admissibile hypotheses to show that the sums are {\em a
priori} finite. Weak admissibility in this context, for instance,
(which is sufficient for the purposes of $\HFp$) is the requirement
that all homologies between the spans of $\liftalphas$ and
$\liftbetas$, when projected down to $\Sigma$, always have both
positive and negative multiplicities.

With these remarks in place, then, we have ``lifted'' a chain complex
$\liftCFp(\Sigma,\liftalphas,\liftbetas,\liftz)$, and also its
corresponding homology group.

We relate this construction directly to the three-manifolds invariants
in two cases.  Fix a Heegaard triple $(\Sigma,\alphas,\betas,z)$ for
an oriented three-manifold $Y$, and also a one-dimensional cohomology
class $\xi\in H^1(\Sigma;\Z)$. Fix lifts $\liftalphas$ and $\liftbetas$
for $\alphas$ and  $\betas$ under the map $\Pi_\xi$.
In this case, $\liftCFp(\liftalphas,\liftbetas)$ splits into summands
indexed by $\spinct\in\SpinC(Y)$, 
where $\liftCFp(\liftalphas,\liftbetas,\spinct)$ is generated by
$\liftx$ with $\spinc_z(\Pi_\xi\liftx)=\spinct$.

There will be two subcases of particular importance to us. 
First, suppose
that $\xi\in H^1(\Sigma,\Zmod{p})$ satisfies the following two
conditions
\begin{itemize}
\item there is an element in the span of the $\alphas$ whose evaluation on 
$\xi$ is $1$
\item the restriction of $\xi$ to the span of the $\betas$ is trivial.
\end{itemize}
In this case, $\liftTa=\Pi_\xi^{-1}(\Ta)$, while $\liftTb$ is one
of the $n$ tori which project to $\Tb$.
\begin{prop}
\label{prop:DisconnectedCase}
Fix a one-dimensional cohomology class $\xi$ with the property that
the induced covering space of $\Ta$ is connected, while the
restriction of $\xi$ to the span of the $\betas$ is trivial. 
Then, for any choice of $\liftbetas$ as above, we have 
a natural identification
$$\liftHFp(\liftalphas,\liftbetas,\spinct)\cong
\HFp(\alphas,\betas,\spinct).$$
\end{prop}

\begin{proof}
According to the assumption, $\Pi_{\xi}^{-1}(\Tb)$ consists of $n$ disjoint
tori. Thus, each intersection point $\x\in\Ta\cap\Tb$ corresponds to a
unique intersection point of $\liftTa$ with $\liftTb$. Moreover, if
$\liftx,\lifty\in\liftTa\cap\liftTb$ lie over $\x,\y\in\Ta\cap\Tb$,
the usual covering space theory gives an identification between
Whitney disks connected $\x$ to $\y$ with Whitney disks connecting
$\liftx$ to $\lifty$. By using the lift of a family of almost-complex
structures over $\Sym^g(\Sigma)$, it is easy to see that the boundary
maps are identified, as well.

Thus, the identification takes place actually on the chain complex
level.
\end{proof}

For the second case, we start again with a pointed Heegaard diagram
$(\Sigma,\alphas,\betas,z)$ for $Y$, and a class $\xi\in
H^1(\Sigma;\Zmod{n})$, only now we suppose that both covering spaces
$\liftTa=\Pi_\xi^{-1}(\Ta)$ and $\liftTb=\Pi_\xi^{-1}(\Tb)$ are
connected (i.e. there are homology classes in the spans of
both $[\alpha_1],...,[\alpha_g]$ and $[\beta_1],...,[\beta_g]$ whose pairing
with $\xi$ is one).

Observe that the covering space $\liftSym$ gives rise to an additive
map $$A \colon \pi_2(\x,\y) \longrightarrow \Zmod{n},$$ in the
following way. For each $\x\in\Ta\cap\Tb$, choose some lift
$\liftx_0\in\liftTa\cap\liftTb$. Next, each $\phi\in\pi_2(\x,\y)$, has
a unique lift starting at ${\widetilde\x}_0$. Its other endpoint ${\widetilde\y}$
(which
depends only on the homotopy class of $\phi$) lies over $\y$; thus we can define
the element
$A(\phi)\in\Zmod{n}$ to be the element with the property that
${\widetilde\y}=A(\phi)\cm \lifty_0$. It is easy to see that this $A$ is additive
under juxtaposition, hence it is an additive assignment of the kind
required to define homology with twisted coefficients $\HFp(Y,\spinc;\Zmod{n})$.
(c.f. Subsection~\ref{HolDiskTwo:subsec:TwistedCoeffs}
of~\cite{HolDiskTwo}). 

\begin{prop}
\label{prop:ConnectedCase}
Let $(\Sigma,\alphas,\betas,z)$ be a Heegaard diagram for $Y$, and fix
a class $\xi$ as above, with the property that the induced covering
spaces for both $\Ta$ and $\Tb$ are connected. Then, for the
additive assignment defined above, there is a natural identification
$$\uHFp(Y,\spinct;\Zmod{n})\cong
\liftHFp(\liftTa,\liftTb;\spinct).$$
\end{prop}

\begin{proof}
We have an identification of $\Zmod{n}$-sets:
$$\liftTa\cap\liftTb\cong (\Ta\cap\Tb)\times\Zmod{n}.$$ Thus, the
generators of $\uCFp(Y,\spinct)$ and $\liftCFp(\Ta,\Tb;\spinct)$ are
identified.  The boundary maps are easily seen to coincide.
\end{proof}

Of course, if we move the $\betas$ by an isotopy which does not cross
the basepoint, there is an induced isotopy amongst the $\liftbetas$,
which does not project to the basepoint -- hence leaving the homology
groups unchanged. But we have more freedom to move the
$\liftbetas$. 

\begin{prop}
\label{prop:IsotopyLift}
Suppose that $\{\zeta_1,...,\zeta_g\}$ are curves in $\liftSigma$
which are simultaneously isotopic to $\liftbetas$ via an isotopy which
never projects to the basepoint $z\in\Sigma$, then the isotopy induces
an identification $$\liftCFp(\liftTa,\liftTb)\cong
\liftCFp(\liftTa,\Tz),$$ where the second chain complex is defined
using the torus $\Tz=\zeta_1\times...\times\zeta_g$.
\end{prop}

\begin{proof}
We adapt the usual proof of isotopy invariance,
noting that the the branched
covering $$\liftSigma^{\times g}\longrightarrow
\Sym^g(\liftSigma)\longrightarrow
\liftSym$$
gives us a bound of the energy of holomorphic disks in $\liftSym$
in terms of the induced homology class in $\liftSigma^{\times g}$.
This energy remains bounded under dynamical boundary conditions,
provided that the $\liftbetas$ are moved by an exact
Hamiltonian in $\liftSigma$. 
\end{proof}

\subsection{The fractional surgeries long exactness}

We can now modify the
proof of the long exact sequence for fractional surgeries
(Theorem~\ref{HolDiskTwo:thm:ExactFrac} of~\cite{HolDiskTwo})
to have a better understanding of
the maps involved. We assume for simplicity that $Y$ is an integer
homology three-sphere, and $K\subset Y$ is a knot
with the canonical framing, so that in particular $b_1(Y_0)=1$.

Consider the Heegaard quadruple
$(\Sigma,\alphas,\betas,\gammas,\deltas,z)$ as in
Section~\ref{HolDiskTwo:subsec:FracSurg} of~\cite{HolDiskTwo}.  In
particular, $(\Sigma,\alphas,\betas,\gammas,z)$ is a Heegaard triple
subordinate to the knot $K\subset Y$ (with its $0$-framing),
$(\Sigma,\alphas,\gammas,\deltas,z)$ is a Heegaard triple 
subordinate to a knot in $Y_0$, which in turn
represents the a canonical cobordism from $Y_0$ to $Y_{1/q}$, 
while $(\Sigma,\alphas,\deltas,\betas)$ represents a cobordism
from the union of $Y_{1/q}$ and a lens space to $Y$.

Let $\xi\in H^1(\Sigma;\Z/q\Z)$ be the Poincar\'e dual of
$\beta_g$. Clearly, $\xi$ evaluates trivially on the span of the
$\betas$ and $\deltas$. However, $\langle \xi, \gamma_g\rangle=1$. It
follows that there is some element in the span of $\alphas$ whose
pairing with $\xi$ is one.

Since in $H_1(\Sigma;\Z)$ we have the relation
$$\beta_g + q \gamma_g = \delta_g,$$
we have the corresponding relation in $H_1(\liftSigma;\Z)$:
$${\widetilde \beta}_g + {\widetilde\gamma}_g = {\widetilde\delta}_g.$$

Observe that we have isomorphisms
$$\HFleq(\liftbetas,\liftgammas)
\cong\HFleq(\liftdeltas,\liftbetas)\cong \HFleq(\#^{g-1}(S^1\times S^2)),$$
in view of Proposition~\ref{prop:DisconnectedCase}, so these first two
groups have (up to sign) canonical top-dimensional generators which we
denote by $\liftTheta_{\beta,\gamma}$ and $\liftTheta_{\delta,\beta}$
respectively. Similarly, we have (up to sign) a canonical top dimensional generator
$\liftTheta_{\gamma,\delta}$ for
$$\HFleq(\liftgammas,\liftdeltas)\cong
\HFleq(\#^{g}(S^2\times S^1),\Zmod{q})\cong \HFleq(\#^{g}(S^2\times S^1)).$$
(The first isomorphism is an application of
Proposition~\ref{prop:ConnectedCase}, and the second is an easy
calculation.)

Define
$$
\begin{CD}
{\widetilde F_1}\colon \HFp(Y)\cong \liftHFp(\liftTa,\liftTb)
@>{\otimes {\liftTheta}_{\beta,\gamma}}>>\liftHFp(\liftTa,\liftTc)\cong
\uHFp(Y_0;\Zmod{n}),
\end{CD}
$$ 
Where the first isomorphism is induced from
Proposition~\ref{prop:DisconnectedCase}, the second isomorphism
is induced from
Proposition~\ref{prop:ConnectedCase}, and the tensor product map 
is defined by counting (zero-dimensional spaces of)
pseudo-holomorphic triangles in $\liftSym$.
Similarly, we define 
$$
\begin{CD}
{\widetilde F_2}\colon 
\uHFp(Y_0;\Zmod{n})\cong
\liftHFp(\liftTa,\liftTc)
@>{\otimes {\liftTheta}_{\gamma,\delta}}>>
\liftHFp(\liftTa,\liftTc)\cong \HFp(Y_{1/q})
\end{CD}
$$ 
Finally, $\liftF{3}$ is defined as it is in~\cite{HolDiskTwo}
(where it is simply denoted $F_3$) by 
$$
\begin{CD}
{\widetilde F_3}\colon 
\HFp(Y_{1/q})\cong
\HFp(\Ta,\Td)
@>{\otimes \Theta_{\delta,\beta}}>>
\HFp(\Ta,\Tb)\cong 
\HFp(Y).
\end{CD}
$$ 

We now have the precise statement of the fractional long exact sequence:

\begin{prop}
\label{prop:FunctorialExactFrac}
The maps ${\widetilde F}_1$, ${\widetilde F}_2$, and ${\widetilde
F}_3$ defined above fit into a long exact sequence:
$$
\begin{CD}
... @>>>
\HFp(Y)	
@>{\liftF{1}}>>
\uHFp(Y_0;\Zmod{q})
@>{\liftF{2}}>>
\HFp(Y_{1/q})
@>{\liftF{3}}>>
...
\end{CD}
$$
\end{prop}

\begin{proof}
Use a pointed Heegaard multi-diagram
$$(\Sigma,\alphas,\betas,\gammas,\deltas,z)$$ 
as above, so that  $(\alphas,\betas)$ describes
$Y$, $(\betas,\gammas)$ describes $Y_0$, $(\alphas,\deltas)$
describes $Y_1$. Also, fix the covering space of $\liftSym$ described above, and
lifts $\liftalphas$, $\liftbetas$, $\liftgammas$, and $\liftdeltas$ of the corresponding curves.

Observe that there is a smooth isotopy taking ${\widetilde\gamma}_g$ arbitrarily
close to the
juxtaposition of ${\widetilde\beta}_g$ with ${\widetilde\delta}_g$,
which does not project to the basepoint $z\in\Sigma$ (see
Figure~\ref{fig:IsotopyLift}). We denote the induced isotopy of
$\liftTc$ by $\Psi_t$ and new isotopic $\gamma$-torus ($\Psi_1(\Tc)$)
by $\liftTcPr$.  The proof of Theorem~\ref{HolDiskTwo:thm:ExactOne}
of~\cite{HolDiskTwo} (using the filtrations when ${\widetilde\gamma}_g$ is a
sufficiently close approximation to the juxtaposition of
${\widetilde\beta}_g$ and ${\widetilde\delta}_g$)
then gives a exactness (in the middle) for the
maps
\begin{equation}
\label{eq:ExactnessAfterIsotopy}
\begin{CD}
\liftHFp(\liftalphas,\liftbetas) @>{\otimes
\liftTheta_{\beta,\gamma'}}>> \liftHFp(\liftalphas,\liftgammasPr) @>{\otimes
\liftTheta_{\gamma',\delta}}>> \liftHFp(\liftalphas,\liftdeltas).
\end{CD}
\end{equation}
Here, $\otimes \liftTheta_{\beta,\gamma'}$ is shorthand for the map
$$\eta\mapsto
\sum_{\spinc\in X_{\alpha,\beta,\gamma'}}\pm {\widetilde f}_{\alpha,\beta,\gamma'}(\eta\otimes
\liftTheta_{\beta,\gamma'},\spinc),$$
(with an appropriate choice of signs) where, of course, ${\widetilde f}_{\alpha,\beta,\gamma'}$ counts
all holomorphic triangles in $\liftSym$. Also, 
$\liftTheta_{\beta,\gamma'}\in\HFleq(\liftbetas,\liftgammas,\spinc_0)$ is a generator which was
explicitly written down in~\cite{HolDiskTwo}. It has an alternative,
more algebraic characterization (up to a
sign) as the generator of $\HFleq(\liftbetas,\liftgammas,\spinc_0)\cong 
\HFleq(\#^{g-1}(S^2\times S^1),\spinc_0)$ 
with maximal (relative) degree.

\begin{figure}
\mbox{\vbox{\epsfbox{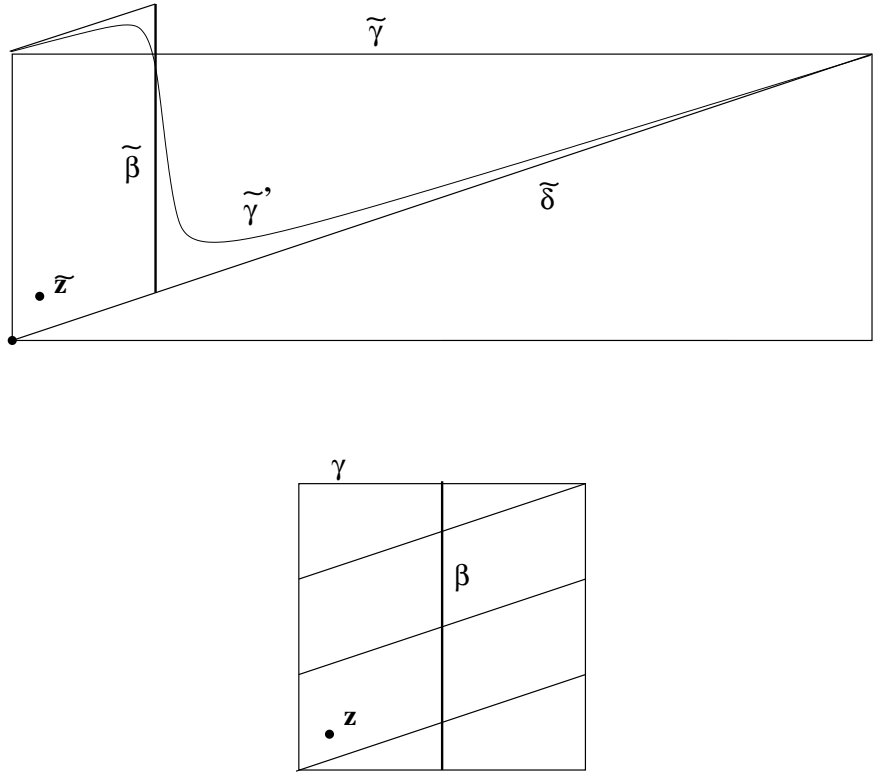}}}
\caption{\label{fig:IsotopyLift}
{\bf{Isotopy of ${\widetilde\gamma}$}} We illustrate the isotopy of
${\widetilde\gamma}$ for the $1/q$ surgery exact sequence, when
$q=3$. The lower picture takes place in a genus one surface, the upper
picture takes place in its triple cover (which unwinds $\gamma$ to
give ${\widetilde\gamma}$. The curve ${\widetilde \gamma}'$
approximates the juxtaposition of the lifts ${\widetilde\beta}$ and
${\widetilde\delta}$ and is
isotopic to ${\widetilde\gamma}$ through an isotopy which does not project to the basepoint.}
\end{figure}

We claim that functoriality of the 
triangle construction in $\liftSym$, together with exactness in the middle 
for Exact Sequence~\eqref{eq:ExactnessAfterIsotopy},
gives exactness in the middle for:
\begin{equation}
\begin{CD}
\liftHFp(\liftalphas,\liftbetas) @>{\otimes
\liftTheta_{\beta,\gamma}}>> \liftHFp(\liftalphas,\liftgammas) @>{\otimes
\liftTheta_{\gamma,\delta}}>> \liftHFp(\liftalphas,\liftdeltas).
\end{CD}
\end{equation}
Observe that the maps here correspond to the maps $\liftF{1}$ and $\liftF{2}$
from the statement of the proposition.

To see this, note that the isotopy of ${\widetilde\gamma}_g$ to
${\widetilde\gamma}_g'$ gives rises to isomorphisms
\begin{eqnarray*}
\liftGamma_{\alpha;\gamma,\gamma'}\colon \liftHFp(\liftalphas,\liftgammas) &\longrightarrow&
\liftHFp(\liftalphas,\liftgammasPr), \\
\liftGamma_{\beta;\gamma,\gamma'}\colon \liftHFp(\liftbetas,\liftgammas) &\longrightarrow&
\liftHFp(\liftbetas,\liftgammasPr), \\
\liftGamma_{\gamma,\gamma';\delta}\colon \liftHFp(\liftgammas,\liftdeltas) &\longrightarrow&
\liftHFp(\liftgammasPr,\liftdeltas)
\end{eqnarray*} 
defined analogously to the maps of isotopies downstairs.
It is straightforward to see that these are functorial for triangles (compare
Theorem~\ref{HolDiskFour:thm:Triangles} of~\cite{HolDiskFour}). In particular, we have that
\begin{eqnarray*}
{\liftF{\alpha,\beta,\gamma'}}(\xi\otimes \liftGamma_{\beta;\gamma,\gamma'}(\eta))&=&
\liftGamma_{\alpha;\gamma,\gamma'}\circ
{\liftF{\alpha,\beta,\gamma}}(\xi\otimes\eta) \\
{\liftF{\alpha,\beta',\gamma}}(\liftGamma_{\alpha;\beta,\beta'}(\xi)\otimes
\liftGamma_{\beta,\beta';\gamma}(\eta)) &=& {\liftF{\alpha,\beta,\gamma}}(\xi\otimes\eta).
\end{eqnarray*}
It follows that the following diagram is commutative:
\begin{equation}
\label{eq:TriangleDiagram}
\begin{CD}
\liftHFp(\liftalphas,\liftbetas) @>{\otimes\liftTheta_{\beta,\gamma}}>> \liftHFp(\liftalphas,\liftgammas) @>{\otimes
\liftTheta_{\gamma,\delta}}>> \liftHFp(\liftalphas,\liftdeltas) \\
@V{\Id}VV @V{\liftGamma_{\alpha;\gamma,\gamma'}}VV @V{\Id}VV \\
\liftHFp(\liftalphas,\liftbetas) @>{\otimes
\liftGamma_{\beta;\gamma,\gamma'}(\liftTheta_{\beta,\gamma}})>> 
\liftHFp(\liftalphas,\liftgammasPr) 
@>{\otimes\liftGamma_{\gamma,\gamma';\delta}(\liftTheta_{\gamma,\delta})}>>
\liftHFp(\liftalphas,\liftdeltas).
\end{CD}
\end{equation} 
Again, the generators $\liftTheta_{\beta,\gamma}$ and
$\liftTheta_{\gamma,\delta}$ are uniquely characterized (up to sign) by
the property that they are top-dimensional dimensional generators. It follows (since
$\liftGamma_{\beta;\gamma,\gamma'}$ and $\liftGamma_{\gamma,\gamma';\delta}$
are isomorphisms on the level of homology) that
\begin{eqnarray*}
\liftGamma_{\beta;\gamma,\gamma'}(\liftTheta_{\beta,\gamma})=\pm \liftTheta_{\beta,\gamma'}
&{\text{and}}&
\liftGamma_{\gamma,\gamma';\delta}(\liftTheta_{\gamma,\delta})=\pm 
\liftTheta_{\gamma',\delta}
\end{eqnarray*}
Thus, Theorem~\ref{HolDiskTwo:thm:ExactOne} establishes
exactness along the bottom row of Diagram~\eqref{eq:TriangleDiagram}, which in turn
implies exactness along the top row of the same diagram.

Of course, by repeating these arguments 
only now isotoping the curve $\delta_g$ to 
$\gamma_g+\beta_g$ (which we can do downstairs in $\Sigma$), we establish 
exactness for the maps:
$$
\begin{CD}
\liftHFp(\liftalphas,\liftgammas) @>{\otimes{\liftTheta_{\gamma,\delta'}}}>> \liftHFp(\liftalphas,\liftdeltasPr) 
@>{\otimes{\liftTheta_{\delta',\beta}}}>> \liftHFp(\liftalphas,\liftbetas).
\end{CD}
$$
(Indeed, this was the isotopy used for the case of fractional surgeries in 
Theorem~\ref{HolDiskTwo:thm:ExactFrac} of~\cite{HolDiskTwo}); and hence, by the
same naturality arguments, we get exactness for 
$$
\begin{CD}
\liftHFp(\liftalphas,\liftgammas) @>{\otimes{\liftTheta_{\gamma,\delta}}}>> \liftHFp(\liftalphas,\liftdeltas) 
@>{\otimes{\liftTheta_{\delta,\beta}}}>> \liftHFp(\liftalphas,\liftbetas).
\end{CD}
$$
\end{proof}

\subsection{Proof of the truncated exact sequence for fractional surgeries}

Proposition~\ref{prop:DegShiftExactFrac} follows easily from
Proposition~\ref{prop:FunctorialExactFrac}.

\vskip.3cm
\noindent{\bf Proof of Proposition~\ref{prop:DegShiftExactFrac}.}
We use the maps $\liftF{1}$ and $\liftF{2}$ as defined in the proof of
Proposition~\ref{prop:FunctorialExactFrac} above. Clearly, the
holomorphic triangles which comprise $\liftF{1}$ (before we move
$\TcPr$ by an isotopy) project down to holomorphic triangles for the
Heegaard triple $(\Sigma,\alphas,\betas,\gammas,z)$. Moreover, it is
easy to see that the Maslov indices of these triangles agree with the
Maslov indices of their projections. Thus, the degree shift of the
$\spinct_0$-component of $\liftF{1}$ must agree with the degree shift
of its projection. But the projection represents the natural cobordism
from $Y$ to $Y_0$, so this degree shift is $\OneHalf$. The same
remarks apply to the map $\liftF{2}$.
\qed
\vskip.3cm

\vskip.3cm
\noindent{\bf{Proof of Theorem~\ref{thm:TruncExactFrac}}}
This follows immediately from Proposition~\ref{prop:DegShiftExactFrac}
in the same way that Theorem~\ref{thm:TruncExact} followed from
Lemma~\ref{lemma:DegShiftExact}.
\qed

\section{The correction term}
\label{sec:CorrTerm}

With the help of the absolute grading, we can define the following
numerical invariant for integer homology three-spheres or, more
generally, rational homology three-spheres and $\SpinC$ structures:

\begin{defn}
Let $Y$ be a rational homology three-sphere.
The {\em correction term $d(Y,\spinct)$} 
is the minimal grading ($\liftGr$) of any non-torsion
element in the image of $\HFinf(Y,\spinct)$ in $\HFp(Y,\spinct)$.
\end{defn}

This invariant is an analogue of a gauge-theoretic invariant defined
by Fr{\o}yshov, see~\cite{Froyshov}. As we shall see (like its
gauge-theoretic analogue), the invariant contains information about
the intersection forms of four-manifolds which bound $Y$.

\begin{prop}
\label{prop:CorrTermOrient}
Let $(Y,\spinct)$ be a rational homology three-sphere. Then, we have that
$$d(Y,\spinct)=d(Y,{\overline\spinct})$$
and
\begin{equation}
\label{eq:DFlipOrientation}
d(Y,\spinct)=-d(-Y,\spinct).
\end{equation}
\end{prop}

\begin{proof}
The conjugation invariance of the correction term follows from the conjugation
invariance of the maps associated to cobordisms 
(c.f Theorem~\ref{HolDiskFour:thm:Conjugation} of~\cite{HolDiskFour}).

We verify Equation~\eqref{eq:DFlipOrientation}.

Consider the natural long exact sequence connecting $\HFm(Y,\spinct)$,
$\HFinf(Y,\spinct)$ and $\HFp(Y,\spinct)$ (Exact
Sequence~\eqref{eq:HFinfExactSequence}).  From this sequence, together
with the fact that $\HFinf(Y,\spinct)\cong \Z[U,U^{-1}]$ (c.f.
Theorem~\ref{HolDiskTwo:thm:HFinfGen} of~\cite{HolDiskTwo}), it
follows that if we define $d^-(Y,\spinct)$ to be the maximal $k$ for
which the map $\iota_k\colon \HFm_k(Y,\spinct)\longrightarrow
\HFinf_k(Y,\spinct)$ is non-trivial, then
\begin{equation}
\label{eq:CompareDs}
d^-(Y,\spinct)=d(Y,\spinct)-2.
\end{equation}

Now there is a duality map $\Duality$
which gives a map from $\HFc$ homology of
$Y$ to $\HF_{\circ}$-cohomology of $-Y$. In its precise graded version,
(c.f. Proposition~\ref{HolDiskFour:prop:GradedDuality} of~\cite{HolDiskFour}),
we obtain a commutative diagram 
$$\begin{CD}
\HFinf_{k}(Y,\spinct)@>{\pi_k}>>
\HFp_{k}(Y,\spinct) \\
@VV{\Duality^\infty}V	@VV{\Duality^+}V \\
\Z\cong \HF_{\infty}^{-k-2}(-Y,\spinct)@>{\iota^{-k-2}}>>
\HF_{-}^{-k-2}(-Y,\spinct)
\end{CD}.$$
(Here, $\Duality[\x,i]=[\x,-i-1]^*$ is the map from the chain complexes 
for $Y$, $(\Sigma,\alphas,\betas,z)$, 
to the cochain complex for $-Y$, $(-\Sigma,\alphas,\betas,z)$.)
The bottom row is the map on cohomology, and the vertical maps are all
isomorphisms. Now, by the universal coefficient theorem in cohomology
and the fact that $\HFinf(Y,\spinct)$ is a free module in each
dimension, it follows that the image of $\iota^{-k-2}$ has rank one if
and only if the map $$\iota_{-k-2}\colon
\HFm_{-k-2}(-Y,\spinct)\longrightarrow \HFinf_{-k-2}(-Y,\spinct)$$ is
non-trivial. In view of this, and Equation~\eqref{eq:CompareDs},
Equation~\eqref{eq:DFlipOrientation} follows.
\end{proof}

We conclude the subsection with the proof of the additivity of the
correction term under the connected sum operation:

\begin{theorem}
\label{thm:AdditivityOfD}
If $(Y,\spinct)$ and $(Z,\spincu)$ are rational homology
three-spheres equipped with $\SpinC$ structures, then
$$d(Y\# Z,\spinct\#\spincu)=
d(Y,\spinct)+d(Z,\spincu).$$
\end{theorem}

The proof occupies the rest of the present subsection.

We define first a natural transformation
$$\Fc_{Y\#Z,\spinct\#\spincu}\colon \HFc(Y,\spinct)\otimes \HFleq(Z,\spincu)\longrightarrow \HFc(Y\# Z,\spinct\#\spincu)$$
as follows.
Observe that for some $g_1$ and $g_2$ there is a pair-of-pants cobordism from
$Y\#(\#^{g_2}(S^2\times S^1))\coprod (\#^{g_1}(S^2\times S^1))\# Z$
to $Y\# Z$. To define this, consider pointed Heegaard diagrams
$(\Sigma_1,\alphas_1,\betas_1,z_1)$ and
$(\Sigma_2,\alphas_2,\betas_2,z_2)$ for $Y$ and $Z$ respectively. Then
the Heegaard triple
$$(\Sigma_1\#\Sigma_2,\alphas_1\alphas_2, \betas_1\alphas_2,
\betas_1\betas_2,z)$$
describes such a cobordism. In the above notation, $\alphas_1\alphas_2$ denotes
the 
union of $\alphas_1$ and $\alphas_2$, thought of now as circles in $\Sigma_1\#\Sigma_2$.
Now, we define $\Fc_{Y\#Z,\spinct\#\spincu}$ to be the composite of 
the map 
$$\HFc(Y,\spinct)\otimes \HFleq (Z,\spincu)
\longrightarrow \HFc(Y\#^{g_2}(S^2\times S^1)) \otimes
\HFleq ((\#^{g_1}(S^2\times S^1))\# Z) $$
induced by the cobordism of one-handles with the map 
$$
\HFc(Y\#^{g_2}(S^2\times S^1)) \otimes
\HFleq ((\#^{g_1}(S^2\times S^1))\# Z) 
\longrightarrow
\HFc(Y\# Z,\spinct\#\spincu)
$$
defined by counting holomorphic triangles in the Heegaard triple.
(Of course, we  perturb the circles to achieve admissibility.)

\begin{prop}
\label{prop:ConnSumTransformation}
The map $\Fc_{Y\# Z,\spinct\#\spincu}$ defined above is independent of
the Heegaard diagrams used for $Y$ and $Z$. Moreover, if $W$ is a
cobordism from $Y$ to $Y'$ equipped with $\SpinC$ structure $\spinc$,
then the following diagram is commutative: $$\begin{CD}
\HFc(Y,\spinct)\otimes \HFleq(Z,\spincu) @>{F_{Y\# Z,\spinct\#\spincu}}>>
\HFc(Y\# Z,\spinct\#\spincu) \\
@V{F_{W,\spinc}\otimes \Id}VV @VV{F_{W\#(Z\times [0,1]),\spinc\#\spincu}}V \\
\HFc(Y',\spinct')\otimes \HFleq (Z,\spincu) @>{F_{Y'\# Z,\spinct'\#\spincu}}>>
\HFc(Y'\# Z,\spinct'\#\spincu).
\end{CD}$$
\end{prop}

\begin{proof}
For simplicity, we consider the case where $W$ is composed entirely of
two-handles. Then, the commutative diagram follows from associativity
of the triangle construction. More specifically, suppose that $W$ is
represented by the Heegaard triple
$(\Sigma_1,\alphas_1,\betas_1,\gammas,z_1)$. Then, by associativity, the following
diagram commutes:
$$
\begin{CD}
\HFc(\alphas_1,\betas_1')\otimes \HFleq(\alphas_2',\betas_2)  \\
@VVV \\
\HFc(\alphas_1\alphas_2,\betas_1'\alphas_2')
\otimes \HFleq(\betas_1'\alphas_2',\betas_1\betas_2)
@>{\otimes}>>
\HFc(\alphas_1\alphas_2,\betas_1\betas_2) \\
@V{\Id\otimes(\cdot\otimes \Theta_{\beta_1\beta_2,\gamma_1\beta_2'})}VV
@VV{\otimes \Theta_{\beta_1\beta_2,\gamma_1\beta_2'}}V \\
\HFc(\alphas_1\alphas_2,\betas_1'\alphas_2')\otimes
\HFleq(\betas_1'\alphas_2',\gammas_1\betas_2') 
@>{\otimes}>> 
\HFc(\alphas_1\alphas_2,\gammas_1\betas_2')
\end{CD}
$$ The first map corresponds to the one-handles. Observe that
$\Theta_{\beta_1\beta_2,\gamma_1\beta_2'}$ is the canonical generator
for $\HFleq(\betas_1\betas_2,\gammas_1\betas_2')$, which describes a
connected sum of $S^1\times S^2$. Thus, going around the above diagram
the two ways corresponds to the two possible compositions.

The case of one- and three-handles follows easily (compare
Section~\ref{HolDiskFour:sec:ECobord}
of~\cite{HolDiskFour}). Independence of the Heegaard diagrams follows
in the same manner.
\end{proof}

\begin{lemma}
\label{lemma:IdMap}
Let $Z=S^3$, and fix the canonical element $\Theta\in\HFleq(S^3,\spinct_0)$.
Then
$$\Fc_{Y\# S^3,\spinct\#\spinct_0}(\cdot\otimes \Theta)\colon  \HFc(Y,\spinct)\longrightarrow
\HFc(Y\# S^3,\spinct\#\spinct_0)\cong \HFc(Y,\spinct)$$
is the identity map.
\end{lemma}

\begin{proof}
This follows from the usual stabilization invariance of the triangle construction.
\end{proof}

\begin{lemma}
\label{lemma:BasicCase}
Let $(Z,\spincu)$ be a rational homology sphere equipped with a
$\SpinC$ structure.  Let $\Theta_Z\in\HFleq(Z,\spincu)$ be an element
whose image in $\HFinf(Z,\spincu)$ under the natural map is a
generator. Then, the map $$\Finf{S^3\#
Z,\spinct_0\#\spincu}(\cdot\otimes \Theta_Z)\colon
\HFinf(S^3,\spinct_0)\longrightarrow
\HFinf(Z,\spincu)$$
is an isomorphism which carries the canonical element of $\HFinf(S^3,\spinct_0)$
to the image of $\Theta_Z$ in $\HFinf(Z,\spincu)$.
\end{lemma}

\begin{proof} By naturality, we have a commutative diagram:
$$\begin{CD}
\HFleq(S^3) @>{\otimes\Theta_Z}>> \HFleq(Z) \\
@V{\iota_{S^3}}VV @VV{\iota_Z}V \\
\HFinf(S^3) @>{\otimes\Theta_Z}>>\HFinf(Z)
\end{CD}$$
It follows from this, and Lemma~\ref{lemma:IdMap} that
$$\Finf{S^3\#Z,\spinct_0\#\spincu}({\iota_{S^3}}(\Theta_{S^3})\otimes \Theta_Z)
=\iota_Z(\Theta_Z).$$
The result follows, with the observation that both
$\HFinf(S^3,\spinct_0)$ and $\HFinf(Z,\spincu)$
are generated as $\Z[U,U^{-1}]$ modules 
by the canonical elements $\iota_{S^3}(\Theta_{S^3})$ and $\iota_Z(\Theta_Z)$
\end{proof}

\begin{prop}
\label{prop:ConnSumShift}
Let $\Theta_Z\in\HFleq(Z,\spincu)$ be an element whose image in
$\HFinf(Z,\spincu)$ under the natural map is a generator. Then, 
$$\Finf{Y\# Z,\spinct\#\spincu}(\cdot\otimes \Theta_Z)\colon 
\HFinf(Y,\spinct)\longrightarrow
\HFinf(Y\#Z,\spinct\#\spincu)$$
is an isomorphism of relatively graded groups. Moreover,
$$\liftGr \left(\Finf{Y\# Z,\spinct\#\spincu}(\xi\otimes \Theta_Z)\right)
= \liftGr(\xi)+\liftGr(\Theta_Z).$$
\end{prop}

\begin{proof}
Each integer homology three-sphere can be obtained from $S^3$ by a
sequence of $\pm 1$ surgeries. When $Y$ is an integer homology sphere,
the result follows from induction on the length of such a sequence,
with Lemma~\ref{lemma:BasicCase} as the base case (when $Y\cong S^3$).

For the inductive step, suppose that the result holds for an integer
homology sphere $Y$. We claim that for each knot $K\subset Y$, the
result also holds for $Y_1=Y_{\pm 1}(K)$. We concentrate on the case
where the sign is $+1$. Then, there is a map of long exact sequences:
\begin{equation}
\label{eq:MapOfLongExacts}
\begin{CD}
...@>>>\HFp(Y)@>>>\HFp(Y_0)@>>>\HFp(Y_1)@>>>... \\
&&@V{\Fp{Y\# Z}(\cdot\otimes\Theta_Z)}VV 
@V{\Fp{Y_0\# Z}(\cdot\otimes\Theta_Z)}VV 
@V{\Fp{Y_1\# Z}(\cdot\otimes\Theta_Z)}VV \\
...@>>>\HFp(Y\# Z,[\spincu])@>>>\HFp(Y_0\# Z,[\spincu])@>>>
\HFp(Y_1\# Z,[\spincu])@>>>... \\
\end{CD}
\end{equation}
We use the convention that if $M$ is a three-manifold, then
$$\HFc(M\# Z,[\spincu])=
\bigoplus_{\{\spinct\in\SpinC(M\# Z)|\spinct|Z=\spincu\}}
\HFc(M\# Z,\spinct),$$
and the vertical maps are obtained by summing the maps defined
in Proposition~\ref{prop:ConnSumTransformation}. Indeed, the
above squares commute by
Proposition~\ref{prop:ConnSumTransformation}. It follows immediately
that $\Finf{Y_0\# Z}(\cdot\otimes \Theta_Z)$ sends the image of the
$\gamma$-action ($\gamma\in H_1(Y_0;\Z)$) in $\HFinf(Y_0,\spinct_0)$ to
the image of the $\gamma$-action in $\HFinf(Y_0\#
Z,\spinct_0\#\spincu)$. It follows from naturality that 
$$\Finf{Y_0\# Z,\spinct_0\#\spincu}(\cdot\otimes \Theta_Z)
\colon \HFinf(Y_0,\spinct_0)\longrightarrow \HFinf(Y_0,\spinct_0\#\spincu)$$ 
is an isomorphism in all degrees. Commutative
Diagram~\eqref{eq:MapOfLongExacts} then forces $\Finf{Y_1\#
Z}(\cdot\otimes\Theta_Z)$ to be an isomorphism, as well.  (The case of
$(-1)$-surgery follows by repeating the above discussion using the
$(-1)$-surgery exact sequence.)

The case of rational homology spheres follows by induction on the rank
of the first homology group, and a comparison of long exact sequences,
parallel to the proof of Theorem~\ref{HolDiskTwo:thm:HFinfGen} of~\cite{HolDiskTwo}
in the case of
rational homology spheres.
\end{proof}

\vskip.3cm
\noindent{\bf{Proof of Theorem~\ref{thm:AdditivityOfD}.}}
Fix $\Theta_Z\in\HFleq(Z,\spincu)$ of degree $d(Z,\spincu)$.
Then, we have the following commutative diagram
$$
\begin{CD}
\HFinf_{k}(Y,\spinct)
@>\Finf{Y\# Z,\spinct\#\spincu}(\cdot\otimes \Theta_Z)>> 
\HFinf_{k+d(Z,\spincu)}(Y\# Z,\spinct\#\spincu) \\
@VVV @VVV \\
\HFp_{k}(Y,\spinct)
@>\Fp{Y\# Z,\spinct\#\spincu}(\cdot\otimes \Theta_Z)>> 
\HFp_{k+d(Z,\spincu)}(Y\# Z,\spinct\#\spincu).
\end{CD}
$$
From this it follows easily that
$$d(Y,\spinct)+d(Z,\spincu)\leq d(Y\# Z,\spinct\#\spincu).$$
By the same reasoning, we have the inequality
$$d(-Y,\spinct)+d(-Z,\spincu)\leq d(-Y\# Z,\spinct\#\spincu).$$
In view of Equation~\eqref{eq:DFlipOrientation}, we can conclude that 
$$d(Y,\spinct)+d(Z,\spincu)= d(Y\# Z,\spinct\#\spincu),$$
as claimed.
\qed

\subsection{Correction terms for lens spaces}

We give an inductive formula for the correction terms of lens spaces.

Let $p$ and $q$ be a pair of relatively prime, positive integers.
The lens space $-L(p,q)$ can be given a Heegaard diagram
$(E,\alpha,\gamma,z)$, where $E$ is an oriented two-manifold with
genus $g=1$. If we let $\alpha$ be the ``horizontal'' circle $S^1\times \{0\}$
and $\beta$ be the vertical one $\{0\}\times S^1$, then $\gamma$ is a smoothly
embedded curve homologous to 
$-q \alpha + p \beta$ (which we can take to be a ``straight'' circle
$$\theta\mapsto (e^{\frac{-2\pi i (\theta+\theta_0)}{q}},
e^{\frac{2\pi i (\theta+\theta_1)}{p}}).$$

\begin{figure}
\mbox{\vbox{\epsfbox{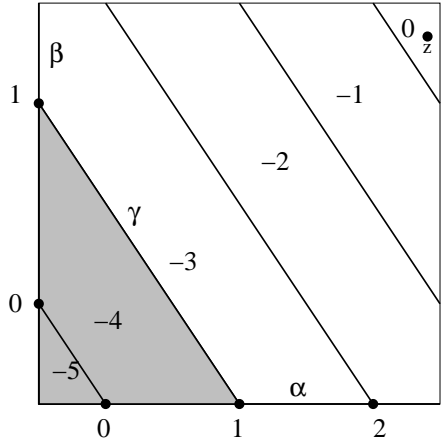}}}
\caption{\label{fig:Lens32}
{\bf{Ordering of $\SpinC$ structures over $-L(3,2)$.}}  This is an
illustration of a Heegaard triple $(E,\alpha,\beta,\gamma,z)$,
representing a cobordism between $-L(3,2)$, $-L(2,1)$, and $S^3$.  The
non-positive integers represent multiplicities of the triply-periodic
domain used in the proof of Proposition~\ref{prop:dLens}, while the
integers labelling the intersection points between $\beta\cap\gamma$
and $\alpha\cap\gamma$ represent the canonical orderings of the
$\SpinC$ structures over $-L(2,1)$ and $-L(3,2)$ respectively. The
darkly-shaded triangle (when taken with multiplicity one) is the
domain associated to the triangle $\psi_1$ in the proof of
Proposition~\ref{prop:dLens}.}
\end{figure}

There is a canonical circular ordering of the $\SpinC$ structures over
$-L(p,q)$ (i.e. a labeling of the $\SpinC$ structures by elements
$i\in\Zmod{p}$), which we can describe as follows. Consider the
pointed Heegaard triple $(E,\alpha,\beta,\gamma,z)$, where the
baspoint is placed so that all the coefficients of the triply-periodic
domain connecting $\alpha$, $\beta$, and $\gamma$ are negative, and
order the intersection points of $\alpha$ with $\gamma$ circularly
(about $\alpha$), so that the $(p-1)^{st}$ one
modulo $p$ is the one adjacent to the
basepoint. (See Figure~\ref{fig:Lens32}.)

\begin{prop}
\label{prop:dLens}
Fix positive, relatively prime integers $p>q$, and also choose
an integer with $0\leq i <p+q$. Then, 
with respect to the above ordering of the $\SpinC$ structures over
$-L(p,q)$, we have the following inductive formula:
\begin{equation}
\label{eq:dLens}
d(-L(p,q),i)=\left(\frac{pq-(2i+1-p-q)^2}{4pq}\right)-d(-L(q,r),j),
\end{equation}
where $r$ and $j$ are the reductions modulo $q$ of $p$ and $i$ respectively.
\end{prop}

\begin{proof}
Observe that the Heegaard triple $(E,\alphas,\betas,\gammas,z)$
represents a cobordism $X_{\alpha,\beta,\gamma}$ with $$\partial
X_{\alpha,\beta,\gamma}=-Y_{\alpha,\beta}-Y_{\beta,\gamma}+Y_{\alpha,\gamma}
= S^3 - Y_{\beta,\gamma}-L(p,q).$$ Indeed, it is easy to see that
$Y_{\beta,\gamma}\cong -L(q,r)$ (and in fact the cobordism obtained
from $X_{\alpha,\beta,\gamma}$ after filling in the $S^3$ is the usual
cobordism between $L(q,r)$ and $L(p,q)$ given by a single two-handle
addition).

Observe that there are $p+q$ triangles $\{\psi_0,...,\psi_{p+q-1}\}$
with non-negative domains but for which $n_z(\psi_0)=0$. We order
these so that $\cald(\psi_0)<\cald(\psi_1)<...<\cald(\psi_{p+q-1})$.
Clearly, with the above conventions, $\psi_i$ connects the $i^{th}$
$\SpinC$ structure (where $i$ is thought of as an element of
$\Zmod{q}$) on $-L(q,r)$ to the $i^{th}$ $\SpinC$ structure (now
regarding $i$ as an element of $\Zmod{p}$) on $-L(p,q)$.  Clearly,
now, the induced map $\Fleq{\alpha,\beta,\gamma}$ determined by
$\psi_i$ takes the tensor product of generators in $\HFleq(S^3)\otimes
\HFleq(-L(q,r),j)$ isomorphically to the generator of
$\HFleq(-L(p,q),i)$ (since there is a unique holomorphic triangle in
the torus representing the $\SpinC$ structure). Thus, as in
Equation~\eqref{eq:DimensionShiftFormula}, it follows that
\begin{equation}
\label{eq:InductionOfD}
d(-L(p,q),i)+d(-L(q,r),j)=\frac{c_1(\spinc_z(\psi_i))^2+1}{4}.
\end{equation}
Strictly speaking, the holomorphic triangle is not the one
corresponding to the two-handle addition:
the diagram $(E,\{\beta\},\{\gamma\},z)$ represents
$L(q,r)$, rather than a connected sum of $S^2\times S^1$.
Thus, we use the variant of this dimension
shift formula as in Proposition~\ref{prop:ConnSumShift}.

Calculating the right-hand-side of the above equation is an
application of Proposition~\ref{HolDiskFour:prop:COneFormula}
of~\cite{HolDiskFour}, a formula which gives the evaluation of the
first Chern class of the $\SpinC$ structure underlying a triangle in terms
of combinatorial data on the Heegaard triple.
Let $\PerDom$ be the generator of the group of triply-periodic
domains, which is nowhere positive. The quantities appearing in
Equation~\eqref{HolDiskFour:eq:COneFormula} of~\cite{HolDiskFour} are easily seen to be:
\begin{eqnarray*}
\EulerMeasure(\PerDom)&=& 0 \\
\#(\partial\PerDom)&=& p+q+1 \\
n_z(\PerDom)&=& 0 \\
\Spider(\psi_i,\PerDom)&=&-p-q+i+1.
\end{eqnarray*}
It now follows from Equation~\eqref{HolDiskFour:eq:COneFormula} of~\cite{HolDiskFour} that
\begin{equation}
\label{eq:COneLens}
\langle
c_1(\spinc_z(\psi_i)), H(\PerDom)\rangle = 2i+1-p-q.
\end{equation}
Clearly,
$H(\PerDom)$ is the generator of the compactly-supported cohomology,
and $H(\PerDom)^2=-pq$, we can combine the above with
Equation~\eqref{eq:InductionOfD} to obtain Equation~\eqref{eq:dLens}.
\end{proof}

\subsection{Correction terms for three-manifolds with $b_1>0$}
\label{subsec:CorrTermBig}

There are versions of the correction term for three-manifolds whose first
Betti number is positive, as well. 
For instance, one could use totally twisted coefficients to define a
rational invariant for three-manifolds equipped with torsion $\SpinC$
structures. However, more structure can be obtained by using untwisted coefficients
(and still a torsion $\SpinC$ structure).
For simplicity, we restrict ourselves presently to the
case where $H_1(Y_0;\Z)\cong\Z$. In this case, there is a unique
$\SpinC$ structure $\spinct_0$ with $c_1(\spinct_0)=0$.

\begin{defn}
Suppose that $H_1(Y_0;\Z)\cong\Z$. Then,
there are two
correction terms $d_{\pm 1/2}(Y_0)$, where $d_{\pm 1/2}(Y_0)$ is the
minimal grading of any
non-torsion element in the image of $\HFinf(Y_0,\spinct_0)$ in
$\HFp(Y,\spinct_0)$ with grading $\pm 1/2$ modulo $2$.
\end{defn}

\begin{prop}
\label{prop:CorrTermOrientBOne}
Let $H_1(Y_0;\Z)\cong \Z$. Then,
$$d_{1/2}(Y_0)-1\leq d_{-1/2}(Y_0).$$
Moreover,
$$d_{\pm 1/2}(Y_0,\spinct)=d_{\pm 1/2}(Y_0,{\overline\spinct}),$$
and
\begin{equation}
\label{eq:DFlipOrientationBOne}
d_{\pm 1/2}(Y_0,\spinct)=-d_{\mp 1/2}(-Y_0,\spinct).
\end{equation}
\end{prop}

\begin{proof}
The first inquality follows from the algebra structure of
$\HFinf(Y_0,\spinct_0)$, together with the grading as given in 
Proposition~\ref{prop:ZModTwoGradingBOne}: if $\xi_0$ is any non-zero element with
$\liftGr(\xi_0)\equiv-\OneHalf\pmod{\Z}$ in $\HFp(Y_0,\spinct_0)$ coming from
$\HFinf(Y_0,\spinct_0)$, then there must be an element $\xi_1$ of
degree one greater (also coming from $\HFinf(Y_0,\spinct_0)$), with
$\theta\cm \xi_1=\xi_0$, where $\theta$ is some element of
$H_1(Y_0,\spinct_0)$.

The other two equations follow exactly as in the case where $b_1(Y)=0$ 
(Proposition~\ref{prop:CorrTermOrient}).
\end{proof}

The following proposition gives a relationship between correction
terms for for integral homology spheres and the correction terms for
zero-surgeries on knots in them. This result will be generalized in
Section~\ref{sec:DefiniteForms} (see especially
Theorem~\ref{thm:IntFormBOneOne} and Corollary~\ref{cor:FracSurg}).

\begin{prop}
\label{prop:CorrTermBound}
Let $K\subset Y$ be a knot in an integral homology three-sphere. Then,
$$ d(Y)-\OneHalf \leq d_{-1/2}(Y_0),$$
and 
$$ d_{+1/2}(Y_0)-\OneHalf \leq d(Y_1).$$
\end{prop}

\begin{proof}
The first inequality follows from the fact that for the induced map $F_1$ on the cobordism from $Y$ to $Y_0$
$$F_1^\infty \colon \HFinf_k(Y)\longrightarrow \HFinf_{k-1/2}(Y_0,\spinct_0)$$ 
is an isomorphism for all even integers $k$,
together with commutativity of the diagram
$$ \begin{CD}
\HFinf_{k}(Y) @>{F_1^\infty}>> \HFinf_{k-1/2}(Y_0,\spinct_0) \\
@VVV @VVV \\
\HFp_{k}(Y) @>{F_1^+}>> \HFp_{k-1/2}(Y_0,\spinct_0) \\
\end{CD}
$$

The second follows similarly, using the fact that
$$F_2^\infty \colon \HFinf_{k+1/2}(Y_0) \longrightarrow \HFinf_{k}(Y_1)	$$
is an isomorphism for all even integers $k$.
\end{proof}

\subsection{Knots in $S^3$}

In light of Proposition~\ref{prop:CorrTermBound}, the correction terms
for homology $S^1\times S^2$ can be used to give obstructions to
obtaining a given three-manifold as zero-surgery on a knot in the
three-sphere. Specifically, since $d(S^3)=0$, we see that if $Y_0$ is
obtained as zero-surgery on a knot in $S^3$, then $$-\OneHalf\leq
d_{-1/2}(Y_0).$$ Moreover, by reflecting the knot and using
Proposition~\ref{prop:CorrTermOrientBOne}, we also obtain the bound
$$d_{1/2}(Y_0) \leq \OneHalf.$$ We will give a generalization of this
observation in Theorem~\ref{thm:IntFormBOneOne} (see especially
Corollary~\ref{cor:NotKnot}).

In another direction, we can think of the correction terms as giving rise to
knot invariants (by considering the two correction terms associated to
the zero-surgery). In fact, we find convenient to package the information as follows:
$$\sigma_+(K)=\frac{d_{-1/2}(S^3_{0}(K))-d_{1/2}(S^3_{0}(K))+1}{2}$$ and
$$\sigma_-(K)=\frac{d_{-1/2}(S^3_{0}(K))+d_{1/2}(S^3_{0}(K))}{2}$$ (where
$S^3_0(K)$ denotes zero-surgery along $K$, given the orientation induced
from $K$, and the correction term is calculated in the unique torsion
$\SpinC$ structure).  

In view of the above results, $\sigma_+(K)$ is a non-negative integer,
while $\sigma_-(K)$ is an integer. Moreover, if
$r\colon S^3\longrightarrow S^3$ is an
orientation-reversing
diffeomorphism of $S^3$ to itself, then
$\sigma_+(r(K))=\sigma_+(K)$, while $\sigma_-(r(K))=-\sigma_-(K)$.

Indeed, the $\sigma_\pm(K)$ could alternatively be defined using only correction terms for 
integral homology three-spheres, in view of the following result:

\begin{prop}
\label{prop:CorrTermEquality}
Let $K\subset S^3$ be an oriented knot in the three-sphere. Then, 
\begin{eqnarray*}
d_{1/2}(S^3_{0}(K))-\OneHalf &=& d(S^3_1(K)) \\
d(S^3_{-1}(K))-\OneHalf &=& d_{-1/2}(S^3_{0}(K))
\end{eqnarray*}
\end{prop}

\begin{proof}
This is a direct consequence of the 
surgery long exact sequence, in view of the structure of
$\HFp(S^3)$.
\end{proof}

\section{The renormalized Euler characteristic and the Casson invariant}
\label{sec:Euler}

Let $Y$ be an integer homology three-sphere (so that it has a unique
$\SpinC$ structure, which we drop from the notation).
We define a renormalized Euler characteristic by
$$\renEuler(Y)=\chi(\HFpRed(Y))-\frac{1}{2}d(Y).$$

Let $\Sigma(2,3,5)$ denote the Poincar\'e homology sphere, oriented as
the boundary of the negative-definite $E8$ plumbing.

\begin{theorem}
\label{thm:EulerCasson}
Let $Y$ be an integer homology three-sphere. Then, 
the renormalized Euler characteristic agrees with Casson's invariant:
$$
\chiRen(Y)=\lambda(Y),
$$
where here  Casson's invariant is normalized so that
$\lambda(\Sigma(2,3,5))=-1$.
\end{theorem}

We will use the surgery long exact sequence for $\HFp$,
and the following observation.

\begin{lemma}
\label{lemma:ChiTrunc}
Let $Y$ be an integral homology three-sphere. Then,
for all sufficiently large $n$, we have that
$$\chiRen(Y)=\chi(\HFp_{\leq 2n-1}(Y))-n.$$
\end{lemma}

\begin{proof}
This follows easily from the structure of $\HFinf(Y)$.
\end{proof}

Recall that when $Y_0$ is an integer homology $S^1\times
S^2$, the Euler characteristic $\chi(\HFp_{\leq 2n+1}(Y,\spinc_0))$
is independent of $n$, provided that $n$ is sufficiently large. We let
$\chiTrunc(\HFp(Y_0,\spinc_0))$ denote this integer.

\begin{prop}
If $K\subset Y$ is a knot in an integral homology three-sphere, and
let $\spinc_0$ denote the $\SpinC$ structure on $Y_0$ with trivial $c_1$. 
Then we
have the following:
\begin{equation}
\label{eq:SurgeryFormula}
\renEuler(Y)-\renEuler(Y_1) =
\chiTrunc(Y_0,\spinc_0)+
\sum_{\spinc\neq\spinc_0}\chi(\HFp(Y_0,\spinc)).
\end{equation}
\end{prop}

\begin{proof}
This follows by taking the Euler characteristic of 
the exact sequence in the form stated in 
Theorem~\ref{thm:TruncExact}, and applying the observation of 
Lemma~\ref{lemma:ChiTrunc}.
\end{proof}

\noindent{\bf{Proof of Theorem~\ref{thm:EulerCasson}.}}
The right-hand-side of Equation~\eqref{eq:SurgeryFormula} can be
calculated from the results of~\cite{HolDiskTwo}. Specifically if we
write the Alexander polynomial of $Y-K$ as $$\Delta_K=a_0+\sum_{i=1}^d
a_i (T^i+T^{-i}),$$ and $$\MT_i=\sum_{j=1}^d ja_{|i|+j},$$ then
according to Theorem~\ref{HolDiskTwo:thm:EulerOne}
of~\cite{HolDiskTwo} (with
Proposition~\ref{HolDiskTwo:prop:PreciseChi} of~\cite{HolDiskTwo} to
pin down the sign, bearing in mind that for an integer homology
$S^1\times S^2$, the mod $2$ reduction of $\liftGr+\OneHalf$ gives the
absolute $\Zmod{2}$ grading on $Y_0$ which is used to determine the sign
of the Euler characteristic, according to Proposition~\ref{prop:ZModTwoGradingBOne}),
for each $i\neq 0$, we have that $$\chi(\HFp(Y_0,\spinc))=-\MT_i,$$
where $c_1(\spinc)$ is $2i$ times a generator of $H^2(Y_0;\Z)$, while
$$\chiTrunc(\HFp(Y_0,\spinc_0))=-\MT_0,$$ for the torsion $\SpinC$ structure
(this case is handled in Theorem~\ref{HolDiskTwo:thm:TruncEuler}).
Plugging these values into Equation~\eqref{eq:SurgeryFormula}, we see
that $\renEuler(Y,\spinc)$ satisfies the same surgery formula as
Casson's invariant.  Moreover, we know that $\renEuler(S^3)=0$. Since
Casson's invariant is characterized by its $+1$ surgery formula and
this normalization, the theorem follows.
\qed

\section{The renormalized complexity and surgeries}
\label{sec:Complexity}

In~\cite{HolDiskTwo}, we defined a numerical invariant for
integer homology three-spheres $$N(Y)=\Rk \HFred(Y).$$ Let $K\subset
Y$ be a knot, and let $Y_1$ be the manifold obtained by $+1$ surgery
along $K$, then Theorem~\ref{HolDiskTwo:thm:Complexity}
of~\cite{HolDiskTwo} gives a bound: 
$$ \max(-\MT_0,0)+2 \sum_{i=1}^{d}|\MT_i(K)| \leq N(Y)+N(Y_1). 
$$ 
Taken with the correction term, however, $N(Y)$ becomes more
effective at distinguishing $Y$ and $Y_1$, as follows:

\begin{theorem}
\label{thm:RenormalizedComplexity} 
Let $Y$ be an integral homology three-sphere and $K\subset Y$ be a knot, then there is a bound:
\begin{equation}
\label{eq:ComplexityBound}
|\MT_0(Y)|+
2 \sum_{i=1}^d |\MT_i(K)| \leq N(Y)+\frac{d(Y)}{2} + N(Y_1)-\frac{d(Y_1)}{2}.
\end{equation}
\end{theorem}

\begin{proof}
By surgery long exact sequence, we have exactness in the middle for
$$
\begin{CD}
\HFred(Y) @>{\Fred{1}}>> \HFred(Y_0) @>{\Fred{2}}>> \HFred(Y_1).
\end{CD}
$$
Indeed, the image of ${\Fred{2}}$ lies in the kernel of the map 
$\Fred{3}\colon\HFred(Y_1)\longrightarrow \HFred(Y)$, 
while the kernel of $\Fred{1}$ contains the image of $\Fred{3}$. 
Thus,
\begin{equation}
\label{eq:BoundHFredZero}
\rk \HFred(Y_0) = \rk \Image \Fred{1} +
\rk \Image \Fred{2}
\end{equation}

We claim that 
\begin{equation}
\label{ineq:KerFThree}
\rk \Image \Fred{2} \leq \rk \HFred(Y_1)-
\left(\frac{d_{-1/2}(Y_0)+\OneHalf-d(Y)}{2}\right).
\end{equation} This is true because
(thanks to the absolute gradings in the exact sequence,
Lemma~\ref{lemma:DegShiftExact}) there is a module $V_1$ of rank
$D_Y=\frac{d_{-1/2}(Y_0)+\OneHalf-d(Y)}{2}$ in the image of
$\HFinfty(Y)$ in $\HFp(Y)$ which maps to zero in $\HFp(Y_0)$. It follows
(from the surgery long exact sequence on $\HFp$)
that $\HFp(Y_1)$ surjects onto $V_1\subset \HFp(Y)$. In fact, the map
$$\Fp{3} \colon \HFp(Y_1) \longrightarrow \HFp(Y)$$ factors as
$$
\begin{CD}
\HFp(Y_1) @>>>\HFred(Y_1) @>>>\HFp(Y),
\end{CD}
$$
since the map $\Finf{3}$ is trivial. 
Moreover, the image of $\Fred{2}$ is clearly contained in the kernel of
this map. But there is a module $W_1$
of rank at least $D_Y$ in $\HFred(Y_1)$ which maps to $V_1$. 
This establishes Inequality~\eqref{ineq:KerFThree}.

Similarly, we claim that 
\begin{equation}
\label{ineq:ImageFOne}
\rk \Image \Fred{1} \leq
\rk \HFred(Y)-\left(\frac{d(Y_1)+\OneHalf-d_{1/2}(Y_0)}{2}\right).
\end{equation} There is a module
$V_2$ or rank 
$D_{Y_1}=\frac{d(Y_1)+\OneHalf-d_{1/2}(Y_0)}{2}$ in the
image of $\HFinf(Y_0)$ in $\HFp(Y_0)$ which maps to zero in $Y_1$.  It
follows (from the surgery 
long exact sequence on $\HFp$) that there is a submodule of $\HFp(Y)$ which surjects onto
$V_2$. The grading of $V_2$ is congruent to $1/2$ modulo $2$ (see
Proposition~\ref{prop:ZModTwoGradingBOne}, so
clearly, $W_2$ consists of elements with grading $1$ modulo $2$ in
$\HFp(Y)$, so $W_2$ injects into $\HFred(Y)$. On the other hand,
$W_2$ maps to zero under $\Fred{1}$. This establishes 
Inequality~\eqref{ineq:ImageFOne}.

We claim that for all sufficiently large $n$,:
$$\rk \HFred(Y_0,\spinc_0) \leq 
\rk \HFp_{\leq 2n+1}(Y_0,\spinc_0)
-\left(\frac{d_{-1/2}(Y_0)-d_{1/2}(Y_0)+1}{2}\right).$$

Combining the Euler characteristic calculations on $Y_0$
(Theorem~\ref{HolDiskTwo:thm:EulerOne} of~\cite{HolDiskTwo}
when the $\SpinC$ structure is non-torsion and 
Theorem~\ref{HolDiskTwo:thm:TruncEuler} of~\cite{HolDiskTwo}
when it is), Equation~\eqref{eq:BoundHFredZero} 
and Inequalities~\eqref{ineq:KerFThree},
and \eqref{ineq:ImageFOne}, we obtain the result claimed.
\end{proof}

It is natural to consider the following class of three-manifolds:

\begin{defn}
\label{def:Invisible}
An integer homology three-sphere $Y$ is said to be 
{\em invisible} if $d(Y)=0$ and $N(Y)=0$. 
\end{defn}

Invisibility is independent of the orientation of $Y$, according to the following:

\begin{prop}
Let $Y$ be an integral homology three-sphere. Then,
$N(-Y)=N(Y)$ and $d(-Y)=-d(Y)$.
\end{prop}

\begin{proof}
The fact that $N(Y)=N(-Y)$ follows from duality between
$\HFp(Y)$ and $\HF_-(-Y)$ (the latter being cohomology)
(see Proposition~\ref{HolDiskOne:prop:Duality} of~\cite{HolDisk}
and also
Section~\ref{HolDiskFour:sec:Duality} of~\cite{HolDiskFour}).
The claim for $d$ was established in
Proposition~\ref{prop:CorrTermOrient}.
\end{proof}

Of course, $S^3$ is invisible. By the additivity of $d$ under connected
sum (Theorem~\ref{thm:AdditivityOfD}) and the following lemma,
the set of invisible three-manifolds is closed under connected sum:

\begin{lemma}
\label{lemma:ConnSumTrivN}
If $N(Y_1)=0$ and $N(Y_2)=0$, then $N(Y_1\# Y_2)=0$, as well. 
\end{lemma}

\begin{proof}
From the long exact sequence connecting $\HFa$ with $\HFp$, it follows
that $N(Y)=0$ if and only if $\HFa(Y)$ has rank one.  The connected sum claim
then follows from the K\"unneth formula for connected sums on $\HFa$
(Proposition~\ref{HolDiskOne:prop:ConnSum} of~\cite{HolDisk}). 
\end{proof}

As we shall see in Section~\ref{sec:SampleCalculations},
$N(\Sigma(2,3,5))=0$. So, it follows from
Lemma~\ref{lemma:ConnSumTrivN} that $\Sigma(2,3,5)\#-\Sigma(2,3,5)$ is
an invisible three-manifold.  Theorem~\ref{thm:RenormalizedComplexity}
has the following consequence for invisible three-manifolds:

\begin{cor}
Let $Y$ be an invisible three-manifold, and $K\subset Y$ be a knot in
$Y$ with non-trivial Alexander polynomial. Then $Y_1$ is not invisible.
\end{cor}

\begin{proof}
The Alexander polynomial of a knot is non-trivial 
if and only if the left-hand-side of Inequality~\eqref{eq:ComplexityBound} is positive; 
so the result follows from the inequlity.
\end{proof}

Using the long exact sequence for $1/n$ surgeries, we have the
following generalization of Theorem~\ref{thm:RenormalizedComplexity}:

\begin{theorem}
\label{thm:RenormalizedComplexityFrac} 
Let $Y$ be an integral homology three-sphere and $K\subset Y$ be a knot, then there is a bound
$$ n\left(|\MT_0(Y)|+
2 \sum_{i=1}^d |\MT_i(K)|\right) \leq N(Y)+\frac{d(Y)}{2} + N(Y_{1/n})-\frac{d(Y_{1/n})}{2}.$$
\end{theorem}

\begin{proof}
The proof proceeds exactly as before, substituting the long exact sequence for $1/n$ surgeries, 
with the additional observation that in each degree $k$,
$\HFinf_k(Y_0,\spinc_0;\Zmod{n})\cong \Z$. Note that we are using Proposition~\ref{prop:DegShiftExactFrac}
in place of Lemma~\ref{lemma:DegShiftExact}.
\end{proof}

This has the following immediate consequence:

\begin{cor}
Let $Y$ be an invisible three-manifold, and $K\subset Y$ be a knot in
$Y$ with non-trivial Alexander polynomial. Then any non-trivial surgery on $K$ gives
a three-manifold which is not invisible.
\end{cor}

\section{Integer surgeries in the three-sphere}
\label{sec:Lens}

We will now consider consequences of the graded exact sequences
to the situation where $K\subset S^3$ is a knot with the property
that $+p$ surgery on $K$ gives a lens space. 

Indeed, as we shall show, when $K\subset S^3$ is a knot for which
$S^3_p(K)\cong L(p,q)$, then, the absolute gradings together with the
long exact sequence for integer surgeries
(Theorem~\ref{HolDiskTwo:thm:ExactP} of~\cite{HolDiskTwo}) determine
the structure of $\HFp(S^3_{0}(K))$. The methods here can be thought
of as elaborations on the proof of
Theorem~\ref{thm:RenormalizedComplexity}. Note that one need consider
only integer surgeries on the knot $K$, since if a non-integral
surgery on a knot in $S^3$ gives rise to a lens space, then the knot
must be a torus knot, according to the ``cyclic surgery theorem'' of
Culler-Gordon-Luecke-Shalen,~\cite{CyclicSurgery}.  

In fact, most of the results we discuss here for integer surgeries
readily generalize to the case where $S^3$ is replaced by an arbitrary
invisible three-manifold (in the sense of
Definition~\ref{def:Invisible}), and the lens space is replaced by an
arbitrary three-manifold $L$ with $\HFpRed(L)=0$ and $H_1(L,\Z)\cong
\Zmod{p}$. However, we state most of our results for $S^3$ and lens
spaces, in the interest of exposition.

We recall now the integral surgeries long exact sequence.
Let $K\subset Y$ be a knot in an integral homology three-sphere and
$p$ is a positive integer, then Theorem~\ref{HolDiskTwo:thm:ExactP}
of~\cite{HolDiskTwo} gives 
a map
$$Q\colon
\SpinC(Y_0)\longrightarrow\SpinC(Y_p)$$
and a long exact sequence of the form:
\begin{equation}
\label{seq:IntegerSurgeries}
\begin{CD}
... @>{F_1}>>  \HFp(Y_0,[\spinct])@>{F_2}>> \HFp(Y_{p},\spinct)@>{F_3}>> \HFp(Y)
@>>> ...
\end{CD},
\end{equation}
where 
$$\HFp(Y_0,[\spinct])=\bigoplus_{\spinct'\in Q^{-1}(\spinct)} \HFp(Y_0,\spinct').$$

To describe $Q$ topologically, recall that the integer surgery
naturally gives rise to a cobordism $W$ between $Y_0$, $Y_p$, and the
lens space $L(p,1)$. The proof of the long exact sequence (see
especially Proposition~\ref{HolDiskTwo:prop:OneQHoClassesCancel}
of~\cite{HolDiskTwo}), gives a preferred $\SpinC$ structure
$\prefspinc$ over $L(p,1)$.
If $\spinc$ is any $\SpinC$ structure over
$Y_0$, then $Q(\spinct)$ is the $\SpinC$ structure over $Y_p$ for which
the triple of $\SpinC$ structures $\spinc$, $Q(\spinct)$ and
$\prefspinc$ extend over $W$.  Indeed, the preferred $\SpinC$ structure
over $L(p,1)$ is characterized by the following:

\begin{prop}
\label{prop:SpinCTZero}
Let $N$ be a neighborhood of a two-sphere $S$ with self-intersection
number $-p$. Then, $\prefspinc$ is the $\SpinC$ structure which extends
to a $\SpinC$ structure $\spinc$ over $N$ with 
$$\langle c_1(\spinc), [S]\rangle = p. $$
\end{prop}

\begin{proof}
The intersection point representing $\prefspinc$ is adjacent to the
basepoint (compare Figure~\ref{HolDiskTwo:fig:OneQSurgery}
of~\cite{HolDiskTwo}). The rest is an application of Proposition~\ref{HolDiskFour:prop:COneFormula} of~\cite{HolDiskFour}
(as in Proposition~\ref{prop:dLens}).
\end{proof}

There are purely algebraic constraints on realizing a given map from
$\SpinC(Y_0)$ to $\SpinC(Y)$ as a map of the type $Q$ above.  Of
course, the two spaces of $\SpinC$ structures are principal
homogeneous spaces for $H^2(Y_0)\cong \Z$ and $H^2(Y_p)\cong \Zmod{p}$
respectively, and the map $Q$ must be equivariant under this action
(given a surjective group homomorphism from $H^2(Y_0)\longrightarrow
H^2(Y_p)$). In addition, both spaces admit actions by $\Zmod{2}$,
given by conjugating the $\SpinC$ structures, and the map $Q$ must
also be equivariant under these $\Zmod{2}$-actions, as well.

In the following statement, recall (c.f. Proposition~\ref{prop:dLens})
that for each positive integer $p$ and each congruence class $i\in\Zmod{p}$
$$d(L(p,1),i)=\frac{(2j-p)^2-p}{4p},$$
where we take $j$ to be the integer in the equivalence class $j\equiv i\pmod{p}$ with
$0\leq j <p$. Note that in that proposition we gave an explicit identification
$$\SpinC(L(p,q))\cong \Zmod{p}.$$

When describing $\SpinC$ structures over the zero-surgery $Y_0$, we will
find it convenient to use an identification
$$\SpinC(Y_0)\cong \Z$$
induced from a choice of generator $H$ for $H_2(Y_0,\Z)$. In particular,
we write
$\HFp(Y_0,i)$ to denote the group for $Y_0$ associated to the $\SpinC$ structure
$\spinct_i\in\SpinC(Y_0)$ with the property that
$$\langle c_1(\spinc_i),[H]\rangle =2i$$
(note that the group $\HFp(Y_0,i)$ is actually independent of the choice
of generator $H$, since the groups are invariant under conjugation).

\begin{theorem}
\label{thm:PSurgeryLens}
Let $K\subset S^3$ be a knot in $S^3$ with the property that $p>0$
surgery on $K$ gives the lens space $L(p,q)$, and  let $Y_0=S^3_0(K)$.
Then, $\HFp(Y_0)$ has the
following structure:
\begin{itemize}

\item The group  $\HFinf(Y_0,0)$ 
surjects onto $\HFp(Y_0,0)$, and $\HFp(Y_0,0)$
contains no torsion.  Thus, $\HFp(Y_0,0)$ is determined by
$d_{\pm 1/2}(Y_0,0)$. In fact,
$$d_{-1/2}(Y_0,0)=-\OneHalf,$$ and
$$d_{1/2}(Y_0,0)=d(L(p,q),Q(0))-d(L(p,1),0))+\OneHalf.$$

\item For each $i$ with $|i|\leq p/2$, 
all non-zero homogeneous elements in $\HFp(Y_0,i)$  have odd grading, 
and in fact we 
have an isomorphism of $\Z[U]$-modules
$$\HFp(Y_0,i)\cong \Z[U]/U^{\ell},$$
where  the integer $\ell=\ell(p,q,Q,i)$
is given by the formula:
\begin{equation}
\label{eq:EqForL}
2\ell=-d(L(p,q),Q(i))+d(L(p,1),i)\geq 0,
\end{equation}

\item For each $i$ with $|i|>p/2$, $\HFp(Y_0,i)=0$.
\end{itemize} 
\end{theorem}

\begin{remark}
Of course, no generality is lost by focusing on the case of $+p$
surgeries. If $-p$ surgery on a knot $K$ gives the lens space
$L(p,q)$, then we can apply the above theorem to the reflection of
$K$, $r(K)$, bearing in mind that $+p$ surgery on $r(K)$ gives the
lens space $-L(p,q)=L(p,p-q)$, and also that $S^3_0(r(K))=-S^3_0(K)$.
\end{remark}

\begin{remark} The methods of this section actually prove a stronger statement:
if $L$ is any three-manifold with $\HFpRed(L)=0$ and $K\subset Y$ is a
knot in an invisible three-manifold with the property that
$Y_p(K)\cong L$, then $\HFp(Y_0)$ is uniquely specified by formulas
similar to those appearing in the statement of
Theorem~\ref{thm:PSurgeryLens}, depending on the correction terms for
$Y$ and $L$ (and the correspondence $Q$). We do not spell these out at
present, since the case of lens space surgeries seems to be the most natural.
\end{remark}

Before turning to the proof, we give some consequences of the above
theorem.

\begin{cor}
\label{cor:AlexLens}
For each pair of relatively prime integers $(p,q)$, there is a finite
set of symmetric Laurent polynomials in a variable $T$ (explicitly determined
by $p$ and $q$) which can arise as the Alexander polynomial of a 
knot $K\subset S^3$ with the property that $S^3_p(K)\cong L(p,q)$.

More precisely, if 
$K\subset S^3$ is a knot with $S^3_p(K)\cong L(p,q)$, then there is a
one-to-one correspondence $$\sigma \colon
\Zmod{p}\longrightarrow \SpinC(L(p,q))$$ 
with the property for each integer $i$, we have that
$$0\leq 2\MT_i(K)=
\left\{\begin{array}{ll}
-d(L(p,q),\sigma(i)) + d(L(p,1),i) &{\text{if $2|i|\leq p$}} \\
0 & {\text{otherwise,}}
\end{array}
\right. $$
and the correspondence $\sigma$ satisfies the following symmetries:
\begin{itemize}
\item $\sigma(-i)={\overline{\sigma(i)}}$
\item there is an isomorphism $\phi\colon \Zmod{p}
\longrightarrow \Zmod{p}$
with the property that 
$$\sigma(i)-\sigma(j)=\phi(i-j).$$
\end{itemize}
\end{cor}

\begin{proof}
If $K$ is as above, the equation for the torsions
$\MT_i(K)$ is an immediate
consequence of Theorem~\ref{thm:PSurgeryLens}, together with the
relationship between the Euler characteristic of $\HFp$ for $Y_0$ and
the torsion invariants for $Y_0$ The symmetry properties of $\sigma$
follow immediately from the corresponding symmetries of $Q$.
Since there are only finitely many different possible choices for 
$Q$ (corresponding to the various choices of $\sigma$), and the Alexander
polynomial of $K$ is uniquely determined by the torsion coefficients
(c.f. Equation~\eqref{eq:defBi}), the first statement in the corollary follows.
\end{proof}

This has the following special cases (stated in the introduction as
Theorems~\ref{intro:SmallLenses}
and \ref{intro:PLens}):

\begin{cor}
\label{cor:SmallLenses}
Suppose that $K\subset S^3$ is a knot with the property that some integer surgery 
along $K$ with coefficient $p$ with $|p|\leq 4$ 
gives a lens space, then 
$\HFp(S^3_{K}(0))\cong \HFp(S^2\times S^1)$ as absolutely graded
groups; in particular,
the Alexander
polynomial of $K$ is trivial. 
\end{cor}

\begin{proof}
This follows from a case-by-case analysis. For each $p$ with $|p|\leq
4$, and each $q\neq 1$, it is easy to see that there is no one-to-one
correspondence $\sigma$ between $\SpinC(L(p,q))$ and $\SpinC(L(p,1))$
for which all differences
$-d(L(p,q),\spinct))+d(L(p,1),\sigma(\spinct))$ are all non-negative,
even integers. Indeed, when $q=1$, the only possible correspondence is
the trivial one (for which all the differences are zero), forcing
$\HFp(S^3_K(0))\cong \HFp(S^2\times S^1)$.

The proof is straightforward, given that the correction terms for
$L(2,1)$ are given by $(-\frac{1}{4},\frac{1}{4})$; the correction
terms for $L(3,1)$ are $(\OneHalf,-\frac{1}{6}, -\frac{1}{6})$, and
those for $L(4,1)$ are $(\frac{3}{4}, 0,-\frac{1}{4},0)$.
\end{proof}

The above result is special to the case where $|p|<5$. For instance,
the lens space $L(5,4)=-L(5,1)$ can be realized as $+5$-surgery on the
right-handed trefoil knot, whose Alexander polynomial is
non-trivial. Indeed, the constraints given above show that any knot
$K$ with the property that $S^3_p(K)=L(5,4)$, the Alexander polynomial
is given by $\Delta_K=T-1+T^{-1}$.  However, we do have the following
statement for general $p$ (see also Section~\ref{subsec:LensAlex} for
a table of possible Alexander polynomials for knots giving $L(p,q)$,
for small values of $p$):

\begin{cor}
\label{cor:AlexP}
Suppose that $K\subset S^3$ is a knot with the property that $+p$ surgery on $K$ gives
the lens space $L(p,1)$. Then, 
$\HFinf(S^3_0(K))\cong\HFinf(S^2\times S^1)$ as absolutely graded
groups; in particular, the Alexander polynomial of $K$ is trivial.
\end{cor}

\begin{proof}
Suppose that $+p$ surgery on $K$ gives the lens space $L(p,1)$.  If
$\HFp(S^3_0(K))\not\cong \HFinf(S^2\times S^1)$, then there must be some
non-zero $\MT_i(K)$, and hence there must be some $i\in\Zmod{p}$ with the
property that $-d(L(p,1),\sigma(i))+d(L(p,1),i)\neq 0$ (according to
Theorem~\ref{thm:PSurgeryLens}). Thus, there must be some
(possibly different) $j\in\Zmod{p}$ so that
$-d(L(p,1),\sigma(j))+d(L(p,1),j)$ is negative. But this contradicts
the non-negativity of the $\MT_i$ from Theorem~\ref{thm:PSurgeryLens}).
\end{proof}

In a different direction, Theorem~\ref{thm:PSurgeryLens} gives
obstructions to realizing a given lens space as integral surgery on a
knot in $S^3$.  For example, the condition of integrality of $\ell$
(as given in Equation~\eqref{eq:EqForL}) could be viewed as an
obstruction to obtaining $L(p,q)$ in this way.
But it is not particularly strong: since it uses the correction term
only modulo $2\Z$, it gives an obstruction only to obtaining $L(p,q)$
on a knot in any homology three-sphere.  Now, there is a complete
characterization of such lens spaces, due to
Fintushel-Stern~\cite{FSLensSurgeries}:

\begin{prop}
The lens space $L(p,q)$ can be obtained as integral 
surgery on a knot in an
integral homology three-sphere if and only if $\pm q$ is a square modulo $p$.
\end{prop}

\begin{proof}
We consider the link obtained by a single unknot with
framing $p/q$, and another knot which links this with linking number
$x$, and which is given the framing $n$. Let $Y$ be the three-manifold
obtained as surgery on this link. Then $|H_1(Y;\Z)|=|np-qx^2|$. Thus, if 
$\pm q$ is a square mod $p$, we can find an $n$ and $x$
such that $Y$ is a homology three-sphere.

Conversely, if $Y$ is a homology three-sphere with a knot on which an
integral surgery gives $L(p,q)$, we can find instead a knot in
$L(p,q)$, on which an integral surgery gives $Y$.  Thus, the above
argument shows that $\pm q$ is a square modulo $p$..
\end{proof}

However, the non-negativity of the $\MT_i$ (coming from the the
absolute $\Zmod{2}$ grading of the $\Z[U]/U^\ell$ as above) gives a
more refined obstruction to realizing a fixed lens space $L(p,q)$ as
surgery on some knot in $S^3$. We content ourselves here with one
infinite family of lens spaces ruled out by this obstruction.

\begin{prop}
\label{prop:FamiliesOfLenses}
Consider the family of lens spaces $L(p,q)$ parameterized by positive
integers $k$ not divisible by four, with $p=2k(3+8k)$, and
$q=2k+1$. These spaces cannot be obtained as integral surgeries on any
knot in $S^3$, though they all arise as integral surgeries on knots in
homology spheres. Moreover, these lens spaces can be obtained by
integral surgeries on two-component links in $S^3$.
\end{prop}

\begin{proof}
We consider the $\SpinC$ structure on $L(2k(3+8k),2k+1)$ labelled by
the integer $k$, according to the ordering given in
Proposition~\ref{prop:dLens}. Since
$$2k(3+8k)\equiv 1\pmod{1+2k},$$
the other lens space appearing in the inductive formula 
is $-L(1+2k,1)$.
Thus after two iterations of Equation~\eqref{eq:dLens}, we get that
$$d(-L(2k(3+8k)k,2k+1),k)=\frac{1-8k}{4}.$$
In the same manner,
$$d(-L(2k(3+8k)k,2k+1),4k(1+2k))=\frac{1+2k}{4}.$$

We claim also that both $\SpinC$ structures in question are spin
structures. To see this, observe that the $\SpinC$ structure labelled
by $i=4k(1+2k)$ extends over the cobordism between $L(p,q)$ and
$L(q,1)$ whose first Chern class is trivial (according to
Equation~\eqref{eq:COneLens}). Now the $\SpinC$ structure labeled $k$
differs from $i$ by $p/2$ times a generator of $H^2(L(p,q);\Z)$ thus
it, too, must come from a $\Spin$ structure.

Thus, the correspondence $Q$ must pair the 
spin correction term $\frac{8k-1}{4}$
for $L(p,q)$ with one of the two possible
correction terms for $L(q,1)$ which come from $\Spin$ structures,
namely, $-1/4$ and $\frac{1-p}{4}$. Pairing with the first is ruled
out by the positivity criterion of Theorem~\ref{thm:PSurgeryLens},
while pairing with the second is ruled out by the integrality
condition of that same theorem, in light of our hypothesis that $k$ is
not divisible by $4$. It follows that $+p$ surgery on a knot
cannot give $L(p,q)$.

To rule out $-L(p,q)$, we observe that the other spin correction term
for this manifold, $\frac{1+2k}{4}$, cannot pair with $\frac{1-p}{4}$
by the integrality criterion (and the hypothesis that $k$ is not
divisible by $4$). On the other hand, pairing it with $-1/4$ is once
again ruled out by the positivity criterion.

On the other hand, these lens spaces all arise as surgeries in
homology spheres.  Specifically, consider the plumbing diagram
consisting of a tree with a central node and three chains of spheres.
The central node has a sphere of square $-1$, the first chain consists
of spheres of square $-2$ and $x$ (to be revealed later), the second
consists of a single sphere labelled with $-8k-1$, and the third
consists of one node with square $-3$, and then a chain of $2k-1$
spheres with self-intersection number $-2$ (see
Figure~\ref{fig:PlumbLens}). When the sphere labelled with $x$ is left
off, the three-manifold described is simply $+1$-surgery on the
(right-handed) torus knot of type $(2,4k+1)$ (as can be seen by
successively blowing down $-1$-spheres in the plumbing diagram). When
$x=0$, the three-manifold described is the lens space $L(p,q)$: the
$x=0$ sphere cancels the $-2$ sphere in the first chain, and we can
then blow down the central $-1$ sphere, to obtain a single chain
$2k+1$ of two-spheres, the first of which is labelled with $-8k$, and
the rest with self-intersection number $-2$.

Indeed, the above procedure applied to the entire plumbing diagram,
i.e. keeping the first unknot (labelled now with $x=0$) allows us to
express the lens spaces in the given family as integral surgeries on
two-component links (one of whose components is the $(2,4k+1)$ torus
knot, with framing $+1$).
 
\begin{figure}
\mbox{\vbox{\epsfbox{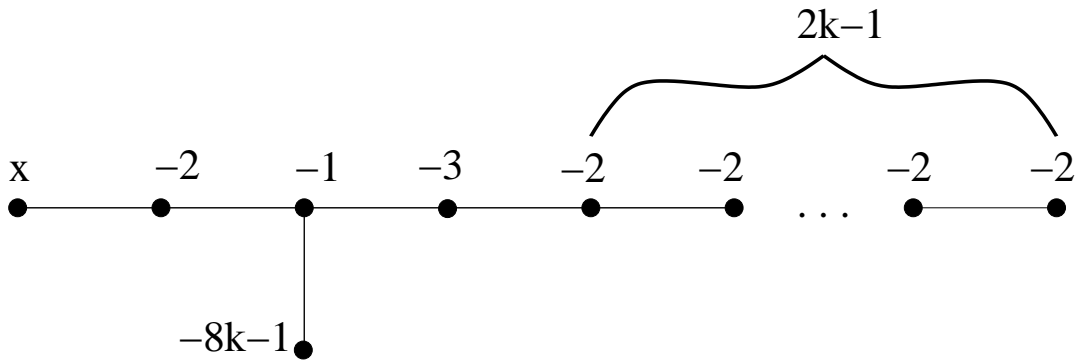}}}
\caption{\label{fig:PlumbLens}
{\bf{Plumbing pictures for $L(2k(3+8k),2k+1)$.}}  This plumbing
picture exhibits $L(2k(3+8k),2k+1)$ as integral surgery
on the homology sphere $-\Sigma(2,4k+1,8k+1)$ (setting $x=0$). 
Blowing down all $(-1)$-spheres, we obtain an description
of the lens space as integral surgery on a two-component link.}
\end{figure}

\end{proof}

\begin{figure}
\label{fig:LensSpaceLink}
\end{figure}

Having seen consequences of Theorem~\ref{thm:PSurgeryLens}, we turn
to its proof, after some lemmas.

For the statement of Theorem~\ref{thm:PSurgeryLens}, it is useful to have an
alternate characterization of $Q$ (at least, up to conjugation).

Fix a knot $K\subset S^3$, and let
$W(Y,K,p)$ denote the cobordism from $S^3$ to $Y_p$. This has $b_2(W)=1$, and indeed, it
has a compactly supported class $\Sigma$ with $\Sigma\cdot\Sigma=p$. 

\begin{lemma}
\label{lemma:IdentifyQ}
Let $\spinct$ be a $\SpinC$ structure over $Y_0$, with
$\langle c_1(\spinct),H\rangle = 2i$,
where $H\in H_2(Y;\Z)$ is a generator.
Then, the $\SpinC$ structure $Q(\spinct)$ extends over $W(Y,K,p)$ as a
$\SpinC$ structure $\spinc$ with
$$\pm \langle c_1(\spinc),[\Sigma]\rangle
\equiv
p+2i\pmod{2p}.$$
\end{lemma}

\begin{proof}
We can juxtapose the cobordisms $W_0$ from $Y$ to $Y_0$ and the
cobordism $W$ from $Y_0\coprod L(p,1)$ to $Y_p$ to obtain a
composite cobordism $X$. By the definition of $Q$, we can find a
$\SpinC$ structure over $W$ whose restriction to $Y_p$ is $Q(\spinct)$
and whose restriction to $L(p,1)$ is the canonical $\prefspinc$.
We can then extend this $\SpinC$ structure from $W$ to obtain
a $\SpinC$ structure $\spinc$ over all of $X$.

Now, the composite cobordism $X$ admits a different
decomposition as the internal connected sum of the canonical cobordism
$W(Y,K,p)$ from $Y$ to $Y_p$ with the null-cobordism of $-L(p,1)$ as a
neighborhood of a sphere $S$ of square $-p$.  (This decomposition is
easily seen by decomposing the Heegaard quadruple
$(\Sigma,\alphas,\betas,\gammas,\deltas,z)$ as a juxtaposition of two
Heegaard triples in two ways.)  Correspondingly, the homology
$H_2(X;\Z)$ is generated by $[\Sigma]$ and $[S]$.  The image of the
generator of $H$ in $X$ is represented by $\pm\PD(\Sigma)-\PD(S)$.
It follows that
$$\langle c_1(\spinc)|Y_0, H\rangle 
= \pm\langle c_1(\spinc), [\Sigma]\rangle
-\langle c_1(\spinc),[S]\rangle 
\equiv \pm \langle c_1(\spinc),[\Sigma]\rangle -p\pmod{2p}.$$
\end{proof}

We also have the following result concerning degree shifts.

\begin{lemma}
\label{lemma:CalcDegrees}
The component of $F_1$ in Exact Sequence~\eqref{seq:IntegerSurgeries} above
which carries $\HFp(Y)$ into the $\spinct_0$-component of
$\HFp(Y_0,[Q(\spinct_0)])$ has degree $-1/2$, while the restriction of
$F_2$ to the $\HFp(Y_0,\spinct_0)$-summand of $\HFp(Y_0,[Q(\spinct_0)])$
has degree $\left(\frac{p-3}{4}\right)$.
\end{lemma}

\begin{proof}
Consider the proof of the integral surgeries long exact sequence
from~\cite{HolDiskTwo}, only now isotope the $\gamma$- rather than the
$\delta$-curves. In so doing, we can realize $F_1$ and $F_2$ as maps
defined by counting holomorphic triangles (where now $F_2$ belongs to
a cobordism between $Y_0$, $Y_p$ and $L(p,1)\#(\#^{g-1}(S^2\times
S^1))$). In particular, $F_1$ is a sum of maps induced by a single
two-handle addition, so it decreases grading by $1/2$.

The map $F_2$ counts pseudo-holomorphic triangles, induced from
$$\HFp(Y_0,\spinct_0)
\otimes \HFleq(-L(p,1)\#^{g-1}(S^1\times S^1),\prefspinc)
\longrightarrow \HFp(Y_p,\spinct),$$
by substituting in the canonical generator of
$\HFleq(-L(p,1)\#^{g-1}(S^1\times S^1),\prefspinc)$.  To calculate the
degree shift of this map, one follows the usual routine.  We let
$(\Sigma,\alphas,\betas,\gammas,\deltas,z)$ be the Heegaard quadruple
representing the integer surgery (so that $Y_{\alpha,\beta}=Y$,
$Y_{\alpha,\gamma}=Y_0$, $Y_{\alpha,\delta}=Y_p$) fix intersection
points $\x\in\Ta\cap\Tb$, $\y\in\Ta\cap\Tc$ (representing $\spinct$),
and $\w\in\Ta\cap\Td$ (representing $Q(\spinct)$), and let
$\Theta_{\beta,\gamma}$, $\Theta_{\beta,\delta}$, and
$\Theta_{\gamma,\delta}$ be intersection points with the canonical
grading (the top non-zero grading in $\HFleq$ of the corresponding
three-manifolds). Indeed, we select the intersection points $\x$,
$\y$, and $\w$ so that there are triangles
$\varphi_1\in\pi_2(\x,\Theta_{\beta,\gamma},\y)$,
$\varphi_2\in\pi_2(\y,\Theta_{\gamma,\delta},\w)$ with
$\Mas(\varphi_1)=\Mas(\varphi_2)=0$. (As usual, in the case where $g=1$, we need to stabilize once to 
achieve this.)

Now we can find alternative triangles
$\psi_1\in\pi_2(\x,\Theta_{\beta,\delta},\y)$ and
$\psi_2\in\pi_2(\Theta_{\beta,\gamma},\Theta_{\beta,\delta},\Theta_{\gamma,\delta})$
with $\Mas(\psi_2)=0$, so that the square obtained by juxtaposing
$\psi_1$ and $\psi_2$ is homotopic to the square obtained by
juxtaposing $\varphi_1$ and $\varphi_2$. Since the the $\SpinC$
structure determined by $\varphi_1+\varphi_2$ factors through the
torsion $\SpinC$ structure on $Y_0$, we can use the previous lemma to
conclude that if $\Sigma$ is a generator for $H_2$ of the cobordism
from $Y$ to $Y_p$ (corresponding to the triple
$(\Sigma,\alphas,\betas,\deltas,z)$), then $$\langle
c_1(\spinc_z(\psi_1)),[\Sigma]\rangle = \pm p.$$ This, together with the
additivity of the Maslov index, which forces $\Mas(\psi_1)=0$, gives that
$$\liftGr(\w)-\liftGr(\w)=\frac{p-5}{4}.$$
Thus, the degree shift of $F_2$ is calculated by
$$\liftGr(\w)-\liftGr(\y)=\liftGr(\w)-\liftGr(\x)+\liftGr(\x)-\liftGr(\y)
=\frac{p-5}{4}+\OneHalf,
$$ (bearing in mind that the degree shift of $F_1$ is $\OneHalf$).
\end{proof}

It will be useful also to have the following result for three-manifolds with 
trivial $\HFpRed$. 

\begin{lemma}
\label{lemma:InvisQSphere}
Let $K\subset Y$ be a knot in an integer homology three-sphere with
$\HFp(Y)=0$, let $p$ be a positive integer, and suppose that
$\HFpRed(Y_p)=0$. Then, for all integers $n\geq p$, $\HFpRed(Y_n)=0$.
\end{lemma}

\begin{proof}
Recall that for a rational homology sphere $L$, it is always the case
that $$|H^2(L;\Z)|\leq \Rk\HFa(L).$$ Moreover, the condition that
$\HFpRed(L)=0$ is equivalent to the condition that $\HFa(L)$ is a free $\Z$-module with
$$|H^2(L;\Z)|
=\Rk
\HFa(L).$$ Thus, to establish
the lemma, it suffices to show that if $\HFa(Y)$ and $\HFa(Y_n)$ are
free $\Z$-modules of rank $1$ and $n$ respectively, then
$\HFa(Y_{n+1})$ is a free $\Z$-module of rank $n+1$.  This, in turn,
follows readily from the general form of the surgery long exact
sequence (Theorem~\ref{HolDiskTwo:thm:GeneralSurgery}
of~\cite{HolDiskTwo}), which specializes to give exactness in: $$
\begin{CD}
...@>>>\HFa(Y)@>>>\HFa(Y_n)@>>>\HFa(Y_{n+1})@>>>\HFa(Y)@>>>...
\end{CD}
$$
It follows immediately that $\HFa(Y_{n+1})$ is a free $\Z$-module 
whose rank satisfies:
$$|H_2(Y_{n+1};\Z)|=n+1\leq \Rk\HFa(Y_{n+1})\leq \Rk\HFa(Y)+\Rk\HFa(Y_n)=1+n.$$
\end{proof}

\vskip.2cm
\noindent{\bf{Proof of Theorem~\ref{thm:PSurgeryLens}.}}
The case of $Q(0)$ is analogous to the case of $+1$ surgeries
considered earlier (in the proof of
Theorem~\ref{thm:RenormalizedComplexity}). Since in this case, the map
$\HFinf(Y_0,[Q(0)])$ to $\HFinf(Y_p,Q(\spinct_0))$ is
surjective, so we still have exactness in the middle for $$
\begin{CD}
\HFpRed(Y) @>>> \HFpRed(Y_0,[Q(0)]) @>>> \HFpRed(Y_p,Q(0)),
\end{CD}
$$
where, under the  identification of $\SpinC(Y_0)\cong \Z$ given earlier,
we have that 
$$\HFpRed(Y_0,[Q(0)])\cong 
\bigoplus_{k\in\Z}\HFpRed(Y_0,kp).$$
Since we assume that
$Y=S^3$ and $Y_p=L(p,q)$, it follows that $\HFpRed(Y_0,[Q(0)])=0$. 
This forces $\HFp(Y_0,kp)=0$ for all $k\neq 0$. Also,
$\HFp(Y_0,0)$ is determined by $d_{\pm 1/2}(Y_0)$:
$$\HFp_k(Y_0,0)=
\left\{
\begin{array}{ll}
\Z & {\text{if $k\equiv -\OneHalf\pmod{2\Z}$ and $k\geq d_{-1/2}(Y_0)$}} \\
\Z & {\text{if $k\equiv \OneHalf\pmod{2\Z}$ and $k\geq d_{1/2}(Y_0)$}} \\
0 & {\text{otherwise}}
\end{array}\right.$$
Moreover, the exact sequence guarantees that $F_1$ maps the image of
$\HFinf(S^3)$ injectively into $\HFinf(Y_0,0)$, while $F_3$
maps the elements coming from $\HFinf(Y_p,Q(0))$ (which in our
present case is all of $\HFp(L(p,q),Q(\spinct_0))$ trivially into
$\HFp(S^3)$.  Thus, it follows that actually $F_1$ maps $S^3$
injectively into $\HFp(Y_0,\spinct_0)$. Since $F_2$ lowers degree by
$1/2$ (Lemma~\ref{lemma:CalcDegrees}), it follows that
$d_{-1/2}(Y_0,\spinct_0)=-1/2$. Again, by exactness, $F_2$ maps the
elements of grading $\equiv 1/2\pmod{2}$ injectively to $L(p,q)$. By
Lemma~\ref{lemma:CalcDegrees}, it then follows that $d_{1/2}(Y_0,\spinct_0)$ is calculated by
the formula claimed.

Having determined $\HFp(Y_0,i)$ for all $i\equiv 0 \pmod{p}$, we
turn to the case of integers $i\in \Z$ with 
the property that $\spinct_0\not\in Q(i)$. In this
case, $\HFp(Y_0,[Q(i)])$ is a finitely $\Z$-module generated, 
so clearly $U^d\HFp(Y_0,[Q(i)])=0$
for sufficiently large powers of $d$. It follows (since $U^d\colon
\HFp(S^3)\longrightarrow \HFp(S^3)$ is surjective for all $d$) that
the image of $\HFp(S^3)$ under $F_1$ is trivial, so that
$\HFp(Y_0,[Q(i)])$ is a $\Z[U]$-submodule of
$\HFp(L(p,q),\spinct)\cong \Z[U^{-1}]$, which is a finitely generated
$\Z$-module. Thus, it follows that there is some integer $\ell$ with
the property that 
\begin{equation}
\label{eq:CyclicModule}
\bigoplus_{j\equiv
i\pmod{p}}\HFp(Y_0,j)\cong\Z[U]/U^\ell.
\end{equation} Indeed, since the module on
the right is a cyclic $\Z[U]$-module, it follows immediately that for
each integer $i$ with $|i|<p$, there is at most one $j\equiv
i\pmod{p}$ with $\HFp(Y_0,j)\neq 0$. To complete the proof, it remains
to show that the integer in this equivalence class is the one with
minimal absolute value, and then to see that $\ell$ is determined as
in Equation~\eqref{eq:EqForL}.

To show the minimality of $j$, we proceed as follows. Observe first
that the above arguments apply to a more general setting: we have
shown that if $K\subset S^3$ is a knot with the property that for some
positive integer $n$, $S^3_n(K)\cong L$, where $L$ is a three-manifold
with $\HFpRed(L)=0$, then for each $i\in\Zmod{n}$, there is at most
one integer $j$ with $j\equiv i\pmod{n}$ with the property that
$\HFp(Y_0,j)\neq 0$. Now, returning to the lens case, note that if
$\HFp(Y_0,j)\neq 0$, then it is also the case that $\HFp(Y_0,-j)\neq
0$ (by the conjugation invariance of the invariants) -- so we can
assume without loss of generality that $j>0$. But both $\HFp(Y_0,j)$
and $\HFp(Y_0,-j)$ lie in the same $2j$-orbit $2j\cm H^2(Y_0;\Z)$; so
it follows immediately that $\HFpRed(Y_{2j})\neq 0$. But then
Lemma~\ref{lemma:InvisQSphere} forces $2j<p$.

Now, it remains to express the integer
$\ell$ from Equation~\eqref{eq:CyclicModule}
in terms of correction terms. To this end,
observe that
the map $F_3$ is realized as a sum of maps belonging to the 
canonical cobordism $W_p$ from
$Y_p$ to $Y$:
$$F_3 = \sum_{\{\spinc\in\SpinC(W_p)
\big|\spinc|Y_p=Q(i)\}}
\pm \Fp{W_p,\spinc},$$ 
where here $W_p=-W(Y,K,p)$ in the notation of
Lemma~\ref{lemma:IdentifyQ}.  
Moreover, each of these component maps
$\Fp{W_p,\spinc}$ is
induced from the corresponding map on
$\HFinf$, $\Finf{W_p,\spinc}$.  Thus, these various maps
differ only by a dimension shift. Moreover, they must all be
isomorphisms, since $\HFp(Y_0,[Q(i)])$ is finitely generated.

To calculate the dimension shift, observe that, 
the cobordism $W_{p}$ from $Y_p$ to $Y$
has $b_2=1$, containing a surface $\Sigma$ with $\Sigma\cm
\Sigma=-p$. Now, if $\spinc$ is any $\SpinC$ structure over
$W_p$, then we have that $$\langle c_1(\spinc),S \rangle =-p+2j,$$ and
indeed the integer $j$ uniquely characterizes $\spinc$. Fixing the
restriction of the $\SpinC$ structure to $Y_p$ fixes the congruence
class of $j \pmod{p}$.  By the dimension formula, the map
$\Finf{W_p,\spinc}$ shifts degree by $$\frac{p-(2j-p)^2}{4p}.$$ When
$j\not\equiv 0\pmod{\Z}$ (and this is equivalent to the assumption that
$\spinct\neq Q(\spinct_0)$), there is a unique maximal such dimension
shift for all $\SpinC$ structures with given restriction to $Y_p$,
which is found by letting $i$ be the representative of its congruence
class with $0\leq i<p$. Indeed, according to
Lemma~\ref{lemma:IdentifyQ}, writing $\spinc|Y_p=Q(\spinct)$, we have
that $$2i\equiv \langle c_1(\spinct),H\rangle \pmod{2p},$$ so the
maximal dimension shift is given by $-d(L(p,1),i)$ (compare
Equation~\eqref{eq:COneLens}). 

It is now easy to see (using the fact that the $\Finf{W_p,\spinc}$ are
all isomorphisms, and since the $\Fp{W_p,\spinc}$ are all induced from
the maps on $\HFinf$) that that the dimension $\ell$ of the kernel of
$F_3$ is given by $$2\ell = d(Y)-d(Y_p,Q(\spinct))+d(L(p,1),i),$$
where $\langle c_1(\spinct), H\rangle \equiv 2i \pmod{2p}$.  In the
present case, since $Y=S^3$, $d(Y)=0$, we see that $\ell$ satisfies
Equation~\eqref{eq:EqForL}.
\vskip.2cm

\section{Calculations} 
\label{sec:SampleCalculations}

We have seen several general results obtained by combining the
absolute gradings with the surgery long exact sequences. In the
present section, we calculate the Floer homologies for a number of
three-manifolds using these techniques. More calculations will be given
in~\cite{HolDiskApp}.

\subsection{Surgeries on torus knots, revisited}
\label{subsec:TorusKnots}

We begin with the trefoil.

The Alexander polynomial for the trefoil is $T-1+T^{-1}$,
so $\MT_0=1$, and all other $\MT_i=0$. Let $Y_0$ denote the manifold
obtained by $0$-surgery on the right-handed trefoil. 

Recall that $+5$ surgery on the right-handed trefoil gives rise to the lens space
$L(5,4)$. Thus, it follows from the long exact sequence (as applied in
Theorem~\ref{thm:PSurgeryLens}) that
\begin{equation}
\label{eq:ZeroSurgeryRHT}
\HFp_k(Y_0,\spinc_0)\cong
\left\{ \begin{array}{ll}
\Z & {\text{if $k\equiv 1/2\pmod{2}$ and $k\geq -3/2$}} \\
\Z & {\text{if $k\equiv -1/2\pmod{2}$ and $k\geq -1/2$}} \\
0 & {\text{otherwise}}
\end{array}\right.
\end{equation}
and $\HFp(Y_0,\spinc)=0$ if $\spinc\neq \spinc_0$. Letting $\gamma\in
H_1(Y_0,\Z)$ be a generator, the action by $\gamma$ is an isomorphism
$\HFp_{k}(Y_0,\spinc_0)\longrightarrow\HFp_{k-1}(Y_0,\spinc_0)$
if $k\equiv 1/2\pmod{2}$ and $k\geq 1/2$, and the action
is trivial otherwise.

Recall that if $p$, $q$, and $r$ are a triple of relatively prime integers, then the Briskorn
variety $V(p,q,r)$ is the locus
$$V(p,q,r)=\{(x,y,z)\in\C^3 \big| x^p+y^q+z^r=0, |x|^2+|y|^2+|z|^2=1\}.$$
The Brieskorn sphere $\Sigma(p,q,r)$ is the homology sphere obtained by $V(p,q,r)\cap S^5$ 
(where $S^5\subset \C^3$ is a standard three-sphere). 
This three-manifold inherits a natural orientation, as the boundary of $V(p,q,r)\cap B^6$ (which in turn
is a manifold away from the origin). 
Now, with these orientation conventions, the three-manifold obtained as $+1$ surgery on the right-handed trefoil
knot is $-\Sigma(2,3,5)$. Another application of the long exact
sequence, this time with $+1$ surgery, shows that 
$$\HFp_k(-\Sigma(2,3,5))
=  \left\{\begin{array}{ll}
\Z & {\text{if $k$ is even and $k\geq -2$}} \\
0 & {\text{otherwise}} 
\end{array}.\right.$$
Moreover, $$U\colon \HFp_{k}(-\Sigma(2,3,5))\longrightarrow
\HFp_{k-2}(-\Sigma(2,3,5))$$ is trivial when $k=0$, otherwise
it is an isomorphism: i.e. $d(-\Sigma(2,3,5))=-2$ and
$\HFred(-\Sigma(2,3,5))=0$.

Now, the Brieskorn sphere $\Sigma(2,3,7)$ is obtained as $-1$ surgery on the 
right-handed trefoil. The exact sequence for $-1$ surgery now reads:
$$\begin{CD}
...@>>>\HFp(\Sigma(2,3,7))@>>> \HFp(Y_0) @>>> \HFp(S^3) @>>>...
\end{CD}$$
Since $\HFp_k(S^3)=0$ for all $k<0$, the generator of
$\HFp_{-3/2}(Y_0,\spinc)$ must come from a generator of
$\HFp_{-1}(\Sigma(2,3,7))$. Thus, we get that
\begin{equation}
\label{eq:Sigma237}
\HFp_k(\Sigma(2,3,7))=\left\{\begin{array}{ll}
\Z & {\text{if $k$ is even and $k\geq 0$}} \\
\Z & {\text{if $k=-1$}} \\
0 & {\text{otherwise}} \end{array}\right.
\end{equation}
Moreover, $\HFred(\Sigma(2,3,7))$ has rank one and $d(\Sigma(2,3,7))=0$.

Indeed, in a similar vein, 
we can consider the manifold $Z_{-n}$
obtained by $-n$-surgery on the right-handed trefoil for any 
negative integer $-n$.
Applying the long exact sequence for surgeries with
negative integer coefficients, we get again that
$\HFpRed(Z_{-n},\spinc)=0$ for each $\spinc\neq Q(\spinc_0)$, while
$$\HFp_k(\Sigma(Z_{-n},Q(\spinc_0))=\left\{\begin{array}{ll}
\Z & {\text{if $k$ is even and $k\geq 0$}} \\
\Z & {\text{if $k=-1$}} \\
0 & {\text{otherwise}} \end{array}\right.,$$ with
$d(Z_{-n},Q(\spinc_0))=0$. This gives an alternate calculation of
Proposition~\ref{HolDiskOne:prop:Trefoil} from~\cite{HolDisk}, in view
of the fact that $Z_{-n}=-Y_n$, and that in~\cite{HolDisk}, the
surgery was performed on the left-handed trefoil. Indeed, even when
$n$ is even, we get an explicit characterization of the $\SpinC$
structure $Q(\spinc_0)$, justifying
Remark~\ref{HolDiskOne:rmk:IdentifySpinStructure} of~\cite{HolDisk}.

To calculate fractional surgeries on the right-handed trefoil, we must first
understand $\uHFp(Y_0;\Zmod{n})$, for a surjective representation
$\Z\cong H^1(Y_0;\Z)\longrightarrow \Zmod{n}$. For this, we apply the
long exact sequence for positive integer surgeries with twisted
coefficients, c.f. Theorem~\ref{HolDiskTwo:thm:ExactPTwist}
of~\cite{HolDiskTwo}. (Actually, there we considered the universal
twisting, Larent polynomials in $T$, whereas here we specialize to
$\Zmod{n}$ twisting, but the proof there
works for any  specialization.) The long exact sequence
in this context reads as follows:
$$
...\longrightarrow 
\uHFp(Y_0,\spinc_0)
\longrightarrow
\HFp(L(5,4),Q(\spinc_0))[\Zmod{n}] 
\longrightarrow
\uHFp(S^3)[\Zmod{n}] 
\longrightarrow ...$$
In the above notation, if $A$ is a $\Z$-module $A[\Zmod{n}]$ denotes
the induced $\Z[\Zmod{n}]$ module $A\otimes_\Z \Z[\Zmod{n}]$. Of
course, as a $\Z$ module, this is simply a direct sum of 
$n$ copies of $A$. 
Recall also that $\uHFinf_k(Y_0,\spinc_0)\cong\Z$ for
all $k\equiv \OneHalf\pmod{1}$.
It follows then that
$$\HFp_k(Y_0,\spinc_0)\cong \left\{ \begin{array}{ll}
\Z & {\text{if $k\equiv 1/2\pmod{2}$ and $k\geq 1/2$}} \\
\Z & {\text{if $k\equiv -1/2\pmod{2}$ and $k\geq -1/2$}} \\
\Z^n\cong \Z[\Zmod{n}] & {\text{if $k=-3/2$}} \\
0 & {\text{otherwise}}
\end{array}\right.$$

Now applying the exact sequence for $1/n$ surgeries 
for positive integers $n$, 
we get that 
$$\HFp_k(Z_{1/n})\cong \left\{
\begin{array}{ll}
\Z & {\text{if $k$ is even and $k\geq 0$}} \\
\Z^n & {\text{if $k=-2$}} \\
0 & {\text{otherwise}}
\end{array}\right.$$
and $d(Z_{1/n})=-2$. Observe that $Z_{1/n}$ is the Brieskorn sphere
$-\Sigma(2,3,6n-1)$. 
Similarly, using the sequence for $-1/n$ surgeries, we get that
$$\HFp_k(Z_{-1/n})\cong \left\{
\begin{array}{ll}
\Z & {\text{if $k$ is even and $k\geq 0$}} \\
\Z^n & {\text{if $k=-1$}} \\
0 & {\text{otherwise}}
\end{array}\right.,$$
with $d(Z_{-1/n})=0$.
Note also that $Z_{-1/n}\cong\Sigma(2,3,6n+1)$.

The key point which facilitated the above calculation was that some positive integral surgery on the
trefoil gives rise to a lens space. More generally, we have the following:

\begin{prop}
\label{prop:CorrTermTorusKnots}
Let $K\subset S^3=Y$ be a knot with the property that some $+p$ surgery on $S^3$ gives a 
lens space, and let $Y_0$ denote the three-manifold obtained by $0$-surgery along
$K$.
Then, for each $i\neq 0$, 
$$\HFp(Y_{0}(K),i)\cong \Z[U]/U^{\MT_i}$$
as a $\Z[U]$ module, which is annihilated by the action of $H_1(Y_0;\Z)$,
where $\MT_i=\MT_i(K)$.
For $i=0$, $\HFp(Y_{0},0)$ is a quotient of $\HFinf(Y_{0},0)$, 
and
\begin{eqnarray*}
d_{-1/2}(Y_0)=-\frac{1}{2} &{\text{and}}&
d_{1/2}(Y_0)=\frac{1}{2}-2\MT_0.
\end{eqnarray*}
Moreover, we have that
\begin{eqnarray*}
d(Y_{1/n})=-2\MT_0
&{\text{and}}&
N(Y_{1/n})=(n-1) \cm \MT_0, + 2n\sum_{i=1}^\infty \MT_i  
\end{eqnarray*}
while
\begin{eqnarray*}
d(Y_{-1/n})=0 
&{\text{and}}&
N(Y_{-1/n})=n\cm \MT_0 + 2n\sum_{i=1}^\infty \MT_i.
\end{eqnarray*}
\end{prop}

\begin{proof}
The statement about $Y_0$ follows from the integral surgeries long
exact sequence, as applied in Theorem~\ref{thm:PSurgeryLens}. The
statements about $Y_{\pm 1/n}$ then follow easily from the fractional
surgeries long exact sequences, as above.
\end{proof}

Note that the above proposition applies to arbitrary torus knots: fix relatively prime positive integers
$p$ and $q$, and let $K_{p,q}$ denote the right-handed $(p,q)$ torus knot. 
It follows from Kirby calculus that some positive surgery of $S^3$ along $K_{p,q}$  always gives a lens space.
Recall also that the Alexander polynomial of $K_{p,q}$ is given by
\begin{equation}
\label{eq:AlexanderTorus}
\Delta_{K_{p,q}}(T)=T^{\frac{p+q-pq-1}{2}}\frac{(1-T)(1-T^{pq})}{(1-T^p)(1-T^q)}
\end{equation}
The $\MT_i(K_{p,q})$ can be calculated from this in the obvious way.

Indeed, there are other knots satisfying the hypothesis of
Proposition~\ref{prop:CorrTermTorusKnots}, including, for
example, the $(-2,3,7)$ pretzel knot (see~\cite{FSLensSurgeries} for
this and more examples).

\subsection{Surgeries on the figure-eight knot}
\label{subsec:Borromeans}

Using the Borromean rings as a stepping-stone
(compare~\cite{FintSternK3}), we can calculate $\HFp$ for fractional
surgeries on the figure eight knot.

Following notation from~\cite{FintSternK3}, 
let  $M\{p,q,r\}$ denote the three-manifold
obtained from $S^3$ by surgeries on the Borromean rings with
coefficients $p$, $q$, and $r$. It is an exercise in Kirby calculus to
see that $M\{p,1,-1\}$ is the manifold obtained by $p$-surgery on the
figure eight knot in $S^3$, while $M\{p,1,1\}$ is $p$-surgery on the
right-handed trefoil. In particular, $M\{-1,1,1\}\cong \Sigma(2,3,7)$.

\begin{prop}
\label{prop:FigureEightZero}
The manifold $M\{-1,0,1\}$, which is zero-surgery on the figure eight
knot has $$\HFp_k(M\{-1,0,1\})\cong \left\{\begin{array}{ll}
\Z & {\text{if $k\equiv \OneHalf\pmod{\Z}$ and $k\geq \OneHalf$}} \\
\Z\oplus \Z &{\text{if $k=-\OneHalf$}} \\
0 & {\text{otherwise}}
\end{array}
\right..$$
Moreover, $d_{-1/2}=-1/2$ and $d_{1/2}=1/2$.
\end{prop}

\begin{proof}
Use the long exact sequence for the three-manifolds
$$S^3\cong M\{-1,\infty,1\},~~~~~~~~~~~~ M\{-1,0,1\}, ~~~~~~~~~~~~
M\{-1,1,1\}\cong \Sigma(2,3,7), $$ 
and Equation~\eqref{eq:Sigma237}.
\end{proof}

\begin{prop}
\label{prop:FigureEight}
Let $E_n$ denote the three-manifold obtained by $1/n$-surgery on
the figure eight knot in $S^3$ (with integral $n>0$). Then, 
$$\HFp_k(E_n)\cong \left\{\begin{array}{ll}
\Z &   {\text{if $k\equiv 0\pmod{2}$ and $k\geq 0$}} \\
\Z^{n} & {\text{if $k=-1$}} \\
0 & {\text{otherwise}}
\end{array}
\right..$$
Moreover, $d(E_n)=0$.
\end{prop}

\begin{proof}
Using the surgery exact sequence for the three-manifolds 
$$M\{-1,0,1\},~~~~~~~~~~~~ M\{-1,\infty,1\}, ~~~~~~~~~~~~
M\{-1,1,1\} $$ with twisted coefficients (in $\Z[\Zmod{n}]$), we see
that 
$$\uHFp_k(M\{-1,0,1\};\Zmod{n})\cong \left\{\begin{array}{ll}
\Z & {\text{if $k\equiv \OneHalf\pmod{\Z}$ and $k\geq \OneHalf$}} \\
\Z\oplus \Z[\Zmod{n}] &{\text{if $k=-\OneHalf$}} \\
0 & {\text{otherwise}}
\end{array}
\right..$$
The result then follows from the $1/n$ surgery exact sequence.
\end{proof}

\subsection{Connected sums}

We remark that more examples can also be constructed using the
connected sum theorem for $\HFm$ (see
Theorem~\ref{HolDiskTwo:thm:ConnSumHFm} of~\cite{HolDiskTwo}). 

For example, it follows easily from Equation~\eqref{eq:Sigma237} (together
with the usual long exact sequence relating $\HFm$ and $\HFp$, Equation~\eqref{eq:HFinfExactSequence}) that
$\HFm(\Sigma(2,3,7))$ is generated as a $\Z[U]$-module by a generator
$\alpha\in\HFm_{-2}(\Sigma(2,3,7))$ (with $U\cm \alpha=0$), and a free
summand generated by an element $\theta\in\HFm_{-2}(\Sigma(2,3,7))$. It follows from the connected
sum theorem that for $$Y=\Sigma(2,3,7)\#\Sigma(2,3,7),$$ $\HFm(Y)$ is
generated as a $\Z[U]$-algebra by elements 
\begin{eqnarray*}
\left(\theta\otimes\theta\right), \left(\alpha\otimes\alpha\right), 
\left(\alpha\otimes\theta\right), \left(\theta\otimes\alpha\right)\in \HFm_{-2}(Y),
&{\text{and}}&\left(\alpha * \alpha\right)\in\HFm_{-3}(Y);
\end{eqnarray*}
and, with the exception of $\left(\theta\otimes\theta\right)$,
all of the other generators are annihilated by $U$. Dualizing again, 
we see that $\HFp(Y)$ 
is generated by three elements in $\HFp_{-1}(Y)$, and one in 
$\HFp_{-2}(Y)$, in addition to the chain of generators
coming from $\HFinf(Y)$ (bearing in mind that $d(Y)=0$).

\subsection{The three-torus}
\label{subsec:T3}

\begin{prop}
\label{prop:T3}
Let $T^3$ denote the three-dimensional torus. Then, we have $H_1(T^3;\Z)$-module isomorphisms:
\begin{eqnarray*}
\HFa(T^3)&\cong &H^2(T^3;\Z)\oplus H^1(T^3;\Z), \\
\HFp(T^3)&\cong & \Big(H^2(T^3;\Z)\oplus H^1(T^3;\Z)\Big)\otimes_\Z \Z[U^{-1}], \\
\HFinf(T^3)&\cong & \Big(H^2(T^3;\Z)\oplus H^1(T^3;\Z)\Big)\otimes_\Z \Z[U,U^{-1}].
\end{eqnarray*}
The absolute grading is symmetric, in the sense that
$\liftGr(H^2(T^3;\Z)\subset \HFa(T^3))=1/2$, 
$\liftGr(H^1(T^3;\Z)\subset \HFa(T^3))=-1/2$.
\end{prop}

Observe that $\HFinf(T^3)$ is smaller than $\HFinf(\#^3(S^1\times
S^2))$. By analogy with Seiberg-Witten theory, this corresponds to the
singular reducible in the character variety, giving rise to a center
manifold picture, compare~\cite{MMR} and~\cite{MMSz}.

\begin{proof}
This, too, is proved by considering surgeries on the Borromean rings,
continuing notation from the previous section. We find it most
convenient to calculate $\HFa$, first.  Since $\HFa(M\{1,1,1\})\cong
\Z$ which is supported in dimension $-2$, and
$\HFa(M\{1,1,\infty\})=\HFa(S^3)\cong \Z$ (supported in dimension
zero), it follows from the surgery exact sequence that 
$\HFa(M\{0,1,1\})\cong \Z\oplus \Z$, where the generators have degree
$-1/2$ and
$-3/2$. Since $b_1(M\{0,0,1\})=2$, it follows from Theorem~\ref{HolDiskTwo:thm:HFinfGen}
of \cite{HolDiskTwo}, that for each $k$, $\HFinf_k(M\{0,0,1\})\cong\Z\oplus\Z$.
Thus, another application of the surgery exact sequence, for the
triple $M\{0,\infty,1\}\cong S^1\times S^2$, and $M\{0,0,1\}$ and
$M\{0,1,1\}$, gives us that $\HFa\{0,0,1\}\cong \Z^2 \oplus \Z^2$,
with two generators in dimension $0$ and two in dimension $-1$. 

Our final surgery exact sequence is applied to the triple
$M\{0,0,\infty\}\cong \#^2(S^1\times S^2)$, $M\{0,0,0\}\cong T^3$, and
$M\{0,0,1\}$. 
Recall that $\HFa(\#^2(S^1\times S^2))=\Z\oplus \Z^2
\oplus \Z$ in gradings $1$, $0$, and $-1$ respectively. 
In this case, the surgery long exact sequence does not uniquely
determine the groups $\HFa(T^3)$, so we proceed as follows.

For simplicity, we will work over a field $\Field$.
Observe first that for each degree $i$, $\HFa_i(T^3,\Field)\cong
\HFa_{-i}(T^3,\Field)$, since $T^3$ admits an orientation-reversing
diffeomorphism. The long exact sequence then implies that
$$\HFa(T^3)=A\oplus A$$ (in degrees $1/2$ and $-1/2$), where $A$ is a
$\Field$-vector space of dimension $\leq 3$. Now, from
the long exact sequence associated to the triple $(\HFa,\HFp,\HFp)$,
and since $\HFa$ is supported in only two consecutive dimensions, it
follows easily that $\HFp(T^3,\Field)\cong (A\oplus A)\otimes_\Field\Field[U^{-1}]$,
and also that $\HFinf(T^3,\Field)\cong (A\oplus A)\otimes_\Field \Field[U,U^{-1}]$.

We argue that the dimension of $A$ can be no smaller than $3$, with the
help of the calculation of $\uHFinf$ in the completely twisted case
(Theorem~\ref{HolDiskTwo:thm:HFinfTwist} of~\cite{HolDiskTwo}). This
gives $$\uHFinf(T^3,\Field[H^1(T^3;\Z)])\cong \Field[U,U^{-1}]$$ as a module over the
ring of Laurent polynomials $\Field[H^1(T^3;\Z)]$. Now, we have an identification
$$\CFinf(T^3,\Field)\cong \uCFinf(T^3,\Field[H^1(T^3;\Z)])\otimes_{\Field[H^1(T^3;\Z)]}\Field,$$
giving rise to a universal coefficients spectral sequence
$$\Tor^i_{\Field[H^1(T^3;\Z)]}(\uHFinf_j(T^3;\Field[H^1(T^3;\Z)]))\Rightarrow \HF_{i+j}(T^3,\Field).$$
Clearly, $\Tor^i_{\Field[H^1(T^3;\Z)]}(\Field)\cong H_*(T^3;\Field)$. Thus, the $E_2$ term in the spectral sequence has the form:
$$\begin{array}{cccc}
\vdots	&	\vdots	&	\vdots	&	\vdots \\
\Field	&	\Field^3	&	\Field^3	&	\Field \\
0	&	0	&	0	&	0 \\
\Field	&	\Field^3	&	\Field^3	&	\Field \\
0	&	0	&	0	&	0 \\
\vdots	&	\vdots	&	\vdots	&	\vdots
\end{array}$$
In particular, the only possible non-trivial differential is $d_3$ from the leftmost to the 
rightmost columns, giving a lower bound of $3$ on the dimension of $A$. 

Together with our previous upper bound, it follows that the dimension of $A$ is three.
The proposition with $\Z$ coefficients then follows easily from the above statement, using
fields $\Q$ and $\Zmod{p}$ for all primes $p$.
\end{proof}

We can calculate $\uHFp(T^3)$ for completely twisted coefficients by
modifying the above techniques.  To state the answer, observe that
there is a canonical map $$\epsilon\colon
\Field[H^1(T^3;\Z)]\longrightarrow \Z,$$ which sends all of $H^1(T^3;\Z)$
to $1$.

\begin{prop}
\label{prop:T3Twist}
There is an identification of $\Z[H^1(T^3;\Z)]$-modules:
$$\uHFp_k(T^3,\spinc_0) \cong \left\{
\begin{array}{ll}
0 		& 	{\text{if $k\equiv 3/2\pmod{2}$ and $k\geq 3/2$}} \\
\Z 		& 	{\text{if $k\equiv 1/2\pmod{2}$ and $k\geq 1/2$}} \\
\ker \epsilon 	& 	{\text{if $k=-1/2$}} \\
0 		& 	{\text{otherwise.}}
\end{array}\right.
$$
\end{prop}

Let $\Laurent(t_1,...,t_b)$ denote Laurent polynomials in $b$
variables, so that if $Y$ has first Betti number $b$, then
$\Z[H^1(Y;\Z)]\cong\Laurent(t_1,...,t_b)$. In this notation, then,
$\ker\epsilon$ for $T^3$ consists of Laurent polynomials
$f(t_1,t_2,t_3)$, with $f(1,1,1)=0$.

\begin{lemma}
Identifying $\Z[H^1(M\{0,1,1\})]\cong\Laurent(t)$, we have an
identification of $\Laurent(t)$-modules: $$\uHFp_k(M\{0,1,1\})\cong
\left\{\begin{array}{ll}
\Z & {\text{if $k\equiv -1/2\pmod{2}$ and $k\geq -\OneHalf$}} \\
\Laurent(t) & {\text{if $k=-3/2$}} \\
0 & {\text{otherwise}}
\end{array}
\right..$$
\end{lemma}

\begin{proof}
Consider the twisted surgery sequence for $\HFa$ connecting
$M\{\infty,1,1\}\cong S^3$, $M\{0,1,1\}$, and $M\{1,1,1\}\cong
-\Sigma(2,3,5)$. It follows that $\uHFa(M\{0,1,1\})\cong
\Laurent(T)\oplus\Laurent(T)$, generated in dimensions $-1/2$ and
$-3/2$.  

In general, it follows from the long exact sequence relating $\uHFa$
and $\uHFp$ that if $\uHFa_k(Y)=0$ for all $k\geq m$, then the map 
$U\colon \uHFp_{i+1}(Y)\longrightarrow \uHFp_{i-1}(Y)$ is an
isomorphism for all $i\geq m$; i.e. $\uHFp_{j}(Y)\cong\uHFinf_{j}(Y)$
for all $j\geq m-1$. Applying this principle to the above calculation
of $\uHFa(M\{0,1,1\})$, and the general calculation of $\uHFinf(Y)$
from Theorem~\ref{HolDiskTwo:thm:HFinfTwist} of~\cite{HolDiskTwo}, the
lemma follows.
\end{proof}

\begin{lemma}
\label{lemma:001twist}
$$\uHFp_k(M\{0,0,1\})\cong
\left\{\begin{array}{ll}
\Z & {\text{if $k\equiv 1 \pmod{2}$ and $k\geq 1$}} \\
0 & {\text{if $k\equiv 0\pmod{2}$ or $k\leq -2$}} \\
\Z\oplus\Laurent(t_1,t_2) & {\text{if $k=-1$}}.
\end{array}
\right.$$
Moreover, the reduced homology group consists of only $\Laurent(t_1,t_2)$ in dimension $-1$.
\end{lemma}

\begin{proof}
Using the long exact sequence for $\uHFa$ connecting
$M\{0,1,\infty\}\cong S^1\times S^2$, $M\{0,1,0\}$, and $M\{0,1,1\}$,
it follows that $\uHFa_k(M\{0,1,0\})=0$ for all $k\neq -1,0$. It then
follows that $\uHFp_k(M\{0,1,1\})$ has the claimed form in all
dimensions except possibly $k=-1$.
Let $f_1$ denote the map 
$$f_1\colon \uHFp(M\{0,1,1\})[t_2,t_2^{-1}]\longrightarrow \uHFp(M\{0,1,\infty\})[t_2,t_2^{-1}].$$
Around $k=-1$, we have the exact sequence reads:
\begin{equation}
\label{eq:SurgTwM001}
\begin{CD}
0\longrightarrow (\Image f_1)\cap\Laurent(t_2)
\longrightarrow\Laurent(t_2) \longrightarrow\uHFp_{-1}(M\{0,0,1\})\longrightarrow\Laurent(t_1,t_2) \longrightarrow0
\end{CD},
\end{equation}
where we use the identification
$\uHFp_{-1/2}(M\{0,1,\infty\})[t_2,t_2^{-1}]\cong \Laurent(t_2)$.
Observe here that the variable $t_1$ corresponds to a generator for
the cohomologies of $S^1\times S^2$ and $M\{0,1,1\}$, while $t_2$
corresponds to a new generator in $H^1(M\{0,0,1\})$. 

We claim that the cokernel of the first map in the above exact sequence is $\Z$. 
To see this, observe that, $U$ induces an isomorphism 
$$\uHFp_{3/2}(M\{0,1,\infty\})[t_2,t_2^{-1}]\cong \uHFp_{~-1/2}(M\{0,1,\infty\})[t_2,t_2^{-1}]$$
which gives an identification  of submodules
$$(\Image f_1)\cap \uHFp_{3/2}(M\{0,1,\infty\})[t_2,t_2^{-1}]\cong (\Image f_1)\cap
\uHFp_{~-1/2}(M\{0,1,\infty\})[t_2,t_2^{-1}],$$
in view of the fact that multiplication by $U$ also induces an
isomorphism $$\uHFp_{\geq 3/2}(M\{0,1,1\})[t_2,t_2^{-1}]\cong
\uHFp_{\geq -1/2}(M\{0,1,1\})[t_2,t_2^{-1}]. $$ Thus, we get an
identification between quotients modules
\begin{eqnarray*}
\frac{\uHFp_{3/2}(M\{0,1,\infty\})[t_2,t_2^{-1}]}{(\Image f_1)\cap
\uHFp_{3/2}(M\{0,1,\infty\})[t_2,t_2^{-1}]}&\cong&
\frac{\uHFp_{-1/2}(M\{0,1,\infty\})[t_2,t_2^{-1}]}{(\Image f_1)\cap
\uHFp_{~-1/2}(M\{0,1,\infty\})[t_2,t_2^{-1}]}\\
&\cong &\uHFp_{1}(M\{0,0,1\}) \\
&\cong& \Z.
\end{eqnarray*}

Thus, the exact sequence in Equation~\eqref{eq:SurgTwM001} shows that 
$$\uHF_{-1}(M\{0,0,1\})\cong \Z\oplus \Laurent(t_1,t_2),$$ as claimed.

\end{proof}

\vskip.2cm
\noindent{\bf{Proof of Proposition~\ref{prop:T3Twist}.}}
From the surgery long exact sequence for $\uHFa$ applied to the triple
$M\{0,0,\infty\}\cong \#^2(S^1\times S^2)$, $M\{0,0,0\}\cong T^3$, and
$M\{0,0,1\}$, it follows that $\uHFa_k(T^3)=0$ for all $k>1/2$, so
$\uHFp_k(T^3)\cong \uHFinf_k(T^3)$ for all $k\geq 1/2$. Again, the
orientation-reversing diffeomorphism of $T^3$ shows that
$\uHFa_{-3/2}(T^3)\cong
\uHFa_{3/2}(T^3)=0$, and thus that $\uHFp_k(T^3)=0$ for all $k\leq -3/2$. 

Thus, it remains to identify $\uHFp_{~-1/2}(T^3)$. 

The long exact sequence gives
\begin{equation}
\label{eq:LExT3Twist}
\begin{CD}
0@>>>\uHFp_{~-1/2}(T^3)@>>> \Laurent(t_3)\oplus \Laurent (t_1,t_2,t_3)
@>{f^+_{-1}}>>\Laurent(t_3)@>>> 0,
\end{CD}
\end{equation}
where the last map is the restriction of
$$f^+\colon \uHFp(M\{0,0,1\})[t_3,t_3^{-1}]\longrightarrow
\uHFp(M\{0,0,\infty\})[t_3,t_3^{-1}]$$
to the part in degree $-1$ (c.f. Lemma~\ref{lemma:001twist}).

We claim that if we further restrict $f^+$ to the summand
$\Laurent(t_3)\subset \Laurent(t_3)\oplus \Laurent (t_1,t_2,t_3)$,
i.e.  $\uHFinf_{-1}(M\{0,0,\infty\})\subset
\uHFp_{~-1}(M\{0,0,\infty\})$, then that restriction is an injection with
cokernel
$\Z$. To see this, observe that the map induced on $\uHFinf$,
$f^\infty$ has some component $g^\infty$ which preserves
$\Z$-degree, and all the other components are translates of
$g^\infty$ by various powers of the $U$-action. 

In view of this, the kernel of $f^+_{-1}$ is identified with the
kernel of the induced surjection $$\Laurent(t_1,t_2,t_3)\longrightarrow \Z\cong
\CoKer \left(f_{-1}|_{\Laurent(t_3)}\right).$$ Any such surjection must carry $t_1$, $t_2$,
and $t_3$ to units in $\Z$, and hence its kernel must be identified, as a 
$\Laurent(t_1,t_2,t_3)$-module, with $\ker \epsilon$.
\qed.

\subsection{The skein exact sequence and some pretzel knots.}

\begin{figure}
\mbox{\vbox{\epsfbox{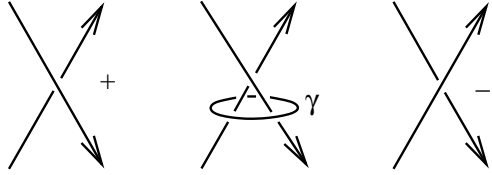}}}
\caption{\label{fig:Skein}
{\bf{Knot crossings for the skein exact sequence.}}  The knot crossing
on the left is a ``positive crossing'' while the one on the right is a
``negative crossing''. In the middle, we have pictured a portion of
$K_+$, together with the curve $\gamma$. (The disk $D$ appearing in
the above discussion is the obvious small disk bounding $\gamma$.)  
Note that the sign of the crossing is independent of the orientation on the knot. }
\end{figure}

Let $Y$ be a three-manifold which is obtained by surgery on a knot
$\KOver\subset S^3$.  Suppose that $D\subset S^3$ is an embedded disk
in $S^3$ which meets $\KOver$ in a pair of intersection points, but
with opposite sign. There is a projection of $\KOver$ for which the
two strands passing through $D$ project to a crossing which is ``positive''
in the usual sense of knot theory (see 
Figure~\ref{fig:Skein}). Now, let
$\KUnder$ be a new knot obtained from $\KOver$ by changing the
over-crossing to an undercrossing.  Let $Y_+=S^3_1(\KOver)$,
$Y_-=S^3_1(\KUnder)$, and let $Y_0$ denote the three-manifold obtained
as a $+1$ surgery on $K_+$, followed by $0$-surgery on the curve
$\gamma=\partial D$. By handlesliding $K$ over $\gamma$, we see that
$Y_0$ could alternatively be thought of as $+1$ surgery on $K_-$
followed by $0$-surgery on $\gamma$. Indeed, if $K'$ is obtained from
$K^+$ by twisting an arbitrary number of times about $D$, the manifold
obtained as $+1$ surgery on $K'$ followed by $0$-surgery on $\gamma$
is diffeomorphic to $Y_0$. (Note that the manifold $Y_0$ has an
alternate description as a sewn-up link complement: let $L_1$ and
$L_2$ be the two components of the link obtained as the self-connected
sum of $K$ using some arc in $D$, then $Y_0$ is an identification
space for $S^3-L_1-L_2$, as observed in ~\cite{Hoste}. However, this
alternate description will not be used in the present discussion.)

There is a long exact sequence of the form
$$
\begin{CD}
.. @>>> \HFp(Y_-)@>>>\HFp(Y_0)@>>>\HFp(Y_+)@>>> ...
\end{CD}
$$ 
Indeed, this is a special case of the general surgery long exact
sequence (Theorem~\ref{HolDiskTwo:thm:GeneralSurgery}
of~\cite{HolDiskTwo}), in view of the fact that $Y_+$ is obtained from
$Y_-$ by a $+1$ surgery on the curve $\gamma$ (as can be seen by
handle-slides over $\gamma$). Long exact sequences of this
kind were introduced by Floer for his instanton homology
(see~\cite{FloerKnots} and \cite{BraamDonaldson}). The analogy with Conway's
``skein relations'' for the Alexander polynomial (and its various
quantum generalizations), should be evident.

Observe also that we have stated a simplified form: the sequence holds
with arbitrary integral surgery coefficient on $\KOver$ (provided that
we perform the same surgery over $K_-$, and the $K_+$ component to
obtain $Y_0$); it also holds if the knot we choose is a single
component of a Kirby calculus link.

Under favorable circumstances, we can use this ``skein exact
sequence'' to calculate $\HFp$ of three-manifolds. We illustrate this for
three-stranded pretzel knots, with odd, positive multiplicities.

\begin{figure}
\mbox{\vbox{\epsfbox{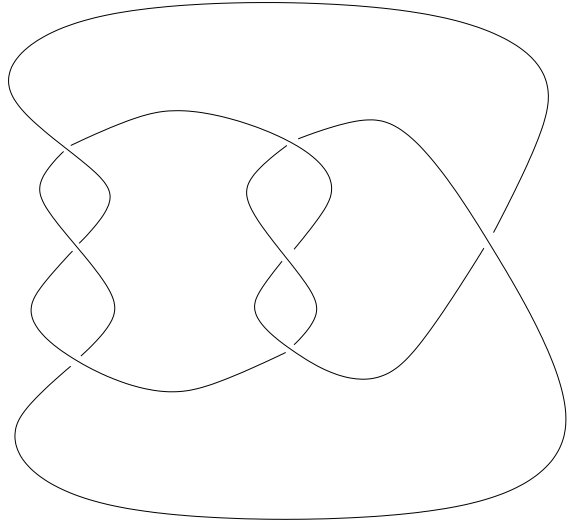}}}
\caption{\label{fig:Pretzel}
{\bf{The (3,3,1) pretzel knot.}}
This is the pretzel knot $(3,3,1)$ with the above notation. }
\end{figure}

Let $a$, $b$ and $c$ be any three integers, and let $P(a,b,c)$
denote the pretzel knot with with three tassels, with $a$, $b$ and $c$
crossings, counted with the sign conventions of
Figure~\ref{fig:Skein}, respectively.
For example, 
$P(1,1,1)$ is the right-handed trefoil, $P(-1,-1,-1)$ is the left-handed trefoil,
$P(-1,1,c)$ is the unknot, and $P(-1,-1,3)$ is the figure eight knot. (See Figure~\ref{fig:Pretzel}.)

\begin{prop}
\label{prop:Pretzel}
Let $\ell$, $m$, and $n$ be three non-negative integers, 
and let $Y(\ell,m,n)$ denote the three-manifold obtained as $0$-surgery on the
pretzel knot $P(2\ell+1,2m+1,2n+1)$.
Then, 
$$\HFp_k(Y)\cong 
\left\{\begin{array}{ll}
\Z^{A+1} & {\text{if $k=-3/2$}} \\
\Z & {\text{if $k\geq -1/2$ and $k\equiv \OneHalf \pmod{\Z}$}} \\
0 & {\text{otherwise,}}
\end{array}
\right.,$$
where $A=m n + \ell n + \ell m + \ell + m + n$. 
In fact, $d_{-1/2}(Y)=-1/2$ and $d_{1/2}(Y)=-3/2$.
\end{prop}

We find it convenient to pass through $\HFa$.
We calculate $\HFa$ for
an auxiliary three-manifold before proceeding to the proof of
Proposition~\ref{prop:Pretzel}.  Let $Y(\ell,m,*)$ denote the
three-manifold which is obtained by $0$ surgery on
$P(2\ell+1,2m+1,1)$, and then $0$-surgery along an unknot $\gamma_3$
which circles the third tassel. In this notation, the skein exact
sequence reads $$
\begin{CD}
...@>>>
\HFa(Y(0,0,-1)) @>>> \HFa(Y(0,0,*)) @>>> \HFa(Y(0,0,0)) @>>>  ... 
\end{CD}
$$ Observe that if we choose $\ell, m$ and $\ell', m'$ so that
$\ell+m=\ell'+m'$, then there is an identification $Y(\ell,m,*)\cong
Y(\ell',m',*)$. To see this, observe that each twist along the second
tassel can be thought of as $-1$-surgery on a standard curve
$\gamma_2$ which goes around that tassel. Handlesliding $\gamma_2$
over $\gamma_3$, we obtain a curve $\gamma_1$, with framing $-1$,
which circles the first
tassel, which we can then remove, hence trading twists along the
second tassel for twists along the first.

\begin{lemma}
\label{lemma:PartPretzel}
Fix non-negative integers $\ell$ and $m$. Then, 
$$\HFa_{k}(Y(\ell,m,*))\cong \left\{
\begin{array}{ll}
\Z^{\ell+m+2} & {\text{if $k=0,-1$}} \\
0 & \text{otherwise}
\end{array}
\right.$$
\end{lemma}

\begin{proof}
We set $\ell=0$, and fix an arbitrary field $\Field$.  We prove that
$\HFa_k(Y(0,m,*))$ (with coefficients in the field $\Field$, which we
suppress from the notation) has rank $\ell+m+2$ in dimensions $k=0$ and $-1$, and has rank $0$ in all other degrees.

To this end, we  claim that there is a skein exact sequence, which reads:
$$
...\longrightarrow
\HFa(Y(0,m,*)) \stackrel{F}{\longrightarrow}
\HFa(T^3) \stackrel{G}{\longrightarrow} 
\HFa(Y(0,m+1,*)) \longrightarrow
 ...,
$$ 
where the middle term $\HFa(T^3)$ was calculated
in Proposition~\ref{prop:T3}. (Though we remind the reader that we are using coefficients in $\Field$.)

To see this, observe that for any $m$, we can use handleslides (over
$\gamma_2$) to pass from the link
$P(1,2m+1,1)\cup\gamma_2\cup\gamma_3$ (all with framing zero) to the
link $P(1,-1,1)\cup\gamma_2\cup\gamma_3$, which is the Borromean
rings. This identifies the middle term with $\HFp(T^3)$.

Moreover, it also follows that $Y(0,0,*)$ is the three-manifold
$M\{0,0,1\}$ in the notation of Section~\ref{subsec:Borromeans}.
Thus, when $m=0$, the lemma was established during the proof of
Proposition~\ref{prop:T3}.

For the inductive step, the inductive hypothesis and the skein exact
sequence sequence clearly give that $\HFa_k(Y(0,m+1,*)$ for all $k\neq
-1, 0$.  Moreover, substituting $H^2(T^3)$ and $H^1(T^3)$ for
$\HFa_{1/2}(T^3)$ and $\HFa_{-1/2}(T^3)$ respectively into the surgery
long exact sequence, we get exactness for: $$
\begin{CD}
0\longrightarrow
H^2(T^3)
\stackrel{G_{\OneHalf}}{\longrightarrow} \HFa_{0}(Y(0,m+1,*))
\longrightarrow \Field^{m+2}
\stackrel{F_0}{\longrightarrow} H^1(T^3)
\stackrel{G_{-\OneHalf}}{\longrightarrow} \HFa_{-1}(Y(0,m+1,*))
\end{CD}
$$

Observe that if $\delta$ is a curve which links $\gamma$ once (i.e. a
copy of the curve along which we perform the surgery to go from $T^3$
to $Y(0,m+1,*)$), then the map $G$ annihilates the image of
$\delta\cm
\HFa(T^3)$ (this follows from the fact that $G$ is induced by cobordisms,
together with naturality of the $H_1$-action); but this image is
two-dimensional inside $H^1(T^3)$.  Thus, the rank of $F_0$
is either two or three.

On the other hand, the following diagram commutes
$$
\begin{CD}
\HFa_{0}(Y(0,m+1,*))@>>>\HFp_0(Y(0,m+1,*)) \\
@V{F_0}VV	@VVV \\
\HFa_{-1/2}(T^3) @>{\cong}>>\HFp_{-1/2}(T^3),
\end{CD}
$$ and since $\HFa_k(Y)=0$ for all $k>0$, it follows that
$$\HFp_0(Y(0,m+1,*))\cong
\HFinf_0(Y(0,m+1,*))\cong \Field^2.$$ This rules out the possibility that $F_0$ surjects.

It follows then that the rank of $\HFa_0(Y(0,m+1,*))$ is $m+3$. This
in turn forces $\HFa_{-1}(Y(0,m+1,*))$ as well (since all other groups
are zero, and the Euler characteristic of $\HFa$ is trivial). Thus, we
have established the lemma for all $m$. Since the field $\Field$ was
arbitrary, the case of integer coefficients follows immediately from
the universal coefficients theorem.

Since $Y(\ell,m,*)\cong Y(0,\ell+m,*)$, the lemma follows for all 
non-negative $\ell$ and $m$.
\end{proof}

\vskip.2cm
\noindent{\bf{Proof of Proposition~\ref{prop:Pretzel}.}}

As usual, we fix an arbitrary field $\Field$. 
We prove for all non-negative $\ell$, $m$, and $n$ that
\begin{equation}
\label{eq:Pretzela}
\HFa_k(Y(\ell,m,n))\cong 
\left\{\begin{array}{ll}
\Field^{A+1} & {\text{if $k=-1/2$ or $-3/2$}} \\
0 & {\text{otherwise,}}
\end{array}
\right.
\end{equation}
where, as before, $A=m n + \ell n + \ell m + \ell + m + n$. 

To establish Equation~\eqref{eq:Pretzela}
for three-manifolds with $\ell=m=0$, we 
proceed as follows.
Observe that $Y(-1,0,n)\cong S^1\times S^2$ (since $P(-1,1,2n+1)$ is
the unknot; indeed, the curve $\gamma_3$ is not linked with
$P(-1,1,2n+1)$), giving
a skein exact sequence
$$
\begin{CD}
...@>>>\HFa(S^1\times S^2)@>{F}>>\HFa(Y(*,0,n))@>{G}>> \HFa(Y(0,0,n)) @>>>...
\end{CD}
$$
In view of Lemma~\ref{lemma:PartPretzel}, we get that
$$\begin{CD}
\HFa_{-1/2}(S^1\times S^2)\cong \Field@>{F_{-1/2}}>>\Field^{n+2} @>{G_{-1}}>> \HF_{-3/2}(Y(0,0,n)) @>>>0
\end{CD}
$$ Letting $\delta$ be a curve linking the surgery curve $\gamma$, we
have once again that $G_{-1}$ annihilates the image of the
$\delta$-action on $\HFa(Y(*,0,n))$. Observe that image of that action
is non-trivial: this follows easily from the fact that the $H_1$
action on $\HFinf(Y(*,0,n))$ is non-trivial. Thus,
$\HF_{-3/2}(Y(0,0,n))\cong \Field^{n+1}$. Next, let $\epsilon$ be a curve
linking the pretzel knot, so that it represents non-trivial homology
classes in both $Y(0,0,n)$ and $Y(*,0,n)$. The map $F$ is equivariant
with respect to action by $\epsilon$, and since the generator of
$\HFa_{-1/2}(S^1\times S^2)$ lies in the image of this action, it
follows that the map $F_{1/2}$ is injective. It then follows immediately
that $\HFa_k(Y(0,0,n))$ has the stated form for all $k$.

Having established Equation~\eqref{eq:Pretzela} for $Y(0,m,n)$ with $m=0$,
we prove the equation with $\ell=0$ and arbitrary $m$, $n$ by induction on $m$.
In this case, the skein exact sequence reads:
$$\begin{CD}
...@>>>\HFa(Y(0,m.n)) @>{F}>> \HFa(Y(0,*,n)) @>{G}>> \HFa(Y(0,m+1,n)) @>>>... 
\end{CD},$$
showing that $\HFa(Y(0,m+1,n))$ is supported in dimensions $k=-3/2$
and $k=-1/2$. Thus $$
\begin{CD}
\HFa_{-1/2}(Y(0,m,n))@>{F_{-1/2}}>> 
\Field^{n+2} @>>> \HFa_{-3/2}(Y(0,m+1,n)) @>>> \Field^{mn+m+n} @>>>0.
\end{CD}
$$

We claim that the image of $F_{-1/2}$ is no more than one dimensional, by
the commutativity of the following diagram:
$$
\begin{CD}
\HFa_{-1/2}(Y(0,m,n)) @>>>\HFp_{-1/2}(Y(0,m,n)) \\
@V{F_{-1/2}}VV @VVV \\
\HFa_{-1}(Y(0,*,n)) @>>>\HFp_{-1}(Y(0,*,n))
\end{CD},
$$
bearing in mind that, since $\HFa_{k}(Y(0,m,n))=0$ for all $k\geq 1/2$,
we have  that $\HFp_{-1/2}(Y(0,m,n))\cong\HFinf_{-1/2}(Y(0,m,n))\cong\Field$;
and also, since $\HF_k(Y(0,*,n))=0$ for all $k<-1$, the natural map
from $\HFa_{-1}(Y(0,*,n))$ to $\HFp_{-1}(Y(0,*,n))$ is an isomorphism.

On the other hand, the image of $F_{-1/2}$ is no less than one
dimensional. Letting $\delta$ be a curve linking $\gamma_2$, it
is easy to see that $\HF_{-1}(Y(0,*n))$ contains an element in the image
$\delta\cm \HF_0(Y(0,*,n))$.

Thus, we have shown that the rank of $\HFp_{-1/2}(Y(0,m+1,n))$ is
$(m+1)n+(m+1)+n$. Since the ranks of $\HFp_{k}(Y(0,m+1,n))$ agree for
$k=-1/2$ and $k=-3/2$, we have established
Equation~\eqref{eq:Pretzela} for three-manifolds of the form
$Y(0,m,n)$ for arbitrary $m$ and $n$.

This same argument is easily modified to give the inductive step (on
$\ell$), establishing Equation~\eqref{eq:Pretzela} for all
non-negative $\ell$, $m$, and $n$.

Going from Equation~\eqref{eq:Pretzela} to the corresponding statement
with coefficients in $\Z$ follows, as usual, from the universal
coefficients theorem, while going from Equation~\eqref{eq:Pretzela} to
the proposition (using $\HFp$) is now a straightforward application of
the long exact sequence relating $\HFa(Y)$ with $\HFp(Y)$.
\qed

\section{Negative-definite intersection forms of smooth four-manifolds}
\label{sec:DefiniteForms}

The aim of the present section is to use the maps associated to
cobordisms to give restrictions on intersection forms of smooth
four-manifolds. 

We begin by giving another proof of the celebrated diagonalizability
theorem of Donaldson~\cite{Donaldson}. The proof parallels the
Seiberg-Witten proof, though the discussion of reducibles is replaced
with the behaviour of the maps on $\HFinf$.  In
Subsection~\ref{subsec:FourBoundary}, we give generalizations for
four-manifolds-with-boundary, which use the correction term $d(Y)$ for
the boundary, in the spirit of Fr{\o}yshov~\cite{Froyshov}. In
Subsection~\ref{subsec:BOneBigForms}, we give further generalizations
for four-manifolds which bound three-manifolds with larger $b_1(Y)$.
As an application of these inequalities, we give another
proof of the ``Thom conjecture'' for $\CP{2}$ in Subsection~\ref{subsec:Thom}.

\subsection{Intersection forms of closed, smooth four-manifolds}

In the present subsection, we give another proof of
of the following result:

\begin{theorem}
\label{thm:Donaldson} (Donaldson)
If $X$ is a smooth, closed, oriented four-manifold with definite
intersection form, then the form is diagonalizable over $\Z$.
\end{theorem}

The theorem follows from two propositions, regarding the map on $\HFinf$ induced by the addition of a two-handle.

We say that $\HFinf(Y)$ is {\em standard} if 
for each torsion $\SpinC$ structure $\spinct_0$,
$$\HFinf(Y,\spinct_0)
\cong \left(\Wedge^{b}H^1(Y;\Z)\right)\otimes_\Z\Z[U,U^{-1}],$$
where $b=b_1(Y)$.

\begin{lemma}
\label{lemma:RkBound}
Let $Y$ be a three-manifold with $b_1(Y)=b>0$, equipped with a torsion
$\SpinC$ structure $\spinct_0$. Then in each dimension $i$, we have that
$$\rk \HFinf_i(Y,\spinct_0)\leq 2^{b-1}.$$
\end{lemma}

\begin{proof}
Consider the universal coefficients spectral sequence whose $E_2$ term
is $$\Tor^i_{\Z[H^1(Y)]}(\Z,\Z)\otimes \uHFinf_j(Y,\spinct_0),$$ and
which converges to $\HFinf_{i+j}(Y,\spinct_0)$.  By the calculation of
$\uHFinf(Y,\spinc_0)$ (Theorem~\ref{HolDiskTwo:thm:HFinfTwist}
of~\cite{HolDiskTwo}), in each degree, the rank of the the $E_2$ term
is $2^{b-1}$. Now, the ranks of the $E_r$ are non-increasing in $r$,
the inequality follows.
\end{proof}

Actually, the rank bound above holds with coefficients in any field ${\mathbb F}$,
in  view of the identification
$\Tor^i_{\Field[H^1(Y;\Z)]}(\Field,\Field)\cong \Wedge^i H_1(Y;\Field)$.

\begin{prop}
\label{prop:ZeroSurgery}
Let $Y$ be closed oriented three-manifold, and $\Knot\subset Y$ a framed
knot, framed so that the cobordism $W(\Knot)$ has
$b_2^-(W(\Knot))=0=b_2^+(W(\Knot))$. Let $\spinc$ be a $\SpinC$ structure on
$W(\Knot)$ whose restrictions to the boundary components 
$Y$ and $Y(\Knot)$, $\spinct$ and $\spinck$, are torsion.
When $\Knot$ represents a non-torsion class in $H_1(Y)$, 
then if $\HFinf(Y,\spinct)$ is standard, then the
induced
map 
$$ F^\infty_{W(\Knot),\spinc}\colon \HFinf(Y,\spinct) 
\longrightarrow \HFinf(Y(\Knot),\spinck)$$ vanishes on the kernel of
the action by $[K]$, inducing an isomorphism on 
$$\HFinf(Y,\spinct)/\Ker[K]
\cong \HFinf(Y(\Knot),\spinck).$$ 
If $\Knot$ represents a torsion class in $H_1(Y)$ and 
$\HFinf(Y(\Knot),\spinck)$ is
standard, then the map $$\HFinf(Y,\spinct) 
\longrightarrow \HFinf(Y(\Knot),\spinck)$$
induces an isomorphism between $\HFinf(Y,\spinct)$ 
and the kernel of the
action by $[L]$ on $\HFinf(Y(\Knot),\spinck)$, where $[L]\in H_1(Y(\Knot))$ is the core of the glued-in
solid torus.
\end{prop}

\begin{proof}
Suppose that $\HFinf(Y)$ is standard, and
$\Knot\subset Y$ is a framed knot whose underlying knot $K$ represents
a non-torsion homology class. Assume first that we are working with $\HFp$ and $\HFinf$ with coefficients
in a field, which we drop from the notation for simplicity.

We find another framing on $K$, denoted $\Knot'$ with the
property that $Y$ and $Y(\Knot)$ fit into a surgery long exact sequence
\begin{equation}
\label{eq:ChangeFraming}
\begin{CD}
...@>>>\HFp(Y(\Knot'))@>{F_1}>>\HFp(Y)@>{F_2}>>\HFp(Y(\Knot))@>>>..., 
\end{CD}
\end{equation}
where we use the general form of the exact sequence, as stated in
Theorem~\ref{HolDiskTwo:thm:GeneralSurgery} of~\cite{HolDiskTwo}. (The
framing $\Knot'$ is one bigger than the framing $\Knot$.)
Here, as usual, $F_1$ and $F_2$ are sums (taken with appropriate signs)
of the maps induced by the two-handle additions.

We have a splitting 
for each torsion $\SpinC$ structure $\spinct$ on $Y$ of
$\HFinf(Y,\spinct)\cong \Ker[K]\oplus
\ker[K]^{\perp}$, where $\Ker[K]=\Image[K]$ is the kernel of the
action by $[K]\in H_1(Y)$ on $\HFinf(Y,\spinct)$ and $\ker[K]^\perp$
is a complementary subspace taken isomorphically to $\Ker[K]$ by
multiplication by $[K]$.

It follows from the fact that the knot $K$ is torsion in $W(\Knot)$
that $\Fp{W(\Knot),\spinc}$ vanishes on the image of the action of
$[K]$ on $\HFinf(Y,\spinc|Y)$, for each $\SpinC$ structure $\spinc$
over $W(\Knot)$; i.e. $\Finf{W(\Knot),\spinc}$ is trivial on the
kernel of $[K]$. In fact, the induced map
$$\Finf{W(\Knot),\spinc}\colon \HFinf(Y,\spinct)/\Ker[K]
\longrightarrow
\HFinf(Y(\Knot),\spinck)$$
is injective, because if it had a kernel element, that would give
$\xi\in\HFinf(Y,\spinct)$ with $K\cm \xi\neq 0$ so that
$\Finf{W(\Knot),\spinc}(\xi)=0$. Moreover, that element $\xi$ would be
in the kernel of $\Finf{W(\Knot),\spinc'}$ for any
$\spinc'\in\SpinC(W(\Knot))$ whose restriction
$\spinc'|Y\cong\spinct$.  (This is clear, by moving the basepoint $z$
in the Heegaard triple representing the cobordism from $Y$ to
$Y(\Knot)$.) Thus, $F_2(\xi)=0$ in the surgery long exact sequence.
By taking preimage of $\xi$ under a sufficiently large power of the
$U$ map, this would give rise to a non-zero element in the image of
$\HFp(Y(\Knot'))$ inside $\HFp(Y)$ which does not lie in the kernel of
the action by $[K]$. But such an element cannot exist, since $[K]$
annihilates the image of $\HFp(Y(\Knot'))$ (as $K$ is torsion in the
cobordism from $Y(\Knot')$ to $Y$).

Since, in each dimension $k$,
the rank of $\HFinf_k(Y,\spinct)/\Ker[K]$ agrees with the upper bound
of the rank of $\HFinf_k(Y(\Knot),\spinct)$ given by Lemma~\ref{lemma:RkBound}
(here we are using the hypothesis that $Y$ has standard $\HFinf$),
it follows that $\Finf{W(\Knot),\spinc}$ is an isomorphism, and that
$\HFinf(Y(\Knot))$ is standard.

To pass from field to $\Z$ coefficients in this case, we observe that
the above proof actually shows that $\Finf{W(\Knot),\spinc}$ always
gives an injection. Applying Lemma~\ref{lemma:RkBound} with
coefficients in $\Zmod{p}$ for each prime, it follows
that $\HFinf(Y(\Knot),\spinck)$ is a free module. Now, the fact that
$\Finf{W(\Knot),\spinc}$ is an isomorphism over $\Z$ follows easily
from the universal coefficients theorem, together with the fact that 
$\Finf{W(\Knot),\spinc}$ is an isomorphism with coefficients in
each $\Zmod{p}$. 

When $\Knot$ represents a torsion class in $Y$, from the previous
arguments (letting $Y(\Knot)$ play the role of $Y$ earlier and
$Y(\Knot')$ play the role of $Y(\Knot)$), we see that $\Ker[L]$ maps
to zero in $Y(\Knot')$. In particular, it follows that there is an
element of $\HFinf(Y,\spinct)$ which maps to an element of $\Ker[L]$,
which generates the submodule $\Ker[L]$ over the ring
$\Z[U,U^{-1}]\otimes_\Z \Wedge^* H_1(Y(\Knot);\Z)$.  From naturality
of the maps of cobordism, it follows that $\Finf{W(\Knot),\spinc}$
surjects onto $\Ker[L]$. Another appeal to Lemma~\ref{lemma:RkBound}
and the hypothesis on the standardness of
$\HFinf(Y(\Knot))$ gives that the map is an
isomorphism (with coefficients in any field, and hence with coefficients
in $\Z$ as above).
\end{proof}

\begin{prop}
\label{prop:NegSurgery}
Let $Y$ be closed oriented three-manifold, and $\Knot\subset Y$ a framed
knot, framed so that the cobordism $W(\Knot)$ has $b_2^-(W(\Knot))=1$. If $\HFinf(Y)$ is standard,
then
for each $\SpinC$ structure $\spinc$ over $W(\Knot)$ whose restriction to
the boundary components $Y$ and $Y(\Knot)$ is torsion, the induced map
$F^\infty_{W(\Knot),\spinc}$ is an isomorphism.
\end{prop}

\begin{proof} 
Let $\Knot'$ denote the same knot, endowed with 
a framing one greater than $\Knot$, i.e. so that 
Exact Sequence~\eqref{eq:ChangeFraming} holds.
Now, there are two cases. Either
$b_1(Y(\Knot'))=b_1(Y)$, or $b_1(Y(\Knot'))=b_1(Y)+1$. In both cases, 
we consider the long exact sequence between $Y(\Knot)$, $Y$, and $Y(\Knot')$
(observing that all three maps are induced by cobordisms).

When $b_1(Y(\Knot'))=b_1(Y)$, we claim that the cobordism from
$Y(\Knot')$ to $Y$ has $b_2^+=1$. Since the associated map shifts the
absolute $\Zmod{2}$ grading by one, and $\uHFinf$ is supported in even
degrees, it follows that the corresponding map on $\HFinf$ must vanish
(c.f. Lemma~\ref{HolDiskFour:lemma:BTwoPlusLemma}
of~\cite{HolDiskFour}). Thus, the image of $\HFp(Y(\Knot'))$ inside
$\HFp(Y)$ is finitely generated. Since $\HFp(Y)$ is infinitely
generated, it follows that the map on $\HFinf$ from $Y$ to $Y(\Knot)$
must be injective. Restricting attention to $\Zmod{p}$ coefficients
where $p$ is any prime, and counting ranks as in
Lemma~\ref{lemma:RkBound}, it follows that the map from $Y$ to
$Y(\Knot)$ is an isomorphism. We can then pass to $\Z$ coefficents as
before.

Assume that $b_1(Y(\Knot'))=b_1(Y)+1$, and again work with $\Zmod{p}$ coefficients.
The knot complement 
gives a natural representation $H^1(Y(\Knot'))\longrightarrow\Z$, 
which in turn gives us a possible twisting of $\HFinf(Y(\Knot'),\Zmod{p}[T,T^{-1}])$.
The twisted long exact sequence reads as follows:
$$
... \longrightarrow \HFp(Y(\Knot))[T,T^{-1}]
\longrightarrow
\uHFp(Y(\Knot'),\Z[T,T^{-1}]) 
\longrightarrow  \HFp(Y)[T,T^{-1}]
\longrightarrow ...,$$
which we can further specialize to $\Zmod{p}$ coefficients.
If the map 
$$\Finf{W(\Knot),\spinc}
\colon\HFinf(Y,\Zmod{p})\longrightarrow \HFinf(Y(\Knot),\Zmod{p})$$
had kernel, it would follow that
$\uHFinf(Y(\Knot'),\Zmod{p}[T,T^{-1}])$ would have (infinitely many)
submodules with non-trivial $T$-action. But this contradicts the fact
that the group $\uHFinf(Y(\Knot'),\Zmod{p})$ with totally twisted
coefficients has a trivial action by the group-ring $\Zmod{p}[H^1(Y;\Z)]$
(this follows easily from Theorem~\ref{HolDiskTwo:thm:HFinfTwist} of~\cite{HolDiskTwo}). It follows that
the map $\Finf{W,\spinc}$
is injective and hence, by Lemma~\ref{lemma:RkBound} (in view
of the fact that $\HFinf(Y,\Zmod{p},\spinct)$ is standard), it is
an isomorphism. Again, since this argument works for all primes $p$,
the statement holds for integral coefficients as well.
\end{proof}

The proof of Theorem~\ref{thm:Donaldson} relies on the following result of Elkies,
see~\cite{Elkies}. Recall that if 
$$Q\colon V\otimes V \longrightarrow \Z$$ is bilinear form over $\Z$, then 
$\xi\in V$ is called a characteristic vector if for all $v\in V$, 
$$Q(\xi,v)\equiv Q(v,v)\pmod{2}.$$ We denote the set of characteristic vectors
for $Q$ by $\Xi(Q)$.

\begin{theorem}(Elkies) Let $Q$ be a negative-definite unimodular 
bilinear form over $\Z$. Then, 
$$0\leq \max_{\xi\in\Xi(Q)} Q(\xi,\xi)+n, $$
with equality if and only if the bilinear form $Q$ is diagonalizable over $\Z$.
\end{theorem}

\vskip.2cm
\noindent{\bf{Proof of Theorem~\ref{thm:Donaldson}.}}
Without loss of generality, we can assume that $b_1(X)=0$ (by
surgering out the one-dimensional homology). We give $X$ a handle
decomposition with a unique zero- and four-handle, and let $W$ be the
associated cobordism from $S^3$ to $S^3$. Decompose $W=W_1\cup W_2
\cup W_3$ into its one-, two-, and three-handles.

We claim that for any
$\SpinC$ structure $\spinc$ over $X$, 
$$\Finf{W,\spinc|W}\colon 
\HFinf(S^3,\spinct_0)\longrightarrow \HFinf(S^3,\spinct_0)$$ is an isomorphism.
To see this, we think of $W_2$ as given by a framed link
$\Link=\Knot_1\cup...\cup \Knot_m$ in $\#^{n_1}(S^2\times S^1)$ (where ${n_1}$
is the number of one-handles in $X$ for our handle-decomposition), and
let $$Y_0=\#^n(S^2\times S^1), Y_1=Y_0(\Knot_1),
Y_2=Y_1(\Knot_2),...,Y_m=\#^{n_3}(S^2\times S^1)$$ 
(where $n_3$ is the number of three-handles in the handle-decomposition).
We claim that since
$b_2^+(X)=0$, the restriction of $\SpinC$ to $Y_i$ is always torsion,
and also $\delta H^1(Y_i;\Z)\subset X$ is trivial (otherwise, we would
have a non-torsion two-dimensional homology class in $X$ coming from $H_2(Y_i)$,
which therefore must have self-intersection number zero). Thus, we can view
$F_{W_2,\spinc}$ as a composite of the maps induced by each individual handle.

Moreover, we claim that we can order the knots so that for
$i=1,...,a$, $b_1(Y_i)$ is decreasing, for $i=a+1,...b$, $b_1(Y_i)$ is
constant (and hence zero), and for $i=a+1,...,m$, and for
$i=b+1,...,m$, $b_1(Y_i)$ is increasing. We call this a {\em standard ordering}.
To achieve a standard ordering, we use two moves.

If $$b_1(Y)<b_1(Y(\Knot_1))>b_1(Y(\Knot_1\cup \Knot_2)),$$ then we can
reorder the knots so that $$b_1(Y)>b_1(Y(\Knot_2))<b_1(Y(\Knot_1\cup
\Knot_2).$$  To see this, observe that the inequality $b_1(Y)<b_1(Y(\Knot_1))$
implies that $K_1$ is a torsion class in $Y$. It is our goal now to rule out the 
possibility that $b_1(Y)\leq b_1(Y(\Knot_2))$. If this inequality were
satisfied, then it would follow that $K_2$ is torsion
in $Y$.  Now, if $K_1$ and $K_2$ were unlinked, then $K_2$ would be
torsion in $Y(\Knot_1)$ as well, contradicting the assumption that
$b_1(Y(\Knot_1))<b_1(Y(\Knot_1\cup\Knot_2))$. If $K_1$ and $K_2$ had
linking number $\ell \neq 0$, then we would be able to find a pair of
surfaces $F_1, F_2\in W(\Knot_1\cup\Knot_2)$ (by capping off the null-homologies of 
$n_1 K_1$ and $n_2 K_2$ in $Y-\Knot_1\cup\Knot_2$) with $F_1\cm F_1=0$,
$F_2\cm F_2=0$, and $F_1\cm F_2=\ell$, which contradicts
$b_2^+(W(\Knot_1\cup\Knot_2))=0$.  

For the second move, observe that if
$$b_1(Y)<b_1(Y(\Knot_1))=b_1(Y(\Knot_1\cup\Knot_2)),$$ then we can
reorder the knots so that
$$b_1(Y)=b_1(Y(\Knot_2))<b_1(Y(\Knot_1\cup\Knot_2)).$$ To see this,
observe that since the first Betti numbers of the three-manifolds do
not drop in the first sequence, it follows that both $K_1$ and $K_2$
represent torsion classes in $H_1(Y)$; moreover it also follows that
$K_1$ and $K_2$ are unlinked. Thus, $K_2$ bounds a Seifert surface
which is disjoint from $K_1$. Since
$b_1(Y(\Knot_1))=b_1(Y(\Knot_1\cup\Knot_2))$, the Seifert framing of
$K_2$ does not agree with the prescribed framing $\Knot_2$. It follows
that $b_1(Y)=b_1(Y(\Knot_2))$, and hence also that
$b_1(Y(\Knot_2))<b_1(Y(\Knot_1\cup\Knot_2))$.

Clearly, by applying the above two moves as necessary, we can arrange for the knots to be
in a standard ordering.

It follows easily from the definition of 
the maps induced by one-handles that 
$$F_{W_1,\spinc_1}(\HFinf(S^3,\spinct_0))\subset
\HFinf(Y_0,\spinc_0)$$ is the subgroup 
$$\Ker[K_1]^\perp\cap...\cap\Ker[K_n]^\perp\subset
\HFinf(Y_0,\spinc_0).$$ Now, since $\HFinf(Y_0,\spinc|Y_0)$ is
standard ($Y_0\cong \#^n(S^2\times S^1)$) successively applying
Proposition~\ref{prop:ZeroSurgery} (in the case where the
$b_1(Y(\Knot))<b_1(Y)$), we see that $\HFinf(S^3)$ is mapped
isomorphically; moreover, by Proposition~\ref{prop:NegSurgery} the
further composite $\Knot_{a+1}\cup...\Knot_{b}$ maps
$\HFinf(S^3,\spinct_0)$ isomorphically to
$\HFinf(Y_b,\spinct|Y_b)$. Successively applying
Proposition~\ref{prop:ZeroSurgery} (this time, in the case where the
knots are homologically trivial), we see that the composite cobordism
carries $\HFinf(S^3,\spinct_0)$ isomorphically to the subgroup of
$\HFinf(Y_m,\spinct|Y_m)$ which is annihilated by the action of
$H_1(Y_m)$.  But it follows easily from the definition of the maps
induced by three-handles that this group is carried isomorphically to
$\HFinf(S^3,\spinct_0)$ under $\Finf{W_3,\spinc|W_3}$.

Since $\Finf{W,\spinc|W}$ is an isomorphism, and
$\HFinf(S^3,\spinct_0)\longrightarrow\HFp(S^3,\spinct_0)$ is surjective, it follows that
$$\Fp{W,\spinc|W}\colon
\HFp(S^3,\spinct_0)\longrightarrow \HFp(S^3,\spinct_0)$$ is surjective;
in particular, we can find some $\xi\in\HFp(S^3)$ so that
$F_{W,\spinc}(\xi)\neq 0$ and $\liftGr(F_{W,\spinc}(\xi))=0$. Thus,
the dimension formula (Equation~\eqref{eq:DimensionShiftFormula})
gives that
$$\liftGr(F_{W,\spinc}(\xi))-\liftGr(\xi)=\frac{c_1(\spinc)^2-2\chi(W)-3\sgn(W)}{4}
= \frac{c_1(\spinc)^2+b_2(X)}{4}\leq 0.$$

This shows that for any characteristic vector for the
intersection form of $H^2(X)$, we have that $$\xi^2 +n \leq 0.$$ 
It follows then from Elkies' theorem cited above that
the intersection form $H^2(X;\Z)$
is diagonalizable.
\qed

\subsection{Intersection forms of definite four-manifolds bounding homology three-spheres}
\label{subsec:FourBoundary}

We give now the generalization of Theorem~\ref{thm:Donaldson} to
four-manifolds bounding rational homology three-spheres.

Let $Y$ be a rational homology three-sphere and $X$ be a four-manifold
which bounds $Y$. The intersection form of $X$ determines a
non-degenerate bilinear form $$Q_X\colon (H_2(X;\Z)/\Tors)\otimes
(H_2(X;\Z)/\Tors)\longrightarrow \Q. $$ More precisely, the image lies
in the subgroup $\frac{1}{|H_1(Y;\Z)|}\Z$.

\begin{theorem}
\label{thm:IntFormQSphere}
Let $Y$ be a rational homology three-sphere, and fix a $\SpinC$
structure $\spinct$ over $Y$. Then, for each smooth, negative-definite
four-manifold $X$ which bounds $Y$, and for each $\SpinC$ structure
$\spinc\in\SpinC(X)$ with  $\spinc|Y\cong \spinct$, we have that
$$c_1(\spinc)^2+\rk (H^2(X;\Z))\leq 4 d(Y,\spinct).$$
\end{theorem}

\begin{proof}
We view $X$ minus a point as a cobordism $W$ from $S^3$ to $Y$, and proceed as
in the proof of Theorem~\ref{thm:Donaldson}, to prove that
for each $\SpinC$ structure $\spinc$ over $X$,
$$\Finf{W,\spinc|W}\colon \HFinf(S^3,\spinct_0)\longrightarrow \HFinf(Y,\spinct)$$
is an isomorphism. Note that now the two-handles give rise to a
cobordism to $Y\#(\#^{n_3}(S^2\times S^1)$ which, again, has standard $\HFinf$. 
The map
$$\Finf{W_3,\spinc}\colon \HFinf(Y\#(\#^{n_3}(S^2\times S^1)),\spinct\#\spinct_0)
\longrightarrow \HFinf(Y,\spinct)$$
induces an isomorphism from $\Ker H_1(\#^{n_3}(S^2\times
S^1);\Z)\subset \HFinf(Y\#(\#^{n_3}(S^2\times
S^1)),\spinct\#\spinct_0)$ onto $\HFinf(Y,\spinct)$, proving the claimed isomorphism.

From this isomorphism, together with the commutative square
$$\begin{CD}
\HFinf_i(S^3,\spinct_0) @>{\Finf{W,\spinc}}>> \HFinf_{d(Y)}(Y,\spinct) \\
@V{\pi}VV @VV{\pi}V \\
\HFp_i(S^3,\spinct_0) @>{\Fp{W,\spinc}}>> \HFp_{d(Y)}(Y,\spinct),
\end{CD}
$$
it follows that 
we can find an element 
$\xi\in\HFp(S^3,\spinct_0)$ with the property that
$\liftGr(\Fp{W,\spinc}(\xi))=d(Y,\spinct)$. 
Thus, we conclude that
$$\liftGr(F_{W,\spinc}(\xi))-\liftGr(\xi)=\frac{c_1(\spinc)^2-2\chi(W)-3\sgn(W)}{4}
= \frac{c_1(\spinc)^2+b_2(X)}{4}\leq d(Y,\spinct).$$
\end{proof}

As a special case, when $Y$ is an integral homology sphere, the
induced bilinear form takes values in $\Z$, and it is unimodular.

\begin{cor}
\label{cor:IntFormZSphere}
Let $Y$ be an integral homology three-sphere, then for each
negative-definite four-manifold $X$ which bounds $Y$, we have the
inequality: $$ Q_X(\xi,\xi)+\rk(H^2(X;\Z)) \leq 4d(Y), $$ for each
characteristic vector $\xi$.
\end{cor}

\begin{proof}
This is an immediate consequence of Theorem~\ref{thm:IntFormQSphere},
with the observation that each characteristic vector for $Q_X$ is the first Chern class of some $\SpinC$ structure over $X$.
\end{proof}

\begin{cor}
\label{cor:NoIntForm}
If $Y$ is an integer homology three-sphere with $d(Y)<0$. Then, there is no
negative-definite four-manifold $X$ with $\partial X=Y$.
\end{cor}

\begin{proof}
This is an immediate consequence of Theorem~\ref{thm:IntFormQSphere} and Elkies' theorem.
\end{proof}

Another consequence of Theorem~\ref{thm:IntFormQSphere} is the
rational homology bordism invariance of $d(Y,\spinct)$. 
Let $(Y_1,\spinct_1)$ and $(Y_2,\spinct_2)$ be a pair of rational homology three-spheres equipped with $\SpinC$ structures. We say that $(Y_1,\spinct_1)$ and $(Y_2,\spinct_2)$ are {\em rational homology
cobordant} if there is a cobordism $W$ from $Y_1$ to $Y_2$ with
$H_*(W;\Q)\cong H_*(S^3\times [0,1];\Q)$, which can be equipped with a
$\SpinC$ structure $\spinc$ with $\spinc|Y_1=\spinct_1$, and
$\spinc|Y_2=\spinct_2$. 

\begin{prop}
\label{prop:HomologyBordismInvariance}
If $(Y_1,\spinct_1)$ and $(Y_2,\spinct_2)$ are rational homology
cobordant (rational homology three-spheres equipped with $\SpinC$
structures), then $d(Y_1,\spinct_1)=d(Y_2,\spinct_2)$.  In particular,
if $(Y,\spinct)$ is a rational homology three-sphere which bounds a
rational homology four-ball $W$, so that $\spinct$ can be extended
over $W$, then $d(Y,\spinct)=0$.
\end{prop}

\begin{proof}
The proof of Theorem~\ref{thm:Donaldson} shows that
if $W$ is a cobordism from $Y_1$ to $Y_2$ with $b_2^+(W)=0$, then
the map 
$$\Finf{W,\spinc}\colon \HFinf(Y_1,\spinct_1)\longrightarrow
\HFinf(Y_2,\spinct_2)$$
is an isomorphism. As before, it follows that
$$d(Y_2,\spinct_1)-d(Y_1,\spinct_2)\geq \frac{c_1(\spinc)^2-2\chi(W)-3\sgn(W)}{4}$$
for any $\SpinC$ structure $\spinc$ over $W$.
When $W$ is a rational homology bordism, then the above formula
gives $d(Y_2,\spinct_2)\geq d(Y_1,\spinct_1)$.  By reversing the orientation of $W$ (and using
Equation~\eqref{eq:DFlipOrientation}), we 
see that $d(Y_1,\spinct_1)\geq d(Y_2,\spinct_2)$, as well.
\end{proof}

Thus, $d$ can be viewed as an obstruction to finding a homology ball
bounding $Y$. Indeed, we have now all the ingredients for
Theorem~\ref{intro:HomologyBordism} stated in the introduction:

\vskip.2cm
\noindent{\bf{Proof of Theorem~\ref{intro:HomologyBordism}}.}
First, we observe that $d$ depends only on the $\SpinC$ cobordism
class of a rational homology sphere $Y$ and $\SpinC$ stucture
$\spinct$; but this was established in
Proposition~\ref{prop:HomologyBordismInvariance} above.  The fact that
$d$ is a homomorphism follows from this, together with the additivity
of $d$ under the connected sum operation, which was established in
Theorem~\ref{thm:AdditivityOfD}.  The fact that $d$ lifts the 
homomorphism $\rho$ (defined in the introduction)
follows immediately from the dimension shift formula for the
absolute grading (Equation~\eqref{eq:DimensionShiftFormula}). Finally, conjugation invariance was established in Proposition~\ref{prop:CorrTermOrient}.
\qed
\vskip.2cm

\subsection{Intersection forms for definite four-manifolds bounding
other three-manifolds}
\label{subsec:BOneBigForms}

Constraints can be given on (semi-definite) intersection forms for
four-manifolds which bound three-manifolds with $b_1(Y)>0$; we will
focus our attention primarily to the case where $H_1(Y;\Z)\cong \Z$.
But first, we set up some terminology.  We have the following easy
consequence of Poincar\'e-Lefschetz duality:

\begin{lemma}
Let $X$ be an oriented four-manifold with boundary $Y$, and let $V$ denote the
image of $H^2(X,Y;\Z)$ in $H^2(X;\Z)/\Tors$. Then, the cup product descends to 
give a non-degenerate bilinear form
$$	Q_X \colon V \otimes V \longrightarrow \Z.		$$
When $H_1(Y;\Z)$ has no torsion, the associated bilinear form is unimodular.
\end{lemma}

\begin{proof}
Non-degeneracy is an immediate consequence of the fact that the
cup-product pairing 
$$\cup\colon  H^2(X;\Z)/\Tors \otimes
H^2(X,Y;\Z)/\Tors\longrightarrow \Z$$ is nondegenerate.

To prove the second claim
(assuming $H_1(Y;\Z)$ has no torsion), we let $K$ denote the kernel of the
natural map $H_1(Y;\Z)\cong H^2(Y;\Z) \longrightarrow H^3(X,Y;\Z)$. Now, 
we get the following (split) short exact sequence
$$\begin{CD}
0@>>> \left(\frac{H^2(X,Y;\Z)}{\delta H^1(Y;\Z)}\right)/\Tors
@>{\iota}>> H^2(X;\Z)/\Tors @>>> K @>>>0,
\end{CD}$$
by considering the Mayer-Vietoris sequence, and using the fact that
$K$ has no torsion. Thus, we can choose bases for
$\left(\frac{H^2(X,Y;\Z)}{\delta H^1(Y;\Z)}\right)/\Tors$, a
Poincar\'e dual basis for its image under $\iota$, which we can then
extend to a basis for $H^2(X;\Z)/\Tors$ by basis vectors which project
to a basis for $K$. With respect to these bases, it is easy to see
that $\iota$ is represented by the intersection matrix $Q_X$,
augmented by a zero matrix. Now, since the cokernel of $\iota$ has no
torsion, $Q_X$ must be unimodular.
\end{proof}

\begin{theorem}
\label{thm:IntFormBOneOne}
Let $X$ be a smooth, oriented four-manifold which bounds a three-manifold $Y$
with $H_1(Y;\Z)\cong \Z$. Let $Q_X$ denote the induced pairing on
$$V=\Image \left( H^2(X,Y;\Z)\longrightarrow H^2(X;\Z)\right)/\Tors,$$ and
suppose that $Q$ is negative definite.  Then, if the restriction map
$H^1(X;\Z)\longrightarrow H^1(Y;\Z)$ is the trivial map, then for each
characteristic vector $\xi$ for $Q_X$, we have that 
\begin{equation}
\label{ineq:IntFormMinusOneHalf}
Q_X(\xi,\xi)+
\rk(V) \leq 4d_{-1/2}(Y)+2;
\end{equation} while if the restriction map
$H^1(X;\Z)\longrightarrow H^1(Y;\Z)$ is non-trivial, then for each 
characteristic vector $\xi$, we have that
\begin{equation}
\label{ineq:IntFormPlusOneHalf}
Q_X(\xi,\xi)+ \rk(V) \leq 4d_{1/2}(Y)-2. 
\end{equation}
\end{theorem}

\begin{remark}
Of course, in the first case, $\rk(V)=\rk H^2(X)-1$, while in the second
case, $\rk(V)=\rk H^2(X)$. Indeed, in the second case, $V\cong
H^2(X,Y)/\Tors\cong H^2(X)/\Tors$. 
\end{remark}

Before turning to the proof, we state the above inequalities in the
case where $X$ has no two-dimensional homology. 

\begin{cor}
\label{cor:NotKnot}
Suppose that $Y$ is a three-manifold with $H_1(Y;\Z)\cong \Z$. Then,
if $Y$ bounds an integral homology $S^2\times D^2$, then  
$d_{-1/2}(Y)\geq -1/2$, while if $Y$ bounds an integral homology $D^3\times S^1$, then
$d_{1/2}(Y)\geq 1/2$.
\end{cor}

\begin{proof}
Apply Inequalities~\eqref{ineq:IntFormMinusOneHalf} and
\eqref{ineq:IntFormPlusOneHalf}, observing that in both applications, the right-hand-side is zero.
\end{proof}

\vskip.2cm
\noindent{\bf{Proof of Theorem~\ref{thm:IntFormBOneOne}.}}
Assume first that the restriction map $H^1(X;\Z)\longrightarrow
H^1(Y;\Z)$ is trivial. We proceed exactly as in the proofs of
Theorem~\ref{thm:Donaldson} and Theorem~\ref{thm:IntFormQSphere}.
First, we surger out all of $b_1(X)$ without affecting its
intersection form.  Then, we puncture $X$ in a point and view the
resulting space as a cobordism $W$ from $S^3$ to $Y$. Order the
two-handles of $W$ as in the proof of Theorem~\ref{thm:Donaldson}, and
observe that since $b_1(Y)<3$, its $\HFinf$ is standard
(c.f. Theorem~\ref{HolDiskTwo:thm:HFinfGen}
of~\cite{HolDiskTwo}). Thus, it follows from this that all the
three-manifolds encountered in the sequence of two-handle additions
have standard $\HFinf$, so Propositions~\ref{prop:ZeroSurgery}
and~\ref{prop:NegSurgery} apply.

In this manner, we show that $$\Finf{W,\spinc|W}\colon
\HFinf(S^3,\spinct_0)\longrightarrow \HFinf(Y,\spinct)$$ is injective,
mapping onto the image of the action by $\gamma\in H_1(Y;\Z)$.
Moreover, we have that $\chi(W)=\chi(X)-1$, $\sgn(W)=-\rk(V)$,
$b_0(X)=1$, $b_1(X)=0$, $b_2(X)=\rk V + 1$, $b_3(X)=\rk
H^1(X,Y;\Z)=0$, and $b_4(X)=0$; thus, the dimension formula implies
that this degree is $$\frac{c_1(\spinc)^2+\rk(V)-2}{4}.$$ In
particular, since the generator of $\HFinf(Y,\spinct_0)$ of degree
$d_{-1/2}$ lies in the image of the $\gamma$ action,
Inequality~\eqref{ineq:IntFormMinusOneHalf} follows.

In the case where the map $H^1(X;\Z)\longrightarrow H^1(Y;\Z)$ is
non-trivial, we proceed as above, except that now we surger out
one-dimensional homology in $X$ until $b_1(X)=1$ and the map in $H^1$
remains non-trivial. In this case, the corresponding map
$\Finf{W,\spinc|W}$ remains injective, only its image is complementary
to the image of the action by $\gamma\in
H_1(Y;\Z)$. Inequality~\eqref{ineq:IntFormPlusOneHalf} follows (bearing in
mind that the first Betti number of the four-manifold is
one). 
\qed
\vskip.2cm

This proof also
gives the
following generalization of Proposition~\ref{prop:CorrTermBound}:

\begin{cor}
\label{cor:FracSurg}
Let $K\subset Y$ be a knot in an integral homology three-sphere.
Then, we have the following inequalities 
(where here $n$ is any positive integer):
$$d_{1/2}(Y_0)-\OneHalf\leq  d(Y_{1/(n+1)})\leq d(Y_{1/{n}}) \leq d(Y) $$
$$d(Y)\leq d(Y_{-1/n})\leq d(Y_{-1/(n+1)}) 
\leq d_{-1/2}(Y_0)+\OneHalf.$$
Furthermore, when $Y\cong S^3$, for all positive integers $n$, we have
that
\begin{eqnarray*}
d_{1/2}(Y_0)-\OneHalf=d(Y_{1/n}) 
\leq 0 \leq 
d(Y_{-1/n})=d_{-1/2}(Y_0)+\OneHalf.
\end{eqnarray*}
\end{cor}

\begin{proof}
Observe first that there are cobordism with $b_2^+=0$ connecting
$Y_0$ to $Y_{1/(n+1)}$ to $Y_{1/n}$ to $Y$ to
$Y_{-1/(n+1)}$ to $Y_{-1/n}$ and back to $Y_0$. 

To construct the cobordism from $Y_0$ to $Y_{1/(n+1)}$, we first take
the knot $K\subset Y$ with $0$-surgery, and then perform surgery along
an additional unknot $L_0\subset Y$ which links $K$ once, given with
framing $-(n+1)$. To go from there on to $Y_{1/n}$, we perform another
surgery along a knot $L_1$ which links $L_0$ once, with framing $-1$.

To go from $Y_{-1}$ to $Y_{-1/n}$, we start with $K$ with framing
$-1$, and we surger along linear plumbing diagram (of length $N$),
with each link given with coefficient $-2$. Surgering along one more
linking circle (with coefficient $-2$) gives the cobordism to
$Y_{-1/(n+1)}$, while framing the linking circle with coefficient $-1$
gives the required cobordism to $Y_0$.

Moreover, in the cobordisms connecting $Y_0$ to $Y_{1/n}$, the image
of $H_1(Y_0)$ is non-trivial, while for the cobordism from $Y_{-1/n}$
to $Y_0$, the one-dimensional homology bounds. The correction terms
are always non-increasing under these cobordisms: for the first and
last cobordisms, we use the proof of Theorem~\ref{thm:IntFormBOneOne},
while for the intermediate steps, we use the version given in
Theorem~\ref{thm:IntFormQSphere} (the intersection forms for these
intermediate forms are obviously diagonalizable, so an appeal to
Elkies' result is unnecessary).

The case where $Y=S^3$ then follows from Proposition~\ref{prop:CorrTermEquality}.
\end{proof}

The most important ingredient in the proof of
Theorem~\ref{thm:IntFormBOneOne}
is that $\HFinf$ of any
three-manifold with $b_1(Y)=1$ is standard. Thus,
Theorem~\ref{thm:IntFormBOneOne} has obvious generalizations to any
three-manifold $Y$ with $b_1(Y)=1$; and there is also a version with
$b_1(Y)=2$, which splits into cases according to the possible images
of $H^1(X;\Z)$ in $H^1(Y;\Z)$. When $b_1(Y)>2$, however, the arguments
run into difficulties, since $\HFinf$ need not be standard (e.g. when
$Y=T^3$), and hence the maps on $\HFinf$ induced by the cobordisms
could be trivial. Indeed, it is also the case that the image of
$H^1(X;\Z)$ in $H^1(Y;\Z)$ cannot be arbitrary: for instance, the cup
product rules out the possibility of a four-manifold $X$ bounding
$T^3$ so that the map $H^1(X,\Z)\longrightarrow H^1(T^3;\Z)$ has
finite cokernel.

However, Theorem~\ref{thm:IntFormBOneOne} generalizes readily to
arbitrary three-manifolds $Y$ with standard $\HFinf$. To keep the
notation simple, we state this only in the case where the restriction
map on $H^1$ is trivial.

\begin{theorem}
\label{thm:IntFormStdHFinf}
Let $Y$ be a three-manifold with standard $\HFinf$, equipped with a
torsion $\SpinC$ structure $\spinct$, and let $d_b(Y,\spinct)$ denote
its ``bottom-most'' correction term, i.e. the one corresponding to the
generator of $\HFinf(Y,\spinct)$ which is in the kernel of the action by
$H_1(Y)$.  Then, for each negative semi-definite four-manifold $W$
which bounds $Y$ so that the restriction map $H^1(W;\Z)\longrightarrow
H^1(Y;\Z)$ is trivial, we have the inequality:
$$c_1(\spinc)^2+b_2^-(W)\leq 4d_b(Y,\spinct)+2b_1(Y)$$
for all $\SpinC$ structures $\spinc$ over $W$ whose restriction to
$Y$ is $\spinct$.
\end{theorem}

\begin{proof}
Follow the proof of Theorem~\ref{thm:IntFormBOneOne}.
\end{proof}

\subsection{The  minimal genus problem in $\CP{2}$}
\label{subsec:Thom}

We give a proof of the Thom conjecture for $\CP{2}$, based on the
theory developed thus far. This result was first proved by
Kronheimer-Mrowka~\cite{KMthom} and
Morgan-Szab{\'o}-Taubes~\cite{MSzT}. The proof we give here is analogous to
a Seiberg-Witten proof given recently by Strle, see~\cite{Strle}.

\begin{theorem}{(Kronheimer-Mrowka, Morgan-Szab{\'o}-Taubes)}
\label{thm:Thom}
Let $\Sigma\subset \CP{2}$ be a smoothly embedded two-manifold, 
which represents $m>0$ times a generator
$H\in H_2(\CP{2};\Z)$. Then, 
\begin{equation}
\label{ineq:AdjCP}
m^2-3m\leq 2g(\Sigma)-2;
\end{equation}
i.e. 
the holomorphic curves in $\CP{2}$ 
minimize genus in their homology class.
\end{theorem}

The proof is based on the results from Section~\ref{sec:DefiniteForms}
on intersection forms (specifically, Theorem~\ref{thm:IntFormStdHFinf}),
together with the following calculation for circle bundles over two-manifolds.

\begin{lemma}
\label{lemma:CorrTermCircleBundle}
Let $Y$ be a circle bundle over a two-manifold, oriented as the boundary of 
a tubular neighborhood $N$ of a two-manifold $\Sigma$ with 
self-intersection number $\Sigma\cdot\Sigma = -n<0$;
and indeed, suppose that $n\geq 2g$, where $g$ denotes the
genus of $\Sigma$.
Let $\spincu$ be the $\SpinC$ structure over $N$ with 
$$\langle c_1(\spincu), [\Sigma] \rangle = -n+2g, $$
and $\spinct$ be its restriction to $Y$. 
Then we have an isomorphism of {\em relatively graded groups}
$$\HFp(Y,\spinct)\cong \HFp(\#^{2g}(S^2\times S^1),\spinct_0),$$
with the bottom-most generator of $\HFp(Y,\spinct)$ in degree 
$$\frac{1}{4} - \frac{g^2}{n} - \frac{n}{4}.$$
\end{lemma}

\begin{proof}
Consider the integral surgeries long exact sequence
(Theorem~\ref{HolDiskTwo:thm:ExactP} of~\cite{HolDiskTwo})
$$
\begin{CD}
...@>>>\HFp(S^1\times \Sigma_g)
@>>> 
\HFp(\#^{2g}(S^1\times S^2)) 
@>>>
\HFp(Y)@>>>...
\end{CD}
$$ 
Recall that the above sequence decomposes, according to $\SpinC$
structures over $Y$,
where we use a sum over all $\SpinC$ structures over $S^1\times \Sigma_g$
in the fiber of $Q\colon \SpinC(S^1\times \Sigma)\longrightarrow \SpinC(Y)$.
If $\spinc$ is a $\SpinC$ structure with $c_1(\spinc)=\ell[\PD(\Sigma_g)]$, then
it is easy to see that $Q(\spinc)$ is the restriction to $Y$
of a $\SpinC$ structure
$\spincu$ over $N$ with $\langle c_1(\spincu),[\Sigma]\rangle = \ell-n$ 
(compare Lemma~\ref{lemma:IdentifyQ}).

Now, by the adjunction inequality for the three-manifold $S^1\times
\Sigma_g$ (c.f. Theorem~\ref{HolDiskTwo:thm:Adjunction}
of~\cite{HolDiskTwo}), $\HFp(S^1\times
\Sigma_g,\spinc)$ is trivial for all $\SpinC$ structures with $Q(\spinc)=\spinct$: it is non-trivial only for those $\SpinC$
structures $\spinc$ for which $c_1(\spinc)=\ell[\PD(\Sigma_g)]$ with
$|\ell|\leq 2g-2$ (and it is trivial for all $\SpinC$ structures whose
first Chern class is not a multiple $\PD[\Sigma_g]$). Thus, the map in the long
exact sequence $$\HFp(\#^{2g}(S^1\times S^2),\spinct_0) \longrightarrow
\HFp(Y,\spinct) $$ is an isomorphism (of relatively graded
groups). In fact, the map can be interpreted as a sum of maps
$$\sum_{\{\spinc\in\SpinC(W)\big|\spinc|_Y\cong \spinct\}} \pm
\Fp{W,\spinc},$$ where $W$ is a single two-handle addition of
$\#^{2g}(S^1\times S^2)$ giving rise to $Y$. The term in this sum
which shifts degree down the least corresponds to the $\SpinC$
structure $\spinc$ with $$\frac{c_1(\spinc)^2+1}{4}=\frac{1}{4}\left(
1-\frac{(2g-n)^2}{n}\right).$$ Since the bottom-most generator of
$\HFp(\#^{2g}(S^1\times S^2),\spinct_0)$ has degree $-g$, the result
follows.
\end{proof}

\vskip.3cm
\noindent{\bf{Proof of Theorem~\ref{thm:Thom}.}}
Suppose that $\Sigma_0\subset \CP{2}$ violates
Inequality~\eqref{ineq:AdjCP}.  By adding handles locally 
if necessary, we can
find another embedded surface $\Sigma\subset \CP{2}$ (representing the
same homology class) with $$m^2-3m=2g(\Sigma).$$ Let $\spinccan$ be
the $\SpinC$ structure over $\CP{2}$ whose first Chern class is
represented by $-3H$ (this is the canonical class of $\CP{2}$).
Then the
restriction of $\spinccan$ to a tubular neighborhood of $\Sigma$
satisfies the hypotheses of Lemma~\ref{lemma:CorrTermCircleBundle}, so
that if $Y$ denotes the boundary of this tubular neighborhood, then
\begin{equation}
\label{eq:BottomCorrTerm}
d_b(Y,\spinccan|Y)=-2+\left(\frac{3m-m^2}{2}\right).
\end{equation}

Let $W$ be the four-manifold with boundary obtained by deleting a
tubular neighborhood of $\Sigma$ from $\CP{2}$. Indeed, according to
the above lemma, $\HFinf$ of $Y$ is standard. Moreover, $H_2(W;\Q)$ is clearly
trivial, as is the map $H^1(W;\Q)\longrightarrow
H^1(Y;\Q)$. Thus, Theorem~\ref{thm:IntFormStdHFinf} applies, and gives the inequality
$$-g=\frac{3m-m^2}{2}\leq d_b(Y,\spinccan|Y),$$
which contradicts Equation~\eqref{eq:BottomCorrTerm}.
\qed

\section{Examples} 
\label{sec:Examples}

In this final section, we illustrate the intersection form results
from Section~\ref{sec:DefiniteForms} by combining them with the
calculations from Section~\ref{sec:SampleCalculations}. In
Section~\ref{sec:IntFormSamples}, we study four-manifolds which bound
surgeries on torus knots, in Section~\ref{subsec:NotSurgeryOnKnot} we
exhibit a homology $S^1\times S^2$ which is not surgery on a single
knot.  Although some of the results contained in the present section
have alternate proofs using more classical gauge-theoretic techniques
(especially the first subsection, which contains some results which
can be found in the work of~\cite{FS},
\cite{Austin},
\cite{Froyshov}), we include the present discussion to  give the reader a
better feel for the theorems in the present paper.
In the final subsection, we illustrate the results of Section~\ref{sec:Lens},
and specifically Corollary~\ref{cor:AlexLens} of that section, by
including a table containing all 
possible symmetric Laurent polynomials in $T$ 
which can arise as Alexander
polynomials of knots in $S^3$ whose $+p$ surgery, for
positive integral $p\leq 26$, gives a lens space.

\subsection{Intersection form bounds}
\label{sec:IntFormSamples}

Continuing notation from Section~\ref{sec:SampleCalculations}, we let
$Y_{p,q}(0)$ denote the three-manifold obtained by zero-surgery 
on the
right-handed $(p,q)$ torus knot.

\begin{prop}
\label{prop:TrefoilBounding}
Let $X_1$ be a four-manifold with $b_2^+(X_1)=0$ and $\partial X_1 =
Y_{2,3}(0)$. Then, the map from $H^1(X_1;\Z)\longrightarrow H^1(Y;\Z)$
is trivial and the intersection form of $X_1$ is
diagonalizable. Similarly, let $X_2$ be a four-manifold with $b_2^+(X_2)=0$
and $\partial X_2=-Y_{2,3}(0)$. Then, if $H^1(X_2;\Z)\longrightarrow
H^1(Y;\Z)$ is non-trivial, then the intersection form of $X_2$ is
diagonalizable. Moreover, if the map on $H^1$ is trivial, then if
$Q_{X_2}$ denotes the intersection form of $X_2$ (on
$V=\Image(H^2(X_2,Y;\Z)\longrightarrow H^2(X_2))$), we have that
$$Q_{X_2}(\xi,\xi)+\rk(V)\leq 8. $$
\end{prop}

\begin{proof}
We have seen in Section~\ref{sec:SampleCalculations} that
\begin{eqnarray*}
d_{-1/2}(Y_{2,3}(0))=-1/2&{\text{and}}&
d_{1/2}(Y_{2,3}(0))=-3/2.
\end{eqnarray*}
Thus, the case where $H^1(X_1;\Z)\longrightarrow H^1(Y_{2,3}(0))$ is
non-trivial is ruled out by Inequality~\eqref{ineq:IntFormPlusOneHalf}
(together with Elkies' theorem), while the diagonalizability of $Q$ in
the other case is forced by Inequality~\eqref{ineq:IntFormMinusOneHalf}
(and another application of Elkies' theorem). 

The correction terms for $-Y_{2,3}(0)$ (which we could alternatively
think of as zero-surgery on the left-handed trefoil knot) are
determined by Equation~\eqref{eq:DFlipOrientationBOne}. The rest of
the proposition then follows from another application of
Theorem~\ref{thm:IntFormBOneOne}. 
\end{proof}

\begin{figure}
\mbox{\vbox{\epsfbox{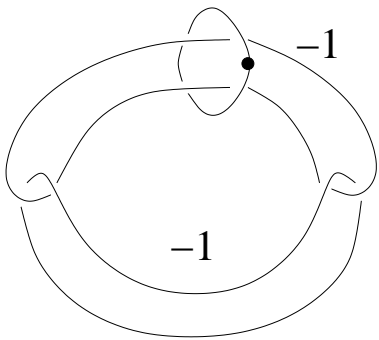}}}
\caption{\label{fig:Tref1}
{\bf{Kirby calculus description for $X_2$.}}  
This is the Kirby calculus description of the four-manifold $X_2$ 
described above, with $\partial X_2 = Y_{2,3}(0)$.}
\end{figure}

\begin{figure}
\mbox{\vbox{\epsfbox{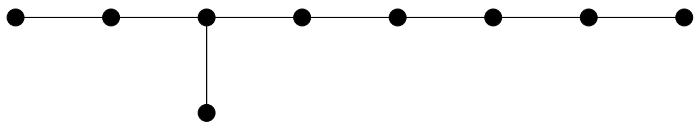}}}
\caption{\label{fig:Tref2}
{\bf{A four-manifold which bounds zero-surgery on the trefoil.}}  The
manifold here ($X_3$) is obtained by plumbing $-2$ spheres: each
vertex represents a sphere of self-intersection $-2$. Each edge
corresponds to an intersection between spheres. This configuration is
simply-connected, and its boundary is $-Y_{2,3}(0)$.}
\end{figure}

Indeed, the inequalities obtained in the above proposition are all
sharp. Let $X_1$ denote four-manifold obtained by attaching a
zero-framed two-handle to the four-ball along a right-handed
trefoil. Let $X_2$ denote the four-manifold whose Kirby calculus
description is given in Figure~\ref{fig:Tref1}:
there is a single one-handle, and a pair of two-handles added with framing $-1$ 
each (along unlinked circles), so that the three circles form the
Borromean rings.
Let $X_3$ denote the four-manifold
obtained by plumbing nine $-2$-spheres as pictured in Figure~\ref{fig:Tref2}.

Now, clearly, $b_2^+(X_1)=0$ and $\partial X_1=Y_{2,3}(0)$ (and its
intersection form is trivially diagonalizable). Moreover, $\partial
X_2=-Y_{2,3}(0)$, and its intersection form is $(-1)\oplus (-1)$
(i.e. it is negative-definite and diagonal), and the map
$H^1(X_2;\Z)\longrightarrow H^1(Y_{2,3}(0);\Z)$ is an
isomorphism. Finally, the intersection form of $X_3$ induced on $V$ is
easily seen to be the negative-definite form $E_8$. According to the inequality
in the above proposition, this is the largest rank of any even 
intersection form which bounds $-Y_{2,3}(0)$.

In fact, Proposition~\ref{prop:TrefoilBounding} admits the following generalization:

\begin{prop} 
In the following statements, $p$ and $q$ can be any pair of
positive, relatively prime integers (both greater than $1$).
\begin{itemize}
\item For all natural
numbers $n$, the manifolds $-\Sigma(p,q,pqn-1)$ cannot bound
a four-manifold $X$ with $b_2^+(X)=0$.
\item The intersection form of any four-manifold $X$ with $b_2^+(X)=0$ which
bounds $\Sigma(p,q,pqn+1)$ is diagonalizable.
\item If $X$ is a four-manifold with $b_2^+(X)=0$ and
$\partial X = Y_{p,q}(0)$, then the map $H^1(X;\Z)\longrightarrow
H^1(Y_{p,q}(0))$ is trivial and the intersection form of $X$ is
diagonalizable.
\end{itemize}
\end{prop}

\begin{proof}
First observe that if $K_{p,q}$ is the $(p,q)$ torus knot, then 
$\MT_0(K)>0$. 
This follows from the fact
that all the non-zero coefficients of the Alexander
polynomial (Equation~\eqref{eq:AlexanderTorus}) are $\pm 1$, and they
come in alternating signs, with the top coefficient $+1$. 

All the above results are direct consequences of this sign, 
the calculations relating
the correction terms with $\MT_0$ 
(Proposition~\ref{prop:CorrTermTorusKnots}), Elkies' result, and
Theorem~\ref{thm:IntFormQSphere} or \ref{thm:IntFormBOneOne} as appropriate.
\end{proof}

\subsection{On manifolds which are not surgery on a knot}
\label{subsec:NotSurgeryOnKnot}

\begin{figure}
\mbox{\vbox{\epsfbox{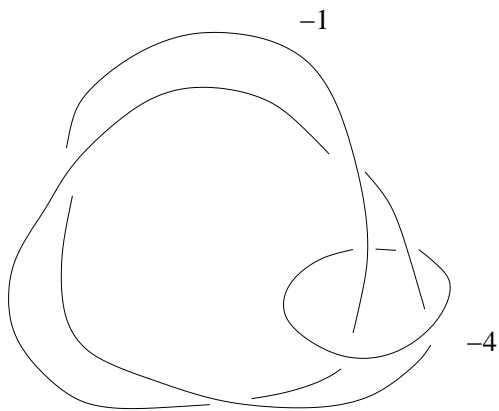}}}
\caption{\label{fig:TwoComponentLink}
{\bf{Kirby calculus description of $Y_{-2}$.}}  
}
\end{figure}

We give a simple illustration of Corollary~\ref{cor:NotKnot}.
Consider the three-manifold $Y_{-2}$ given by the Kirby calculus
description pictured in Figure~\ref{fig:TwoComponentLink}.

\begin{figure}
\mbox{\vbox{\epsfbox{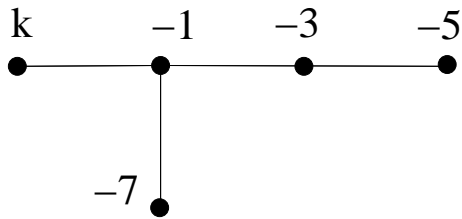}}}
\caption{\label{fig:PlumbSeif}
{\bf{Plumbing description for $Y_{k}$}}}
\end{figure}

This three-manifold can alternately be
given as a plumbing as in Figure~\ref{fig:PlumbSeif}, substituting in
$k=-2$.

From the plumbing description, it is clear that $Y_{-2}$ is obtained
from from $L(49,40)$ (which is the manifold obtained by leaving off
the vertex labelled by $k$) by attaching a single two-handle with
framing $-2$. Since $w_2$ of the plumbing diagram can be represented
by the sum of the Poincar\'e duals of the $(-7)$-sphere and the $(-5)$-sphere,
it follows easily that $W$ is a spin cobordism. Let
$\spinc_0$ denote the $\SpinC$ structure on $W$ with
$c_1(\spinc_0)=0$.  We claim that the map $\Finf{W,\spinc_0}$ is
non-trivial. This follows easily from the fact that the map
$\Fp{W,\spinc_0}$ is a summand in the map $F_2$ belonging to a surgery
long exact sequence: $$
\begin{CD}
...@>>>\HFp(L(49,40)) @>>> \HFp(Y_{-2}) @>>> \HFp(L(49,44))@>>>...,
\end{CD}
$$ where $L(49,44)$ is the manifold which is boundary of the plumbing
pictured in Figure~\ref{fig:PlumbSeif} with $k=-1$. Using
Proposition~\ref{prop:dLens}, we see that the correction term for
$L(49,40)$ in its spin structure is $-2$, and as usual one can see
that the $\Finf{W,\spinc_0}$ decreases degree by $1/2$. It follows
that $d_{-1/2}(Y_{-2})\leq -5/2$. (In fact, it follows from another
glance at the above exact sequence that $\Fp{W,\spinc_0}$ can have no
kernel, and hence that $d_{-1/2}(Y_{-2})=-5/2$.)

Thus, according to Corollary~\ref{cor:NotKnot}, it follows that
$Y_{-2}$ is not obtained as the zero-surgery on a one-component knot
in $S^3$ (and, indeed, that $Y_{-2}$ does not bound any homology
$S^2\times D^2$).

\subsection{Alexander polynomials of knots giving lens spaces.}
\label{subsec:LensAlex}

Results from Section~\ref{sec:Lens} (c.f. Corollary~\ref{cor:AlexLens})
give, for each pair of relatively prime integers $(p,q)$,
an explicitly calculable, 
finite list of symmetric polynomials in $T$ which contains
the Alexander polynomials of all knots $K\subset S^3$ with the property
that $S^3_p(K)\cong L(p,q)$.  It is interesting to compare the
with the conjecture of Berge~\cite{Berge}.

This list is determined as follows. Given $(p,q)$, we list
the correction terms $d(L(p,q),i)$ and $d(L(p,1),i)$ for $i=0,...,p-1$, 
as
determined in Proposition~\ref{prop:dLens}. For each multiplicative generator
generator $h\in(\Zmod{p})^*$, and each $c\in \Zmod{p}$
we define a corresponding sequence $t_i$
(indexed by integers $i$) by the formula $$0\leq 2\MT_i=
\left\{\begin{array}{ll}
-d(L(p,q),c+h\cm i) + d(L(p,1),i) &{\text{if $2|i|\leq p$}} \\
0 & {\text{otherwise.}}
\end{array}
\right. 
$$
If the numbers $t_i$ are all non-negative integers, then there is 
a Laurent polynomial in $T$
$$\Delta(T)= a_0 + \sum_{i=1}^{\infty}a_i(T^i+T^{-i})$$
whose coefficients are given by
$$a_i=\left\{\begin{array}{ll}
t_{i-1}-2t_i+t_{i+1} & {\text{if $i\neq 0$}} \\
1 + t_{-1}-2t_0+t_1  & {\text{if $i=0$}}
\end{array}\right.$$
(compare Equation~\eqref{eq:defBi}).  Let ${\mathcal F}(p,q)$ denote
the set of all symmetric Laruent polynomials obtained in this
manner. In particular, observe that the number of polynomials in
${\mathcal F}(p,q)$ is bounded by $p$ times the number of integers
integers in $\{1,...,p-1\}$ which are relatively prime to $p$ (note
that this is a crude estimate, which can easily be improved, using the
conjugation symmetry).

It follows from Corollary~\ref{cor:AlexLens} that if $K\subset S^3$ is
any knot with the property that $S^3_p(K)\cong L(p,q)$, then its
Alexander polynomial $\Delta_K$ appears in the list ${\mathcal
F}(p,q)$.

We include here a table of ${\mathcal F}(p,q)$ for all integers $p\leq 26$.  To
conserve space, we do not display ${\mathcal F}(p,q)$ when either:
\begin{itemize}
\item $q=1$ (for
then ${\mathcal F}(p,1)=\{1\}$, according to
Theorem~\ref{intro:PLens}),
\item ${\mathcal F}(p,q)=\emptyset$,
\item  we have already displayed ${\mathcal F}(p,q')$, where $q$ and $q'$ are related by
$q\cm q'\equiv 1\pmod{p}$,
because in this case $L(p,q')\cong L(p,q)$, so that
${\mathcal F}(p,q)={\mathcal F}(p,q')$.
\end{itemize}

\newpage
Note that in this list, 
${\mathcal F}(21,16)$ consists of  {\em two} polynomials (while all
other sets displayed consist of a single polynomial).
\begin{tiny}
\begin{eqnarray*}
{\mathcal F}(5,4)&=&
\{ -1 + {T^{-1}} + T\} 
\\
{\mathcal F}(7,4)&=&
\{ -1 + {T^{-1}} + T\} 
\\
{\mathcal F}(9,7)&=&
\{ 1 + T^{-2} - {T^{-1}} - T + T^2\} 
\\
{\mathcal F}(10,9)&=&
\{ 1 + T^{-6} - {2}\,{T^{-5}} + T^{-4} + T^{-3} - T^{-2} - T^2 + T^3 + T^4 - 
   2\,T^5 + T^6\} 
\\
{\mathcal F}(11,9)&=&
\{ 1 + T^{-3} - T^{-2} - T^2 + T^3\} 
\\
{\mathcal F}(11,4)&=&
\{ 1 + T^{-2} - {T^{-1}} - T + T^2\} 
\\
{\mathcal F}(13,12)&=&
\{ -3 + T^{-6} - {2}\,{T^{-5}} + {2}\,{T^{-4}} - T^{-3} + {2}\,{T^{-1}} + 
   2\,T - T^3 + 2\,T^4 - 2\,T^5 + T^6\} 
\\
{\mathcal F}(13,10)&=&
\{ -1 + T^{-3} - T^{-2} + {T^{-1}} + T - T^2 + T^3\} 
\\
{\mathcal F}(13,9)&=&
\{ 1 + T^{-3} - T^{-2} - T^2 + T^3\} 
\\
{\mathcal F}(14,11)&=&
\{ -1 + T^{-4} - T^{-3} + {T^{-1}} + T - T^3 + T^4\} 
\\
{\mathcal F}(15,4)&=&
\{ -1 + T^{-3} - T^{-2} + {T^{-1}} + T - T^2 + T^3\} 
\\
{\mathcal F}(16,9)&=&
\{ -1 + T^{-4} - T^{-3} + {T^{-1}} + T - T^3 + T^4\} 
\\
{\mathcal F}(17,16)&=&
\{ -3 + T^{-9} - {2}\,{T^{-8}} + T^{-7} + T^{-6} - {2}\,{T^{-4}} + T^{-3} + 
   T^{-2} + \\
&& {T^{-1}} 
 + T + T^2 + T^3 - 2\,T^4 + T^6 + T^7 - 2\,T^8 + T^9\} 
\\
{\mathcal F}(17,15)&=&
\{ 1 + T^{-7} - {2}\,{T^{-6}} + {2}\,{T^{-5}} - T^{-4} + T^{-2} - 
   {T^{-1}} - T + T^2 - T^4 + 2\,T^5 - 2\,T^6 + T^7\} 
\\
{\mathcal F}(17,13)&=&
\{ 1 + T^{-4} - T^{-3} + T^{-2} - {T^{-1}} - T + T^2 - T^3 + T^4\} 
\\
{\mathcal F}(18,13)&=&
\{ 1 + T^{-5} - T^{-4} + T^{-2} - {T^{-1}} - T + T^2 - T^4 + T^5\} 
\\
{\mathcal F}(19,17)&=&
\{ 1 + T^{-9} - {2}\, {T^{-8}} + T^{-7} + T^{-6} - T^{-5} + T^{-3} - T^{-2} - 
   T^2 + T^3 - T^5 + T^6 + T^7 - 2\,T^8 + T^9\} 
\\
{\mathcal F}(19,16)&=&
\{ -1 + T^{-6} - T^{-5} + T^{-2} + T^2 - T^5 + T^6\} 
\\
{\mathcal F}(19,11)&=&
\{ 1 + T^{-5} - T^{-4} + T^{-2} - {T^{-1}} - T + T^2 - T^4 + T^5\} 
\\
{\mathcal F}(19,5)&=&
\{ 1 + T^{-4} - T^{-3} + T^{-2} - {T^{-1}} - T + T^2 - T^3 + T^4\} 
\\
{\mathcal F}(20,9)&=&
\{ 1 + T^{-6} - T^{-5} + T^{-3} - T^{-2} - T^2 + T^3 - T^5 + T^6\} 
\\
{\mathcal F}(21,16)&=&
\{ -1 + T^{-5} - T^{-4} + T^{-3} - T^{-2} + {T^{-1}} + T - T^2 + T^3 - 
   T^4 + T^5, \\
&& -1 + T^{-6} - T^{-5} + T^{-2} + T^2 - T^5 + T^6\} 
\\
{\mathcal F}(22,9)&=&
\{ 1 + T^{-6} - T^{-5} + T^{-3} - T^{-2} - T^2 + T^3 - T^5 + T^6\} 
\\
{\mathcal F}(23,18)&=&
\{ -1 + T^{-7} - T^{-6} + T^{-4} - T^{-3} + {T^{-1}} + T - T^3 + T^4 - 
   T^6 + T^7\} 
\\
{\mathcal F}(23,16)&=&
\{ -1 + T^{-7} - T^{-6} + T^{-3} - T^{-2} + {T^{-1}} + T - T^2 + T^3 - 
   T^6 + T^7\} 
\\
{\mathcal F}(23,6)&=&
\{ -1 + T^{-5} - T^{-4} + T^{-3} - T^{-2} + {T^{-1}} + T - T^2 + T^3 - 
   T^4 + T^5\} 
\\
{\mathcal F}(25,24)&=&
\{ -5 + {2}\, {T^{-{12}}} - {3}{T^{-{11}}} + {2}\, {T^{-9}} + T^{-8} - 
   {4}{T^{-7}} + {2}\, {T^{-6}} + {2}\, {T^{-5}} - {2}\, {T^{-4}} - T^{-3} + 
   {2}\, {T^{-2}} + {2}\, {T^{-1}} \\
&& + 2\,T + 2\,T^2 - T^3 - 2\,T^4 + 2\,T^5 + 
   2\,T^6 - 4\,T^7 + T^8 + 2\,T^9 - 3\,T^{11} + 2\,T^{12}\} 
\\
{\mathcal F}(25,21)&=&
\{ 1 + T^{-9} - {2}\, {T^{-8}} + {2}\, {T^{-7}} - T^{-6} + T^{-4} - T^{-3} + 
   T^{-2} - {T^{-1}}  - T + T^2 - T^3 + T^4 - T^6 + 2\,T^7 - 2\,T^8 + T^9\} 
\\
{\mathcal F}(25,19)&=&
\{ 1 + T^{-6} - T^{-5} + T^{-4} - T^{-3} + T^{-2} - {T^{-1}} - T + T^2 - 
   T^3 + T^4 - T^5 + T^6\} 
\\
{\mathcal F}(25,16)&=&
\{ 1 + T^{-8} - T^{-7} + T^{-4} - T^{-3} + T^{-2} - {T^{-1}} - T + T^2 - 
   T^3 + T^4 - T^7 + T^8\} 
\\
{\mathcal F}(25,14)&=&
\{ -1 + T^{-7} - T^{-6} + T^{-4} - T^{-3} + {T^{-1}} + T - T^3 + T^4 - 
   T^6 + T^7\} 
\\
{\mathcal F}(26,25)&=&
\{ 1 + {3}\,{T^{-14}} - {5}\,{T^{-{13}}} + T^{-12} + T^{-11} + T^{-10} + 
   T^{-9} - {3}\,{T^{-8}} + T^{-6} + T^{-5} + T^{-4} - {2}\,{T^{-3}} - 
   T^{-2} + {T^{-1}} \\
&& + T - T^2 - 2\,T^3 + T^4 + T^5 + T^6 - 3\,T^8 + T^9 + 
   T^{10} + T^{11} + T^{12} - 5\,T^{13} + 3\,T^{14}\} 
\\
{\mathcal F}(26,23)&=&
\{ -1 + T^{-12} - {2}\,{T^{-{11}}} + T^{-10} + T^{-9} - T^{-8} + T^{-5} - 
   T^{-4} + T^{-2} \\
&& + T^2 - T^4 + T^5 - T^8 + T^9 + T^{10} - 2\,T^{11} + T^{12}
   \} 
\end{eqnarray*}
\end{tiny}

\commentable{
\bibliographystyle{plain}
\bibliography{biblio}
}

\end{document}